\documentclass[a4paper,11pt,reqno, english]{amsart}

\usepackage{anysize,color}
\usepackage[utf8]{inputenc}
\usepackage[dvipsnames]{xcolor}
\usepackage{float}
\usepackage[colorlinks=true,urlcolor=blue,linkcolor={magenta},citecolor={orange}]{hyperref}
\usepackage{amsfonts,amssymb,amsmath}
\usepackage{amsthm}
\usepackage{url}
\usepackage{paralist}
\usepackage{tikz}
    \usetikzlibrary{decorations.markings,decorations.pathreplacing,arrows,shapes.callouts,calc,positioning,arrows.meta}
\usepackage{rotating}
\usepackage[mathscr]{euscript}
\usepackage{charter}
\usepackage{verbatim}
\usepackage{tikz-cd}
\usepackage{mathtools}
\usepackage{tabularx}
\IfFileExists{soul.sty}{\usepackage{soul}}{}
\usepackage[arrow,curve,matrix,tips,2cell]{xy}
  \SelectTips{eu}{10} \UseTips
  \UseAllTwocells
\IfFileExists{bbold.sty}{\usepackage{bbold}}{}
\usepackage{booktabs}

\usepackage{mathabx}

\newtheorem*{theorem*}{Theorem}
\newtheorem{theorem}{Theorem}[section]
\newtheorem{proposition}[theorem]{Proposition}
\newtheorem{lemma}[theorem]{Lemma}

\newtheorem{corollary}[theorem]{Corollary}
\makeatletter\let\c@equation\c@theorem\makeatother

\theoremstyle{definition}
\newtheorem{definition}[theorem]{Definition}
\theoremstyle{remark}

\newcommand\R{\mathbb R}

\newcommand\Z{\mathbb Z}

\newcommand\FF{\mathbb F}
\newcommand\sym{\mathfrak S}
\newcommand\K{\mathscr{K}}
\newcommand\LL{\mathscr{L}}

\newcommand\Tv{\mathrm{Tv}}
\newcommand\E{\mathbb{E}}

\usepackage{scalerel}
\usepackage{array}
\usepackage{placeins}

\DeclareMathOperator{\Ind}{Ind}
\newcommand\cart{\mathbin{\Box}}
\DeclareMathOperator{\St}{St}
\DeclareMathOperator{\conv}{\mathrm{conv}}

\DeclareMathOperator{\pr}{\mathrm{pr}}
\DeclareMathOperator{\supp}{\mathrm{supp}}

\DeclareMathOperator{\relint}{\mathrm{relint}}

\begin{document}

\title[Tverberg graphs beyond prime powers]{Tverberg graphs beyond prime powers}

\author[Blagojevi\'c]{Matija  Blagojevi\'{c}}
\thanks{The research by Matija  Blagojevi\'{c} leading to these results has
received funding from the Studienstiftung des deutschen Volkes.}
\address{}
\email{matib03@zedat.fu-berlin.de}

\author[Blagojevi\'c]{Pavle V. M. Blagojevi\'{c}}
\thanks{The research by Pavle V. M. Blagojevi\'{c} leading to these results has
received funding from the Serbian Ministry of Science, Technological Development and Innovation. It was also supported by G\"unter Ziegler's group at the Freie Universit\"at Berlin ``Arbeitsgruppe Diskrete Geometrie und Topologische Kombinatorik''.}
\address{Inst. Math., FU Berlin, Arnimallee 2, 14195 Berlin, Germany\hfill\break
	\mbox{\hspace{4mm}}Mat. Institut SANU, Knez Mihailova 36, 11001 Beograd, Serbia}
\email{blagojevic@math.fu-berlin.de}

\begin{abstract}
A Tverberg partition of a map $\Delta_n\longrightarrow\mathbb{R}^d$ is a collection of $r$ pairwise disjoint faces of the simplex $\Delta_n$ whose images have a point in common, and by Tverberg's theorem every affine map has one when $n\geq(d+1)(r-1)$.
We call a graph $H$ on the vertex set of the simplex $\Delta_n$ a Tverberg graph if every map $\Delta_n\longrightarrow\mathbb{R}^d$ has a Tverberg partition into $r$ faces which covers all vertices and separates the endpoints of every edge of $H$.
Beyond matchings, such graphs were known only when $r$ is a prime power.
We show that for every $r\geq3$ and $n\geq(d+1)(r-1)$, with no additional vertices, many graphs are Tverberg graphs for affine maps, for instance forests of maximum degree at most $\frac{b(r)}{2}$, where $b(r)$ is the largest prime power dividing $r$.
More generally, for graphs whose blocks are chordal graphs or cycles, we bound the number of additional vertices depending on properties of the graph, and when $r$ is a prime power these results hold for continuous maps and confirm a conjecture of Engstr\"om for this class.
We also extend a local criterion of Engstr\"om to all $r$: every graph with $|N^2(v)|+2|N(v)|<b(r)$ for all vertices $v$ is a Tverberg graph for affine maps.
\end{abstract}

\maketitle

\bigskip
\section{Introduction and statement of the main results}
\label{sec : intro}
\label{sec : introduction}
\bigskip

The famous 1966 theorem of Tverberg~\cite{Tverberg1966} says that for integers $r\geq2$, $d\geq1$ and $N:=(d+1)(r-1)$ every affine map $a\colon\Delta_N\longrightarrow\R^d$ has a {\em Tverberg partition}, that is, a collection $(F_1,\dots,F_r)$ of pairwise disjoint faces of $\Delta_N$ with
$
	a(F_1)\cap\cdots\cap a(F_r)\neq\emptyset
$.

\medskip
Once a Tverberg partition is known to exist, the next natural question to ask is: How many Tverberg partitions can we expect to exist for an arbitrary affine map $a$?
Permuting the faces of a Tverberg partition gives $r!$ Tverberg partitions with the same faces, so it is reasonable to count unordered Tverberg partitions  $\{F_1,\dots,F_r\}$ instead.
This number of unordered Tverberg partitions of an affine map $a$ is denoted by $\Tv(a)$.
In his notes from 1979, Sierksma~\cite{Sierksma1979} considered the following affine map $a\colon\Delta_N\longrightarrow\R^d$, the so-called $(d+1)$-clouds configuration: split the $N+1=(d+1)(r-1)+1$ vertices of $\Delta_N$ into $d+1$ groups of $r-1$ vertices and one further vertex, map the vertices of the $i$-th group to points close to the $i$-th vertex of a $d$-dimensional simplex in $\R^d$, and map the remaining vertex to a point in the interior of this simplex.
This map has exactly $((r-1)!)^d$ unordered Tverberg partitions, and Sierksma's ``Dutch cheese problem'' asks whether it is extremal, consult also~\cite{Reay1982} and~\cite[Sec.\,4.1]{BaranySoberon2018}.
In other words, {\em Sierksma's conjecture} asserts that every affine map $a\colon\Delta_N\longrightarrow\R^d$ satisfies
\[
	\Tv(a)\ \geq\ ((r-1)!)^d.
\]

\medskip
Let $\Tv_r(d)$ be the minimum of $\Tv(a)$ over all affine maps $a\colon\Delta_N\longrightarrow\R^d$.
The $(d+1)$-clouds configuration gives $\Tv_r(d)\leq((r-1)!)^d$, so Sierksma's conjecture says that $\Tv_r(d)=((r-1)!)^d$.
If the conjecture holds, then the bound is attained by many affine maps, as shown by White~\cite{White2017}, by Bukh, Loh and Nivasch~\cite{BukhLohNivasch2017}, and by P\'or~\cite{Por2018}.
Apart from the elementary cases $r=2$ and $d=1$, the conjecture has been confirmed only for $r=3$ and $d=2$, that is, $\Tv(a)\geq4$ for every affine map $a\colon\Delta_6\longrightarrow\R^2$.
This case is due to Hell~\cite[Thm.\,6]{Hell2008constraints}, see also~\cite[Rem.\,3.19]{Hell2006thesis}, and an elementary proof was recently given by Sober\'on~\cite{Soberon2026}.

\medskip
Lower bounds for $\Tv(a)$ which grow with $d$ were first obtained for prime powers.
If $r=p^\alpha$ is a prime power, then every affine map $a\colon\Delta_N\longrightarrow\R^d$ satisfies
\begin{equation}
	\label{eq : VZ-Hell}
	\Tv(a)\ \geq\ \frac{1}{(r-1)!}\Big(\frac{r}{\alpha+1}\Big)^{\lceil\frac{N}{2}\rceil}.
\end{equation}
This bound is due to Vu\v{c}i\'c and \v{Z}ivaljevi\'c~\cite{VucicZivaljevic1993} for primes and the extension to prime powers is due to Hell~\cite{Hell2007number}.
Both results hold for arbitrary continuous maps $\Delta_N\longrightarrow\R^d$ as well.
For affine maps in general position Hell showed that $\Tv(a)\geq(r-d)!$ for every $r$~\cite{Hell2008Birch}, and that $\Tv(a)\geq\min\{(r-1)!,\,c_{d,r}(r-d)!\}$ with an integer constant $c_{d,r}\geq2$ when $r>2$ is a prime power and $d\geq2$~\cite[Thm.\,5]{Hell2008constraints}.
Until recently, for $r$ not a prime power and $d\geq r-1$ no lower bound better than $\Tv(a)\geq1$ was known.
A recent paper~\cite[Thm.\,1.4]{BlagojevicSierksma} shows that for every $r\geq3$, every decomposition $r=kq+s$ with $q=p^\alpha$ a prime power, $k\geq1$ and $s\geq0$, and every affine map $a\colon\Delta_N\longrightarrow\R^d$
\begin{equation}
	\label{eq : BB}
	\Tv(a)\ \geq\ \frac{q}{r!}\Big(\frac{q}{\alpha+1}\Big)^{\lceil\frac{(r-1-2s)(d+1)}{2}\rceil}.
\end{equation}
For $k=1$ and $s=0$ this is the bound~\eqref{eq : VZ-Hell}.
Taking for $q$ the largest prime power which divides $r$ and $s=0$, the exponent becomes $\lceil\frac{N}{2}\rceil$, and so for every fixed $r\geq3$ the number $\Tv_r(d)$ grows at least exponentially in $d$, see~\cite[Cor.\,1.5]{BlagojevicSierksma}.

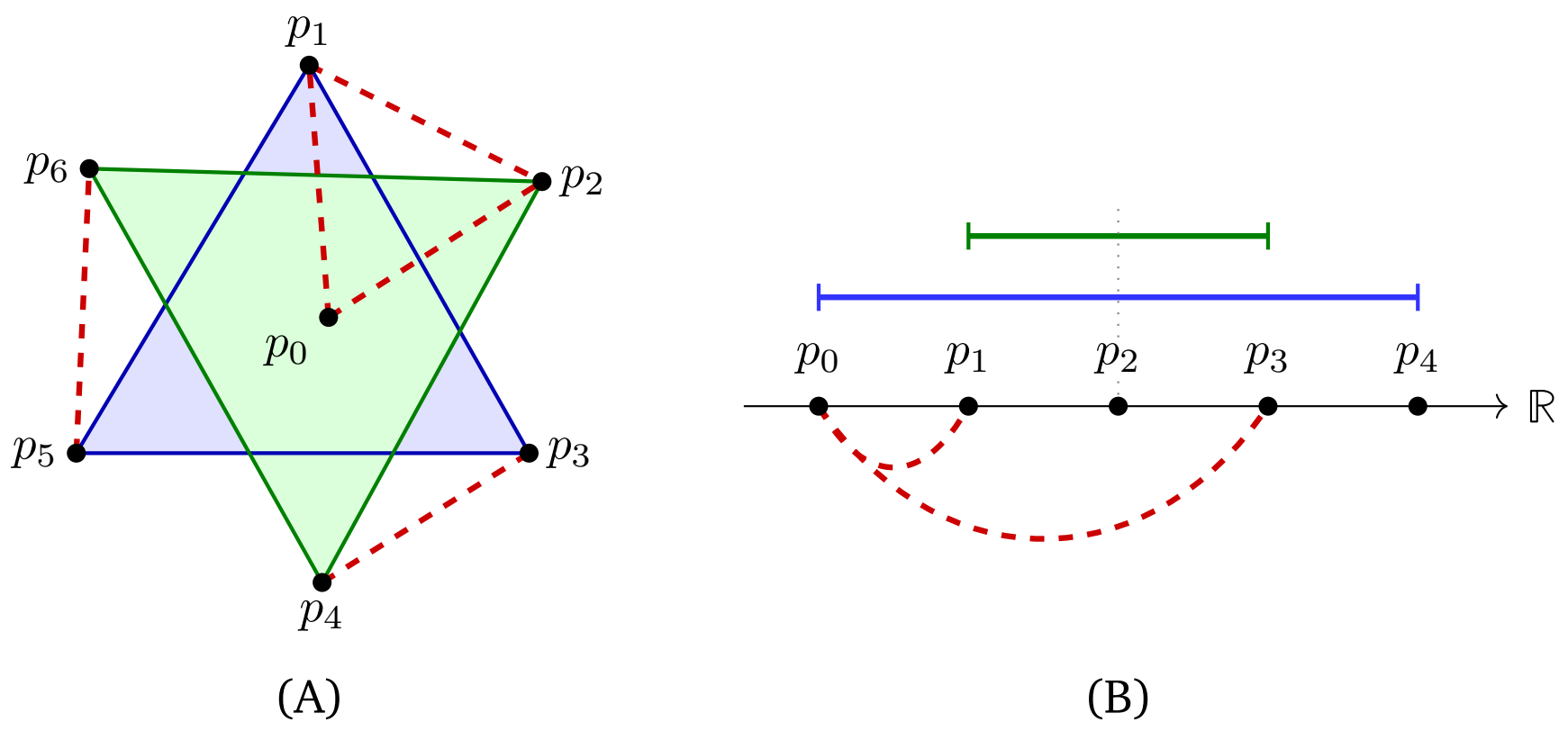
\begin{figure}[ht]
	\centering
	\begin{tikzpicture}[scale=0.95,
			pt/.style={circle,fill=black,inner sep=1.4pt},
			lab/.style={font=\small}]
		\coordinate (p0) at (0.15,0.05);
		\coordinate (p1) at (0,2);
		\coordinate (p2) at (1.8,1.1);
		\coordinate (p3) at (1.7,-1.0);
		\coordinate (p4) at (0.1,-2.0);
		\coordinate (p5) at (-1.8,-1.0);
		\coordinate (p6) at (-1.7,1.2);
		\fill[blue!12] (p1)--(p3)--(p5)--cycle;
		\fill[green!14] (p2)--(p4)--(p6)--cycle;
		\draw[blue!70!black,thick] (p1)--(p3)--(p5)--cycle;
		\draw[green!50!black,thick] (p2)--(p4)--(p6)--cycle;
		\draw[red!80!black,very thick,dashed] (p0)--(p1) (p0)--(p2) (p1)--(p2) (p3)--(p4) (p5)--(p6);
		\foreach \i in {0,...,6}{\node[pt] at (p\i) {};}
		\node[lab,below left] at (p0) {$p_0$};
		\node[lab,above] at (p1) {$p_1$};
		\node[lab,right] at (p2) {$p_2$};
		\node[lab,right] at (p3) {$p_3$};
		\node[lab,below] at (p4) {$p_4$};
		\node[lab,left] at (p5) {$p_5$};
		\node[lab,left] at (p6) {$p_6$};
		\node[font=\small] at (0,-2.9) {(A)};
	\end{tikzpicture}
	\qquad
	\begin{tikzpicture}[x=11mm,y=10mm,
			pt/.style={circle,fill=black,inner sep=1.4pt},
			lab/.style={font=\small}]
		\draw[->] (-0.5,0)--(4.6,0) node[right,font=\small] {$\R$};
		\draw[blue!80,very thick] (0,0.8)--(4,0.8);
		\draw[blue!80,thick] (0,0.7)--(0,0.9) (4,0.7)--(4,0.9);
		\draw[green!50!black,very thick] (1,1.25)--(3,1.25);
		\draw[green!50!black,thick] (1,1.15)--(1,1.35) (3,1.15)--(3,1.35);
		\draw[gray,dotted] (2,0)--(2,1.5);
		\draw[red!80!black,very thick,dashed] (0,0) .. controls (0.3,-0.6) and (0.7,-0.6) .. (1,0);
		\draw[red!80!black,very thick,dashed] (0,0) .. controls (0.8,-1.3) and (2.2,-1.3) .. (3,0);
		\foreach \i in {0,...,4}{\node[pt] at (\i,0) {}; \node[lab,fill=white,inner sep=1pt] at (\i,0.35) {$p_{\i}$};}
		\node[font=\small] at (2,-2.15) {(B)};
	\end{tikzpicture}
	\caption{\small Two Tverberg partitions which respect a graph $H$ (red, dashed), vertices of $\Delta_n$ and their images are both denoted by $p_j$.
		(A) An affine map $a\colon\Delta_6\longrightarrow\R^2$, $r=3$, $d=2$.
		The faces $\{p_0\}$, $\{p_1,p_3,p_5\}$ (blue) and $\{p_2,p_4,p_6\}$ (green) form a Tverberg partition, since both triangles contain $p_0$, and every edge of the graph $H$, a triangle together with two disjoint edges, joins vertices of different faces.
		(B) An affine map $a\colon\Delta_4\longrightarrow\R$, $r=3$, $d=1$.
		Every Tverberg partition consists of the median $p_2$ and two segments each of which contains it, so the faces are $\{p_2\}$, $\{p_0,p_4\}$ (blue) and $\{p_1,p_3\}$ (green), or $\{p_2\}$, $\{p_0,p_3\}$ and $\{p_1,p_4\}$.
		For the path $H$ with the edges $p_0p_1$ and $p_0p_3$ only the first one respects $H$.
		The path with the edges $p_3p_0$ and $p_0p_4$ is respected by none of them, so this path with two edges is not a Tverberg graph for $r=3$, $d=1$ and $n=N=4$.}
	\label{fig : Tverberg graph}
\end{figure}

\medskip
The bounds~\eqref{eq : VZ-Hell} and~\eqref{eq : BB} are obtained by counting Tverberg partitions which satisfy additional constraints on pairs of vertices of $\Delta_N$.
Constrained versions of Tverberg's theorem ask for Tverberg partitions with additional combinatorial properties. For example, we may prescribe pairs of vertices which should lie in different faces, or we can ask for the faces to contain at most one vertex of each colour, for a given colouring of the vertex set of the simplex $\Delta_N$.
In his 2008 paper Hell~\cite{Hell2008constraints} proposed to encode such constraints in the language of graphs, by asking that the two vertices of every edge of a graph lie in different faces.
Such constraint graphs are called  Tverberg graphs, consult~\cite{Engstrom2011} and see Figure~\ref{fig : Tverberg graph} for an illustration.
More precisely:

\medskip
\begin{definition}
	\label{def : Tverberg graph}
	Let $r\geq2$, $d\geq1$ and $n\geq N:=(d+1)(r-1)$ be integers. Let $\mathfrak{F}$ be a family of either affine or continuous maps  $\Delta_n\longrightarrow\R^d$ from the $n$-dimensional simplex to $\R^d$.
	
	\smallskip\noindent
	A graph $H$ on the vertex set of $\Delta_n$ is a {\em Tverberg graph} for $\mathfrak{F}$ if every map $f\in \mathfrak{F}$ possesses a Tverberg partition $(F_1,\dots,F_r)$, that is, a collection of pairwise disjoint faces of $\Delta_n$ with $f(F_1)\cap\dots\cap f(F_r)\neq\emptyset$, such that the faces $F_1,\dots,F_r$ cover all vertices of $\Delta_n$ and the two vertices of every edge of $H$ lie in different faces.
\end{definition}

\medskip
A Tverberg partition in which the two vertices of every edge of $H$ lie in different faces, that is, in which no face contains an edge of $H$, is said to {\em respect} $H$.

\medskip
The first version of Tverberg's theorem with constraints on pairs of vertices appeared in the work of Vu\v{c}i\'c and \v{Z}ivaljevi\'c~\cite{VucicZivaljevic1993}.
For $r$ an odd prime, they considered the Tverberg partitions in which the two vertices of every pair of a fixed matching lie in different faces. Using the connectivity of joins of $1$-spheres and Dold's theorem they showed that every continuous map $\Delta_N\longrightarrow\R^d$ has such a partition, and obtained their lower bound by counting.
In the language of Definition~\ref{def : Tverberg graph}, matchings are Tverberg graphs for $n=N$ and $r$ an odd prime.
Hell~\cite[Thm.\,2]{Hell2008constraints} extended this to prime powers $r\geq 3$, and recent work~\cite[Cor.\,1.8]{BlagojevicSierksma} to all $r\geq3$, for affine maps.
Of course, if $r$ is not a prime power, Tverberg partitions of continuous maps may fail to exist.

\medskip
Colour constraints go back to the coloured Tverberg problem of B\'ar\'any and Larman~\cite{BaranyLarman1992}, which was motivated by the work of B\'ar\'any, F\"uredi and Lov\'asz on halving planes~\cite{BaranyFurediLovasz1990}.
B\'ar\'any and Larman asked for the smallest number of points in $\R^d$ such that every colouring of these points by $d+1$ colours, in which every colour class has at least $r$ points, admits $r$ pairwise disjoint rainbow sets, that is, sets with at most one point of each colour, whose convex hulls have a common point.
This number is at least $(d+1)r$, and B\'ar\'any and Larman conjectured that it equals $(d+1)r$.
They proved the conjecture for $d\leq2$, and for $r=2$ with a topological argument of Lov\'asz.

\medskip
\v{Z}ivaljevi\'c and Vre\'cica~\cite{Zivaljevic1992} modified the problem and asked instead for the smallest size of the colour classes.
They showed, using the connectivity of chessboard complexes~\cite{BjornerLovaszVrecicaZivaljevic}, that for $r$ a prime and $d+1$ colour classes of at least $2r-1$ points each there are $r$ pairwise disjoint rainbow faces whose images have a common point, even for continuous maps.
The extension to prime powers is due to \v{Z}ivaljevi\'c~\cite{Zivaljevic1998}, and by Bertrand's postulate colour classes of $4r-3$ points suffice for every $r$~\cite[Sec.\,4.3]{BlagojevicZiegler2017}.
Since in the problem of B\'ar\'any and Larman a single colour class may remain small, these bounds on the size of the colour classes do not imply that a finite number of points suffices.
On the other hand, the conjecture of B\'ar\'any and Larman is equivalent to the statement that colour classes of exactly $r$ points suffice, because $r$ points of each class can be selected.
Blagojevi\'c, Matschke and Ziegler~\cite[Thm.\,2.1]{Blagojevic2009} proved the optimal coloured Tverberg theorem: if $r$ is a prime and the vertices of $\Delta_N$ are coloured so that every colour class has at most $r-1$ elements, then every continuous map $\Delta_N\longrightarrow\R^d$ has a Tverberg partition into rainbow faces.
Applied with $r+1$ parts to a cone over the simplex, whose apex forms an additional colour class, it implies the conjecture of B\'ar\'any and Larman whenever $r+1$ is a prime~\cite[Cor.\,2.3]{Blagojevic2009}.
Together with Bertrand's postulate this shows that $d+1$ colour classes of $2r-2$ points suffice for every $r$~\cite[Cor.\,2.4]{Blagojevic2009}, which improves the bound $4r-3$ mentioned above.
Recently Sober\'on~\cite[Thm.\,1]{Soberon2026BL} proved the conjecture of B\'ar\'any and Larman, even for continuous maps, whenever $r$ is a prime, building on the optimal coloured Tverberg theorem.
The first open case is $r=8$ and $d=3$, and apart from the cases $d\leq2$, $r=2$, $r$ a prime and $r+1$ a prime, it is not even known whether a finite number of points suffices in the problem of B\'ar\'any and Larman.
In the same paper Sober\'on~\cite[Sec.\,3]{Soberon2026BL} showed that the optimal coloured Tverberg theorem does not extend to prime powers: for $r=4$ and every $d\geq3$ there exist an affine map $\Delta_N\longrightarrow\R^d$ and a colouring of its $N+1=3d+4$ vertices by $d+1$ classes of size $3$ and one class of size $1$ without a Tverberg partition into rainbow faces.
A colouring is encoded by the disjoint union of the complete graphs on its colour classes: a face is rainbow exactly when it contains no edge of this graph, so the Tverberg partitions into rainbow faces are those which respect this graph.
Moreover, a vertex outside the faces of a rainbow Tverberg partition can be added to a face which contains no vertex of its colour, since its colour class has at most $r-1$ elements.
Hence, for $r$ prime and $n=N$, the disjoint unions of complete graphs $K_{\ell}$ with ${\ell}\leq r-1$ are Tverberg graphs for continuous maps.
By the example of Sober\'on, for $r=4$, $d\geq3$ and $n=N$ the disjoint union of $d+1$ triangles and an isolated vertex is not a Tverberg graph, not even for affine maps.

\medskip
Sch\"oneborn and Ziegler~\cite[Thm.\,5.8]{SchoenebornZiegler2005} showed, with winding numbers, that for $d=2$ and $r$ prime certain pairs of vertices can be forced into different faces.
Then Hell~\cite{Hell2008constraints} introduced constraint graphs and gave the first list of examples beyond matchings and colourings.
He proved~\cite[Thm.\,3]{Hell2008constraints} that for every prime power $q>2$ and $n=N$ the complete graphs $K_\ell$ with $2\ell<q+2$, the stars $K_{1,\ell}$ with $\ell<q-1$, paths (when $q>3$), cycles (when $q>4$), and disjoint unions of these are Tverberg graphs for continuous maps.
With computer-generated configurations he also observed that a path with two edges is not a Tverberg graph for $q=3$, and that the cycle with four vertices is not one for $q=4$.

\medskip
Engstr\"om~\cite[Thm.\,3.7]{Engstrom2011} gave a local criterion for Tverberg graphs: for a prime power $q$ and $n=N$, a graph is a Tverberg graph for continuous maps if every vertex $v$ satisfies $q>|N^2(v)|+2|N(v)|$, where $N(v)$ is the set of neighbours of $v$ and $N^2(v)$ the set of vertices at distance two from $v$.
In particular, every graph with maximum degree $\Delta$ where $\Delta(\Delta+1)<q$ is a Tverberg graph~\cite[Cor.\,3.8]{Engstrom2011}.
His proof has two steps.
First, by a topological criterion~\cite[Thm.\,2.2]{Engstrom2011}, a graph is a Tverberg graph if the complex of its partial proper colourings with $q$ colours is $(N-1)$-connected.
Second, the local condition implies the required  connectivity~\cite[Thm.\,3.5]{Engstrom2011}.

\medskip
Further, Engstr\"om~\cite[Conj.\,3.9]{Engstrom2011} conjectured that for a prime power $q$ already on the $(d+1)(q-1)+1$ vertices of $\Delta_N$ the condition $q>K\Delta$ suffices for a graph with maximum degree $\Delta$ to be a Tverberg graph.
Here $K$ is some constant.
Engstr\"om and Nor\'en~\cite{EngstromNoren2014} related Tverberg graphs to graph colouring and proved a version of this conjecture with more vertices, or more precisely, that $q>K_\epsilon\Delta$ suffices when the number of vertices is increased by the factor $1+\epsilon$ and $d$ and $\Delta$ are large~\cite[Cor.\,3.4]{EngstromNoren2014}.

\medskip
In 2014 Blagojevi\'c, Frick and Ziegler~\cite{Blagojevic2014} developed the constraint method, which derives constrained versions directly from the topological Tverberg theorem.
This gives short proofs of the coloured Tverberg theorem of \v{Z}ivaljevi\'c and Vre\'cica and of the generalised van Kampen--Flores theorem, and the same idea, combined with the work of Mabillard and Wagner~\cite{Mabillard2014}, underlies the counterexamples to the topological Tverberg conjecture~\cite{Frick2015,Blagojevic2019}.

\medskip
As far as we know, for $n=N$ and apart from the matchings of~\cite[Cor.\,1.8]{BlagojevicSierksma}, all Tverberg graphs with more than one edge have been obtained for prime powers $r=q$ and continuous maps.
Beyond prime powers, the counterexamples to the topological Tverberg conjecture~\cite{Blagojevic2019,Frick2015} show that for $r$ not a prime power and $d$ large enough there exist continuous maps $\Delta_N\longrightarrow\R^d$ without any Tverberg partition, and so for such $d$ and $n=N$ no graph is a Tverberg graph for continuous maps.
For affine maps and every $r\geq3$, a single edge is a Tverberg graph by~\cite[Cor.\,1.7]{Blagojevic2025}, and every matching is a Tverberg graph by the affine matching theorem~\cite[Cor.\,1.8]{BlagojevicSierksma}.

\medskip
In this paper we prove that several classes of graphs are Tverberg graphs with respect to affine maps $\Delta_n\longrightarrow\R^d$ and every $r\geq3$.
To describe them we recall some notions from graph theory.
A graph is {\em chordal} if every cycle of length at least four has a chord, that is, an edge joining two vertices of the cycle which are not consecutive on it.
The {\em blocks} of a graph are its maximal connected sub-graphs which have no cut vertex, where a {\em cut vertex} of a connected graph is a vertex whose removal disconnects it.
So every block is a single vertex, an edge, or a graph with at least three vertices which stays connected after the removal of any vertex.
Clearly, any two blocks have at most one vertex in common.
For a set $X$ of vertices of a graph, the {\em induced sub-graph} on $X$ has the vertex set $X$ and all edges of the graph with both endpoints in $X$. 
An {\em induced cycle} is an induced sub-graph which is a cycle.
A {\em rooted} graph is a graph together with a vertex in every connected component, the {\em root} of this component.
For more details on these notions see for example~\cite[Sec.\,3.1 and Sec.\,5.5]{Diestel2017}.

\begin{figure}[ht]
	\centering
	\begin{tikzpicture}[x=8mm,y=8mm,v/.style={circle,fill=black,inner sep=1.3pt},c/.style={circle,draw,fill=white,inner sep=1.5pt},lab/.style={font=\small}]

		\begin{scope}
			\node[v] (r) at (0,0.8) {};
			\node[v] (a1) at (-0.8,0) {}; \node[v] (a2) at (0,0) {}; \node[v] (a3) at (0.8,0) {};
			\node[v] (b1) at (-1.2,-0.8) {}; \node[v] (b2) at (-0.4,-0.8) {}; \node[v] (b3) at (0.8,-0.8) {};
			\draw (r)--(a1) (r)--(a2) (r)--(a3) (a1)--(b1) (a1)--(b2) (a3)--(b3);
			\node[lab] at (0,-1.7) {(A)};
		\end{scope}

		\begin{scope}[xshift=30mm]
			\foreach \i in {0,...,5}{\node[v] (h\i) at ({60*\i+90}:0.95) {};}
			\draw (h0)--(h1)--(h2)--(h3)--(h4)--(h5)--(h0);
			\draw (h0)--(h2) (h0)--(h3) (h0)--(h4);
			\node[lab] at (0,-1.7) {(B)};
		\end{scope}

		\begin{scope}[xshift=57mm]
			\foreach \i in {0,...,5}{\node[v] (z\i) at ({60*\i+90}:0.95) {};}
			\draw (z0)--(z1)--(z2)--(z3)--(z4)--(z5)--(z0);
			\node[lab] at (0,-1.7) {(C)};
		\end{scope}

		\begin{scope}[xshift=82mm,yshift=-1mm]
			\node[v] (A) at (0,-0.45) {};
			\node[c] (B) at (0.9,-0.45) {};
			\node[v] (C) at (0.9,0.45) {};
			\node[c] (D) at (0,0.45) {};
			\node[v] (P1) at (1.45,-1.15) {};
			\node[c] (P2) at (2.3,-0.95) {};
			\node[c] (P3) at (2.5,-0.05) {};
			\node[v] (P4) at (1.65,0.25) {};
			\node[v] (T1) at (3.3,0.45) {};
			\node[v] (T2) at (3.35,-0.45) {};
			\node[v] (L1) at (-0.6,1.1) {};
			\node[v] (L2) at (2.85,-1.55) {};
			\draw (A)--(B)--(C)--(D)--(A) (A)--(C) (B)--(D);
			\draw (B)--(P1)--(P2)--(P3)--(P4)--(B);
			\draw (P3)--(T1)--(T2)--(P3);
			\draw (D)--(L1) (P2)--(L2);
			\node[c] at (B) {}; \node[c] at (D) {}; \node[c] at (P2) {}; \node[c] at (P3) {};
			\node[lab] at (1.4,-1.6) {(D)};
		\end{scope}

		\begin{scope}[xshift=127mm]
			\node[v] (u1) at (-0.5,0.7) {}; \node[v] (u2) at (0.5,0.7) {};
			\node[v] (w1) at (-0.9,-0.7) {}; \node[v] (w2) at (0,-0.7) {}; \node[v] (w3) at (0.9,-0.7) {};
			\draw (u1)--(w1)--(u2) (u1)--(w2)--(u2) (u1)--(w3)--(u2);
			\node[lab] at (0,-1.7) {(E)};
		\end{scope}
	\end{tikzpicture}
	\caption{\small Graphs in $\mathscr B$ and a graph which is not in $\mathscr B$.
		(A) A tree.
		(B) A chordal graph: every cycle of length at least four has a chord.
		(C) A cycle.
		(D) A graph with five blocks, the complete graph $K_4$, a cycle of length five, a triangle and two edges. Its cut vertices are drawn as circles.
		(E) The complete bipartite graph $K_{2,3}$ is a single block which is neither chordal, since its cycles of length four have no chord, nor a cycle, so it does not belong to $\mathscr B$.}
	\label{fig : classes}
\end{figure}
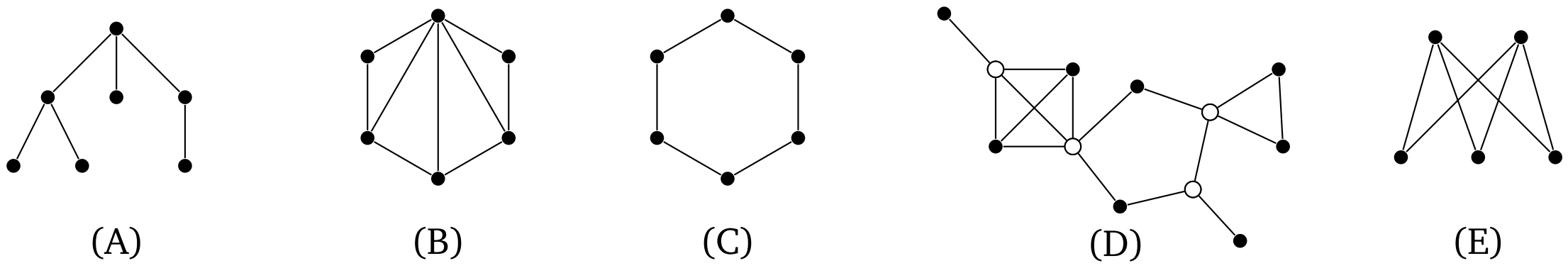

\medskip
Let $\mathscr B$ denote the class of graphs, each of whose blocks is a chordal graph or a cycle.
The class $\mathscr B$ contains:
\begin{compactitem}[\quad --]
	\item Forests, these are graphs without cycles. Their connected components are called {\em trees}. For instance the {\em star} $K_{1,k}$ is a tree in which there exists a vertex, called  the centre, which is adjacent to all of the other $k$ vertices.
	\item Chordal graphs, in particular complete graphs.
	\item Cycles.
	\item Graphs constructed from chordal graphs and cycles by successively identifying one vertex of a new piece with one vertex of the graph built so far.
\end{compactitem}
Disjoint unions of such graphs also belong to $\mathscr B$.
For an illustration see Figure~\ref{fig : classes}.

\medskip
We now state the results precisely.
They are expressed through the largest prime power dividing $r$ and through a deficit of the graph, which measures by how much the degrees of its vertices are too large for the available colours. We begin with the relevant definitions.

\medskip
For an integer $r\geq2$ let $b(r)$ denote the largest prime power which divides $r$.
Given integers $\ell\geq0$ and $b\geq1$, let
\[
	\nu(\ell,b)\ :=\ \min\Big\{\ell,\,b,\,\Big\lfloor\frac{\ell+b+1}{3}\Big\rfloor\Big\}
\]
be the number which governs the connectivity of chessboard complexes.
Recall that the {\em chessboard complex} $\Delta_{\ell,b}$ is the simplicial complex on the vertex set $[\ell]\times[b]$ whose faces are the sets which contain at most one vertex from every row and every column.
It is $(\nu(\ell,b)-2)$-connected~\cite[Thm.\,1.1]{BjornerLovaszVrecicaZivaljevic}.
This bound is best possible as demonstrated by Shareshian and Wachs~\cite{ShareshianWachs2007}.

\medskip
Now, let $H\in\mathscr{B}$ and $b\geq1$.
For a connected component $C$ of $H$ and a vertex $u$ of $C$ set
\[
	\varepsilon_b(C,u)\ :=\ \sum_{w\in C}\max\{0,\,2\deg_H(w)+\chi_u(w)-b\},
\]
where $\chi_u(w)=1$ if $w=u$ or if $w$ lies on an induced cycle of length at least four, and $\chi_u(w)=0$ otherwise.
Note that only the first condition depends on $u$.
Furthermore, let 
\[
\varepsilon_b(H):=\sum_C\varepsilon_b(C),
\] 
where the sum runs over the connected components $C$ of $H$, and $\varepsilon_b(C)$ is given by the first of the following cases which applies:
\[
	\varepsilon_b(C)\ :=\
	\begin{cases}
		\ell-\nu(\ell,b)                & \text{if $C$ is a complete graph $K_\ell$},         \\[2pt]
		0                               & \text{if $C$ is a star $K_{1,k}$ with $k\leq b-2$}, \\[2pt]
		\min_{u\in C}\varepsilon_b(C,u) & \text{otherwise}.
	\end{cases}
\]

\medskip
To illustrate this definition, Table~\ref{tab : deficits} lists the values of $\varepsilon_b(C)$ for the basic connected graphs in the class $\mathscr B$.

\begin{table}[ht]
	\centering\small
	\begin{tabular}{@{}lll@{}}
		\toprule
		connected graph $C$ & $\varepsilon_b(C)$ & $\varepsilon_b(C)=0$ iff\\
		\midrule
		complete graph $K_\ell$ & $\ell-\nu(\ell,b)$ & $b\geq2\ell-1$\\[2pt]
		star $K_{1,k}$, $k\geq2$ & $0$ if $k\leq b-2$, else $2k-b$ & $k\leq b-2$\\[2pt]
		any other tree $T$ & $\sum_{w}\max\{0,\,2\deg_T(w)-b\}$ & $\Delta(T)\leq\frac{b}{2}$\\[2pt]
		cycle $C_m$, $m\geq4$ & $m\cdot\max\{0,\,5-b\}$ & $b\geq5$\\[2pt]
		$C\in\mathscr B$ of maximum degree $\Delta$ & $\leq\sum_{w}\max\{0,\,2\deg_C(w)+1-b\}$ & if $b\geq2\Delta+1$\\
		\bottomrule
	\end{tabular}\vspace{5pt}
	\caption{\small The numbers $\varepsilon_b(C)$ for connected graphs $C$ in $\mathscr B$ and $b\geq3$. For trees and stars the minimum over $u$ is attained at a leaf, and for cycles at any vertex.}
	\label{tab : deficits}
\end{table}

\medbreak
Our first result gives the following criterion for the graphs of the class $\mathscr B$ to be Tverberg graphs.

\medskip
\begin{theorem}
	\label{th : main}
	Let $r\geq2$, $d\geq1$ and $n\geq N:=(d+1)(r-1)$ be integers, and let $H\in\mathscr B$ be a graph on the vertex set of $\Delta_n$ whose maximum degree is at most $r-1$.
	If
	\[
		n\ \geq\ N+\max\{\varepsilon_{r}(H),\,\varepsilon_{b(r)}(H)-1\},
	\]
	then $H$ is a Tverberg graph for affine maps $\Delta_n\longrightarrow\R^d$.
\end{theorem}

\medskip
Since $b(r)\geq3$ for every $r\geq3$, Theorem~\ref{th : main} gives Tverberg graphs with many edges for every $r\geq3$ without additional vertices.

\medskip
\begin{corollary}
	\label{cor : no extra vertices}
	Let $r\geq3$, $d\geq1$ and $n\geq N:=(d+1)(r-1)$ be integers, and let $H$ be a graph on the vertex set of $\Delta_n$ each of whose connected components is
	\begin{compactenum}[\rm\quad (A)]
		\item a graph in $\mathscr B$ whose maximum degree is at most $\frac{b(r)-1}{2}$, or
		\item a tree of maximum degree at most $\frac{b(r)}{2}$, or
		\item a star $K_{1,k}$ with $k\leq b(r)-2$.
	\end{compactenum}
	Then $H$ is a Tverberg graph for affine maps $\Delta_n\longrightarrow\R^d$.
\end{corollary}

\medskip
For instance, for $r=10$ we have $b(10)=5$, and so every graph of maximum degree at most $2$, that is every disjoint union of paths and cycles, is a Tverberg graph for affine maps, and so are disjoint unions of such graphs with stars $K_{1,3}$.
If $r=2p$ with an odd prime $p$, then the graphs in $\mathscr B$ of maximum degree at most $\frac{p-1}{2}$ are Tverberg graphs, among them forests of this maximum degree and complete graphs on $\frac{p+1}{2}$ vertices.
Table~\ref{tab : degrees} lists the bounds for $r\leq30$ not a prime power.
Theorem~\ref{th : main} allows, in addition to the components of Corollary~\ref{cor : no extra vertices}, one further component $C$ with $\varepsilon_{b(r)}(C)=1$ and $\varepsilon_r(C)=0$.
For example, if $b(r)$ is even and $r\neq b(r)$, this component can be the complete graph on $\frac{b(r)}{2}+1$ vertices. So for $r=12$ a triangle together with a matching is a Tverberg graph for affine maps.

\begin{table}[ht]
	\centering\small
	\begin{tabular}{@{}l|ccccccccccccc@{}}
		\toprule
		$r$                             & $6$ & $10$ & $12$ & $14$ & $15$ & $18$ & $20$ & $21$ & $22$ & $24$ & $26$ & $28$ & $30$ \\
		$b(r)$                          & $3$ & $5$  & $4$  & $7$  & $5$  & $9$  & $5$  & $7$  & $11$ & $8$  & $13$ & $7$  & $5$  \\
		\midrule
		(A) maximum degree              & $1$ & $2$  & $1$  & $3$  & $2$  & $4$  & $2$  & $3$  & $5$  & $3$  & $6$  & $3$  & $2$  \\
		(B) maximum degree of trees    & $1$ & $2$  & $2$  & $3$  & $2$  & $4$  & $2$  & $3$  & $5$  & $4$  & $6$  & $3$  & $2$  \\
		(C) number of leaves of stars & $1$ & $3$  & $2$  & $5$  & $3$  & $7$  & $3$  & $5$  & $9$  & $6$  & $11$ & $5$  & $3$  \\
		\bottomrule
	\end{tabular}\vspace{5pt}
	\caption{\small The bounds of Corollary~\ref{cor : no extra vertices} for $r\leq30$ not a prime power.}
	\label{tab : degrees}
\end{table}

\medskip
For a disjoint union of complete graphs, being a Tverberg graph means that the vertex sets of the complete graphs, viewed as colour classes, can be made rainbow: every face of the Tverberg partition contains at most one vertex of every class.
In this way Theorem~\ref{th : main} gives a coloured Tverberg theorem for every $r$, compare~\cite{Zivaljevic1992,Blagojevic2009}.

\medskip
\begin{corollary}
	\label{cor : coloured}
	Let $r\geq2$, $d\geq1$ and $n\geq N:=(d+1)(r-1)$ be integers, and let the vertex set of $\Delta_n$ be partitioned into colour classes $C_1,\dots,C_m$.
	If
	\[
		\sum_{i=1}^m\nu(|C_i|,r)\ \geq\ N+1
		\qquad\text{and}\qquad
		\sum_{i=1}^m\nu(|C_i|,b(r))\ \geq\ N,
	\]
	then every affine map $\Delta_n\longrightarrow\R^d$ has a Tverberg partition $(F_1,\dots,F_r)$ into rainbow faces.
	If, in addition, $|C_i|\leq r$ for all $1\leq i\leq m$, then the faces $F_1,\dots,F_r$ can be chosen to cover all vertices of $\Delta_n$.
	When $r$ is a prime power, the same holds for continuous maps $\Delta_n\longrightarrow\R^d$.
\end{corollary}

\medskip
Both conditions of Corollary~\ref{cor : coloured} hold when every colour class has at most $\frac{b(r)+1}{2}$ elements, since then $\nu(|C_i|,b(r))=\nu(|C_i|,r)=|C_i|$ and $\sum_{i=1}^m|C_i|=n+1\geq N+1$.
Thus, for every colouring of the $N+1$ vertices of $\Delta_N$ with colour classes of at most $\frac{b(r)+1}{2}$ elements, every affine map $\Delta_N\longrightarrow\R^d$ has a Tverberg partition into rainbow faces which uses all vertices.
For example, for $r=10$ and $d=2$ this applies to $28$ points in the plane coloured with classes of at most three points.
If $r=q$ is a prime power and there are $d+1$ colour classes of at least $2q-1$ elements each, then $\nu(|C_i|,q)=q$ and $(d+1)q\geq N+1$, so Corollary~\ref{cor : coloured} contains the coloured Tverberg theorem of \v{Z}ivaljevi\'c and Vre\'cica~\cite{Zivaljevic1992,Zivaljevic1998}.
On the other hand, when $r$ is not a prime power, the second condition forces $m\geq\frac{N}{b(r)}>d+1$ colour classes, because $\nu(\ell,b)\leq b$ and $r\geq2b(r)$.
Some restriction on the sizes of the colour classes is necessary for $n=N$.
Indeed, for $r=4$ and $n=N$ the conditions of Corollary~\ref{cor : coloured} hold exactly when all colour classes have at most two elements, while by the example of Sober\'on~\cite[Sec.\,3]{Soberon2026BL} classes with three elements cannot be allowed in general.

\medskip
The proof of Theorem~\ref{th : main} works for every prime power $q$ dividing $r$ in place of $b(r)$.
If $q$ is small compared with $r$, one more vertex can be saved.
For integers $b\geq1$ and $m\geq b$ let
\[
	\varepsilon_{b,m}(H):=\sum_C\varepsilon_{b,m}(C),
\]
where the sum again runs over the connected components $C$ of $H$, and
\[
	\varepsilon_{b,m}(C)\ :=\
	\begin{cases}
		\min\{\,\ell-\nu(\ell,b),\ \min_{u\in C}\varepsilon_m(C,u)\,\} & \text{if $C$ is a complete graph $K_\ell$},         \\[2pt]
		0                                                              & \text{if $C$ is a star $K_{1,k}$ with $k\leq m-2$}, \\[2pt]
		\min_{u\in C}\varepsilon_m(C,u)                                & \text{otherwise}.
	\end{cases}
\]
This is the analogue of $\varepsilon_b(H)$ for lists which are unions of parts of size $b$ and have at least $m$ elements, see Proposition~\ref{prop : condition B}. In particular, $\varepsilon_{b,b}(H)\leq\varepsilon_b(H)$.

\medskip
\begin{theorem}
	\label{th : small q}
	Let $r\geq2$, $d\geq1$ and $n\geq N:=(d+1)(r-1)$ be integers, let $q$ be a prime power dividing $r$ with $r\geq(d+2)q$, and let $H\in\mathscr B$ be a graph on the vertex set of $\Delta_n$ whose maximum degree is at most $r-1$.
	If
	\[
		n\ \geq\ N+\max\{\varepsilon_{r}(H),\,\varepsilon_{q}(H)-2,\,\varepsilon_{q,\,r-dq}(H)-1\},
	\]
	then $H$ is a Tverberg graph for affine maps $\Delta_n\longrightarrow\R^d$.
\end{theorem}

\medskip
Compared with Theorem~\ref{th : main} for the prime power $q$, the second term is smaller by one, and the third term is small, since the lists it refers to have at least $r-dq\geq2q$ elements.
The proof rests on a refined form of the $G$-box criterion, Theorem~\ref{th : refined}, in which the full connectivity is needed only for the $G$-boxes whose strata on the fixed sphere are closed, and on Theorem~\ref{th : closed strata}, which shows that, for affine maps whose image spans $\R^d$, in these $G$-boxes every column contains at least $\frac{r}{q}-d$ orbits of the group.
For example, for $r=20$, $q=4$ and $d\leq3$ the complete graph $K_5$ has $\varepsilon_{20}(K_5)=0$, $\varepsilon_4(K_5)=5-\nu(5,4)=2$ and $\varepsilon_{4,20-4d}(K_5)\leq\varepsilon_{4,8}(K_5)=\varepsilon_8(K_5,u)=1$, so $K_5$ is a Tverberg graph for affine maps with $n=N$, while Theorem~\ref{th : main} needs $n\geq N+1$ for every prime power dividing $20$.
Table~\ref{tab : complete} lists all cases with $r\leq30$ in which Theorem~\ref{th : small q} gives a larger complete Tverberg graph with $n=N$ than Theorem~\ref{th : main}. There is none for $r\in\{6,10,14,18,22,26\}$.
For disjoint unions the deficits add up, and the gain is one vertex in total. As an illustration, for $r=12$, $q=3$ and $d\leq2$ two triangles together with a matching are a Tverberg graph for affine maps with $n=N$.

\begin{table}[ht]
	\centering\small
	\begin{tabular}{@{}cccccc@{}}
		\toprule
		$r$  & $b(r)$ & Theorem~\ref{th : main}: $K_\ell$ for $\ell\leq$ & $q$        & dimensions           & Theorem~\ref{th : small q}: $K_\ell$ for $\ell\leq$ \\
		\midrule
		$12$ & $4$    & $4$                                              & $3$ or $4$ & $d=1$                & $5$                                                 \\
		$15$ & $5$    & $4$                                              & $3$; \ $5$ & $d\leq2$; \ $d=1$    & $5$; \ $6$                                          \\
		$20$ & $5$    & $4$                                              & $4$; \ $5$ & $d\leq3$; \ $d\leq2$ & $5$; \ $6$                                          \\
		$21$ & $7$    & $5$                                              & $7$        & $d=1$                & $7$                                                 \\
		$24$ & $8$    & $6$                                              & $8$        & $d=1$                & $7$                                                 \\
		$28$ & $7$    & $5$                                              & $7$        & $d\leq2$             & $7$                                                 \\
		$30$ & $5$    & $4$                                              & $3$; \ $5$ & $d\leq7$; \ $d\leq4$ & $5$; \ $6$                                          \\
		\bottomrule
	\end{tabular}\vspace{5pt}
	\caption{\small Complete graphs $K_\ell$ which are Tverberg graphs for affine maps with $n=N$ by Theorem~\ref{th : main} and by Theorem~\ref{th : small q}, for $r\leq30$ not a prime power. In each row the largest $\ell$ is given, for Theorem~\ref{th : small q} for the stated $q$ and dimensions.}
	\label{tab : complete}
\end{table}

\medskip
For prime powers the method also applies to continuous maps.

\medskip
\begin{corollary}
	\label{cor : prime powers}
	Let $q$ be a prime power, $d\geq1$ and $n\geq N:=(d+1)(q-1)$ be integers, and let $H\in\mathscr B$ be a graph on the vertex set of $\Delta_n$ whose maximum degree is at most $q-1$.
	If $n\geq N+\varepsilon_q(H)$, then $H$ is a Tverberg graph for continuous maps $\Delta_n\longrightarrow\R^d$.
	In particular, if $n\geq N$ and $q\geq3$, then every graph whose connected components are as in Corollary~\ref{cor : no extra vertices}, with $b(r)$ replaced by $q$, is a Tverberg graph for continuous maps $\Delta_n\longrightarrow\R^d$.
\end{corollary}

\medskip
Corollary~\ref{cor : prime powers} contains all examples of Hell~\cite[Thm.\,3]{Hell2008constraints}: complete graphs $K_\ell$ with $2\ell<q+2$, stars $K_{1,\ell}$ with $\ell<q-1$, paths for $q>3$, cycles for $q>4$, and their disjoint unions.
In turn, Corollary~\ref{cor : no extra vertices} extends these examples to affine maps and every $r\geq3$, with $q$ replaced by $b(r)$.
For $K_\ell$ with $\ell\leq q<2\ell-1$ Corollary~\ref{cor : prime powers} needs only $\ell-\lfloor\frac{\ell+q+1}{3}\rfloor$ additional vertices, for instance $\lceil\frac{q-1}{3}\rceil$ additional vertices for $K_q$.
If $q$ is a prime and $\ell\leq q-1$, no additional vertices are needed, by the optimal coloured Tverberg theorem~\cite{Blagojevic2009}. By the example of Sober\'on~\cite[Sec.\,3]{Soberon2026BL} this does not extend to $q=4$, where Corollary~\ref{cor : prime powers} needs one additional vertex for every triangle.
In the case of forests the bound $q\geq2\Delta$ is linear in the maximum degree, while Engstr\"om's corollary~\cite[Cor.\,3.8]{Engstrom2011} requires $\Delta(\Delta+1)<q$.
Moreover, Corollary~\ref{cor : prime powers} confirms the conjecture of Engstr\"om~\cite[Conj.\,3.9]{Engstrom2011} with $K=2$, for prime powers $q$ and for all graphs in $\mathscr B$.

\medskip
Our last result extends Engstr\"om's local criterion to all $r$, for affine maps, and it applies to graphs which are not in $\mathscr B$.

\medskip
\begin{theorem}
	\label{th : local}
	Let $r\geq2$, $d\geq1$ and $n\geq N:=(d+1)(r-1)$ be integers, and let $H$ be a graph on the vertex set of $\Delta_n$ such that
	$
		|N^2(v)|+2|N(v)|\ <\ b(r)
	$
	for every vertex $v$ of $H$.
	Then $H$ is a Tverberg graph for affine maps $\Delta_n\longrightarrow\R^d$.
	If $r$ is a prime power, then $H$ is a Tverberg graph for continuous maps $\Delta_n\longrightarrow\R^d$.
\end{theorem}

\medskip
For prime powers $r=q$ and $n=N$ this is Engstr\"om's criterion~\cite[Thm.\,3.7]{Engstrom2011}.
Since $|N^2(v)|\leq\Delta(\Delta-1)$, every graph with maximum degree $\Delta$ and $\Delta(\Delta+1)<b(r)$ is a Tverberg graph for affine maps.
For instance, the complete bipartite graph $K_{3,3}$ is a Tverberg graph for affine maps when $b(r)\geq9$, for example for $r=18$.
Further, the graph of the $3$-cube is a Tverberg graph for affine maps when $b(r)\geq10$, for example for $r=22$, and so is the Petersen graph when $b(r)\geq13$, for example for $r=26$.
None of these graphs belongs to $\mathscr B$.
On the class $\mathscr B$ the bounds of Corollary~\ref{cor : no extra vertices} are better, since they are linear in the maximum degree.

\medskip
The remaining sections are organised as follows.
Section~\ref{sec : preliminaries} recalls the deleted join and the $G$-box criterion from~\cite{BlagojevicSierksma}, proves a refined form of the criterion, and shows that closed strata on the fixed sphere force many orbits in every column.
In Section~\ref{sec : list colourings} we give lower bounds for the connectivity of complexes of partial list colourings of graphs in $\mathscr B$, of stars, and of complete graphs whose lists are unions of parts of a partition of the colours, and a version for lists of Engstr\"om's local connectivity bound.
Section~\ref{sec : proofs} contains the proofs of the results stated above.

\bigskip
\section{The \texorpdfstring{$G$}{G}-box criterion and its refinement}
\label{sec : preliminaries}
\bigskip

In this section we recall relevant notions and results from~\cite{BlagojevicSierksma}.
Additionally, we sharpen the $G$-box criterion into the form needed in the proofs of our main results.

\medskip
We use two general notions throughout the paper.
For a simplicial complex $X$ on a vertex set $W$ and a subset $U\subseteq W$, the sub-complex of $X$ {\em induced} on $U$ is
\[
	X[U]\ :=\ \{F\in X : F\subseteq U\},
\]
the complex of all faces of $X$ whose vertices lie in $U$.
A space is {\em homologically $\tau$-connected} over a field $\FF$ if its reduced homology with coefficients in $\FF$ vanishes in all degrees $i$ with $-1\leq i\leq\tau$. In particular, it is non-empty if $\tau\geq-1$.

\medskip
\subsection{The \texorpdfstring{$G$}{G}-box criterion}
\label{sec : deleted join}

Let $r\geq 2$, $d\geq 1$ and $n\geq0$ be integers, let $N:=(d+1)(r-1)$, and let $p_0,\dots,p_n$ be the vertices of the simplex $\Delta_n$.
The vertices of the $i$-th copy of $\Delta_n$ in the join $(\Delta_n)^{*r}$ are denoted by $p_{0,i},\dots,p_{n,i}$, and we set
\[
	\K:=(\Delta_n)^{*r}_{\Delta(2)}\ \cong\ [r]^{*(n+1)}=Q_0*\cdots*Q_n,\qquad Q_j:=\{p_{j,1},\dots,p_{j,r}\}.
\]
Here $(\Delta_n)^{*r}_{\Delta(2)}$ denotes the $r$-fold $2$-wise deleted join of the simplex $\Delta_n$ with the vertex set $V:=Q_0\sqcup\cdots\sqcup Q_n$.
For more on deleted joins and the isomorphism above we refer to~\cite[Sec.\,5.5 and Sec.\,6.3]{Matousek2003} and~\cite[Sec.\,3.2.2]{BlagojevicZiegler2017}.
A typical face of the simplicial complex $\K$ is a set $J\subseteq V$ with $|J\cap Q_j|\leq1$ for all $j$.
It can be presented as $J=F_1*\cdots*F_r$, where $F_i$ is the face of $\Delta_n$ spanned by the vertices $p_j$ with $p_{j,i}\in J$.
Thus a face $J$ of $\K$ is the same as a partial assignment of labels $i\in[r]$ to the vertices of $\Delta_n$: the vertex $p_j$ receives the label $i$ if $p_{j,i}\in J$, and no label if $J\cap Q_j=\emptyset$.
In particular, the facets of $\K$ are the assignments of labels to all vertices.
For example, the partition $(\{p_2\},\{p_0,p_4\},\{p_1,p_3\})$ of Figure~\ref{fig : Tverberg graph}(B) is the facet $\{p_{0,2},p_{1,3},p_{2,1},p_{3,3},p_{4,2}\}$, which labels $p_2$ by $1$, the vertices $p_0$ and $p_4$ by $2$, and $p_1$ and $p_3$ by $3$.

\medskip
For an affine or continuous map $f\colon\Delta_n\longrightarrow\R^d$, a face $F_1*\cdots*F_r$ of $\K$ is a {\em Tverberg face} of $f$ if $(F_1,\dots,F_r)$ is a Tverberg partition of $f$.
Since enlarging the faces $F_i$ enlarges the sets $f(F_i)$, every face of $\K$ which contains a Tverberg face of a map $f$ is again a Tverberg face of $f$, compare~\cite[Lem.\,2.6]{BlagojevicSierksma}.

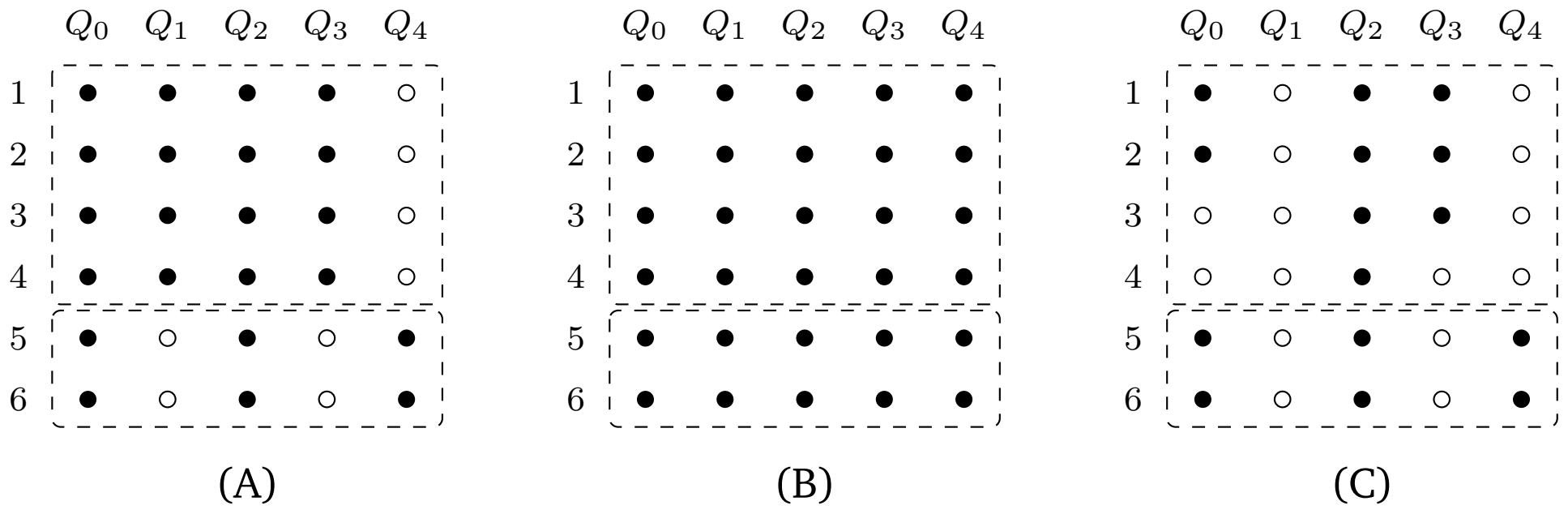
\begin{figure}[ht]
	\centering
	\begin{tikzpicture}[x=6.5mm,y=5mm,
			pt/.style={circle,fill=black,inner sep=1.4pt},
			hollow/.style={circle,draw,fill=white,inner sep=1.3pt},
			lab/.style={font=\scriptsize}]
		\foreach \X/\P/\T in {0/{1/1,1/0,1/1,1/0,0/1}/(A),7/{1/1,1/1,1/1,1/1,1/1}/(B),14/{2/1,0/0,1/1,3/0,0/1}/(C)}{
			\begin{scope}[xshift=\X*6.5mm]
				\draw[dashed,rounded corners=2pt] (-0.45,-0.55) rectangle (4.45,-4.45);
				\draw[dashed,rounded corners=2pt] (-0.45,-4.55) rectangle (4.45,-6.45);
				\foreach \i in {1,...,6}{\node[lab,anchor=east] at (-0.55,-\i) {$\i$};}
				\foreach \A/\B [count=\j from 0] in \P {
					\node[lab] at (\j,0.1) {$Q_{\j}$};
					\foreach \i in {1,2,3,4}{
						\ifnum\A=1 \node[pt] at (\j,-\i) {};
						\else\ifnum\A=0 \node[hollow] at (\j,-\i) {};
						\else\ifnum\A=2 \ifnum\i<3 \node[pt] at (\j,-\i) {}; \else \node[hollow] at (\j,-\i) {}; \fi
						\else \ifnum\i=4 \node[hollow] at (\j,-\i) {}; \else \node[pt] at (\j,-\i) {}; \fi\fi\fi\fi}
					\foreach \i in {5,6}{\ifnum\B=1 \node[pt] at (\j,-\i) {}; \else \node[hollow] at (\j,-\i) {}; \fi}}
				\node[font=\small] at (2,-7.4) {\T};
			\end{scope}}
	\end{tikzpicture}
	\caption{\small Subsets of $V$ for $r=6$, $n=4$ and the group $G\cong(\Z/2)^2$ generated by $(1\,2)(3\,4)(5\,6)$ and $(1\,3)(2\,4)$, whose orbits $\{1,2,3,4\}$ and $\{5,6\}$ are framed, and filled vertices belong to the set.
		(A) and (B) are $G$-boxes: every column contains $\{1,2,3,4\}$, $\{5,6\}$ or $[6]$, and (B) is the $G$-box $V$.
		(C) is not a $G$-box: the columns $Q_0$ and $Q_3$ contain only a part of the orbit $\{1,2,3,4\}$, and $Q_1$ contains no vertex at all.}
	\label{fig : G-box}
\end{figure}

\medskip
A subgroup $G$ of the symmetric group $\sym_r$ acts on the labels, and hence on $V$ by $g\cdot p_{j,i}:=p_{j,g(i)}$ and consequently on $\K$.
An {\em elementary abelian group} is a group isomorphic to $(\Z/p)^\alpha$ for a prime $p$ and some $\alpha\geq0$.
For an elementary abelian group $G\leq\sym_r$, a set $U\subseteq V$ is a {\em $G$-box} if every $U\cap Q_j$, viewed as a subset of $[r]$, is a non-empty union of $G$-orbits.
As an example, let $r=6$ and let $G\cong(\Z/2)^2$ be generated by the permutations $(1\,2)(3\,4)(5\,6)$ and $(1\,3)(2\,4)$. Its orbits are $\{1,2,3,4\}$ and $\{5,6\}$.
Then a $G$-box chooses, independently in every column, one of the sets $\{1,2,3,4\}$, $\{5,6\}$ and $[6]$, as illustrated in Figure~\ref{fig : G-box}.

\medskip
Now we quote the $G$-box criterion~\cite[Thm.\,3.10]{BlagojevicSierksma}, whose refinement will be introduced and proved in Section~\ref{sec : refined}.

\begin{theorem}[$G$-box criterion]
	\label{th : G-box}
	Let $r\geq2$, $d\geq1$ and $n\geq0$ be integers, $N:=(d+1)(r-1)$, $p$ a prime, $G\leq\sym_r$ an elementary abelian group isomorphic to $(\Z/p)^\alpha$ for some $\alpha\geq0$, and let $a\colon\Delta_n\longrightarrow\R^d$ be an affine map.
	Assume that $\LL\subseteq\K$ is a $G$-invariant sub-complex such that $\LL$ is homologically $(N-1)$-connected over $\FF_p$, and $\LL[U]$ is homologically $(N-2)$-connected over $\FF_p$ for every $G$-box $U$.
	
	\smallskip\noindent
	Then $\LL$ contains a Tverberg face of $a$.
	If $G$ has only one orbit, then the same holds for every continuous map $\Delta_n\longrightarrow\R^d$.
\end{theorem}

Theorem~\ref{th : G-box} and its refinement, Theorem~\ref{th : refined}, do not relate $n$ to $N$.
In our applications $n\geq N$, and for $n<N$ the first assumption on $\LL$ says that $\LL$ is acyclic over $\FF_p$, because $\dim\LL\leq n\leq N-1$.

\medskip
We outline the proof of Theorem~\ref{th : G-box} given in~\cite[Sec.\,2 and Sec.\,3]{BlagojevicSierksma} because its ingredients are used for the proof of our main topological result, Theorem~\ref{th : refined}.
 
\begin{proof}[Proof outline]
Let $W_r:=\{y\in\R^r : y_1+\cdots+y_r=0\}$ and $\E:=W_r^{\oplus(d+1)}$, viewed as the space of real $(d+1)\times r$ matrices whose rows sum to zero, with the standard scalar product and unit sphere $S(\E)\cong S^{N-1}$.
The {\em test map} $A\colon(\Delta_n)^{*r}\longrightarrow\E$ of the affine map $a$ is affine on every simplex of the join, and on the vertices it is given as follows: the $\ell$-th column of the matrix $A(p_{j,i})$ is $(\delta_{i\ell}-\tfrac1r)\,\widehat{a}(p_j)$, where $\widehat{a}(x):=(1,a(x))\in\R^{d+1}$ and $\delta_{i\ell}$ is the Kronecker delta.
A face $J$ of $\K$ is a Tverberg face of $a$ if and only if $0\in A(J)$, see~\cite[Lem.\,2.6]{BlagojevicSierksma}.	

\medskip
The group $G$ acts on $\E$ by permuting the columns, $A$ is $G$-equivariant, and $S(\E^G)$ denotes the unit sphere of the fixed subspace $\E^G$, the so-called fixed sphere.
For $y\in S(\E)$ let
$
	U_{a,y}\ :=\ \{v\in V : \langle y,A(v)\rangle\leq0\}.
$
Every $U_{a,y}$ meets every column $Q_j$, and if $y\in S(\E^G)$, then $U_{a,y}$ is a $G$-box~\cite[Lem.\,2.7 and Lem.\,3.2]{BlagojevicSierksma}.

\medbreak
The incidence space of $\LL$ is
\[
	T_{a,\LL}\ :=\ \{(y,\lambda)\in S(\E)\times\LL : \supp(\lambda)\subseteq U_{a,y}\},
\]
where $\supp(\lambda)$ is the face of $\LL$ whose relative interior contains $\lambda$.
It is a compact $G$-space with the equivariant projections $\pr_1\colon T_{a,\LL}\longrightarrow S(\E)$ and $\pr_2\colon T_{a,\LL}\longrightarrow\LL$.
The fibre of $\pr_1$ over $y\in S(\E)$ is just $\LL[U_{a,y}]$.

\medskip
Suppose now that $G$ has at least two orbits and that $\LL$ contains no Tverberg face of $a$.
All cohomology is equivariant cohomology of Borel constructions with coefficients in $\FF_p$, and $\pi_X\colon X_G\longrightarrow\mathrm{B}G$ denotes the projection of the Borel construction of a $G$-space $X$.
The proof proceeds in six steps through a chase inside the commutative diagram
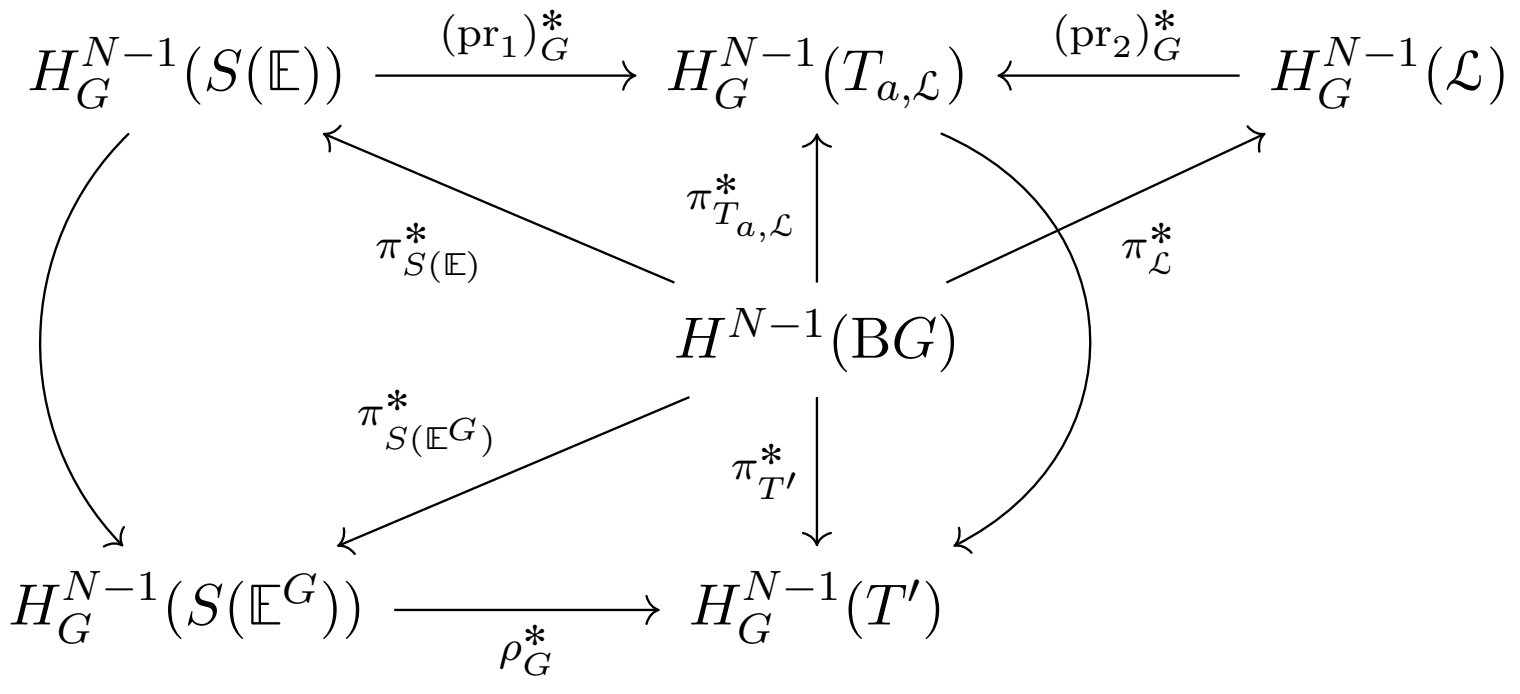
\[
\begin{tikzcd}[column sep=4.2em,row sep=2.6em]
	H^{N-1}_G(S(\E)) \arrow[r,"(\pr_1)_G^*"] \arrow[dd,bend right=45] & H^{N-1}_G(T_{a,\LL}) \arrow[dd,bend left=65,looseness=1.3] & H^{N-1}_G(\LL) \arrow[l,"(\pr_2)_G^*"'] \\
	& H^{N-1}(\mathrm{B}G) \arrow[ul,"\pi_{S(\E)}^*"] \arrow[u,"\pi_{T_{a,\LL}}^*"] \arrow[ur,"\pi_{\LL}^*"'] \arrow[dl,"\pi_{S(\E^G)}^*"'] \arrow[d,"\pi_{T'}^*"'] & \\
	H^{N-1}_G(S(\E^G)) \arrow[r,"\rho_G^*"'] & H^{N-1}_G(T') &
\end{tikzcd}
\]
in which the curved arrows are induced by the inclusions $S(\E^G)\subseteq S(\E)$ and $T'\subseteq T_{a,\LL}$ introduced in Step~4.

\medskip
	{\em Step 1:}
The fibres of $\pr_2$ are contractible, since no face of $\LL$ is a Tverberg face, and by the equivariant Vietoris--Begle theorem~\cite[Lem.\,3.5]{BlagojevicSierksma} the map $(\pr_2)_G^*\colon H^*_G(\LL)\longrightarrow H^*_G(T_{a,\LL})$ is an isomorphism.

\medskip
	{\em Step 2:}
Since $\LL$ is homologically $(N-1)$-connected, $\pi_{ \LL}^*\colon H^{N-1}(\mathrm{B}G)\longrightarrow H^{N-1}_G( \LL)$ is an isomorphism which implies that $\pi_{T_{a,\LL}}^*\colon H^{N-1}(\mathrm{B}G)\longrightarrow H^{N-1}_G(T_{a,\LL})$ is also an isomorphism.

\medskip
	{\em Step 3:}
There is a class $\Theta\in H^{N-1}_G(S(\E))$, the equivariant Poincar\'e dual of a fixed point, whose restriction to the fixed sphere is $e\times[S(\E^G)]^*\in H^{N-1}(\mathrm{B}G\times S(\E^G))$, where $e\neq0$ is the Euler class of the representation $(\E^G)^\perp$, or $e=1$ if $G$ is trivial, see~\cite[Lem.\,3.7]{BlagojevicSierksma}.
	By Step~2 there is $\omega\in H^{N-1}(\mathrm{B}G)$ with $(\pr_1)_G^*\Theta=\pi_{T_{a,\LL}}^*\omega$.	

\medskip
	{\em Step 4:}
Let $T':=\pr_1^{-1}(S(\E^G))$ and let $\rho\colon T'\longrightarrow S(\E^G)$ be the restriction of $\pr_1$.
	Restricting the identity of Step~3 to $T'$ gives $\rho_G^*(e\times[S(\E^G)]^*)=\rho_G^*(\omega\times1)$ in $H^{N-1}_G(T')$.

\medskip
	{\em Step 5:}
The fibres of $\rho$ are the complexes $\LL[U_{a,y}]$ for $y\in S(\E^G)$, which are homologically $(N-2)$-connected by assumption, because every $U_{a,y}$ is a $G$-box.
	By the equivariant Vietoris--Begle theorem $\rho_G^*$ is injective in degree $N-1$, and so $e\times[S(\E^G)]^*=\omega\times1$.

\medskip
	{\em Step 6:}
Since $\dim S(\E^G)\geq1$, the two sides of this equation lie in different summands of the K\"unneth decomposition of $H^{N-1}(\mathrm{B}G\times S(\E^G))$, so both vanish. 
This contradicts $e\neq0$.	

\medskip
If $G$ has only one orbit, then $\E^G=0$, and the proof uses only the equivariant map $\frac{A}{\|A\|}\colon\LL\longrightarrow S(\E)$ together with the non-vanishing of the Euler class of $\E$.
This is the argument of \"Ozaydin and Volovikov~\cite{Oezaydin1987,Volovikov1996-1}, and it applies to continuous maps as well~\cite[Sec.\,3.4]{BlagojevicSierksma}.
\end{proof}

\medskip
In the following, the {\em orbits} of a subgroup $G\leq\sym_r$ are always its orbits on the set of labels $[r]$. 
They form a partition of $[r]$.
Recall that $b(r)$ denotes the largest prime power which divides $r$.
The following lemma provides elementary abelian groups all of whose orbits have $b(r)$ elements, and it shows that no elementary abelian group does better.

\medskip
\begin{lemma}
	\label{lem : group}
	Let $r\geq2$ be an integer.
	\begin{compactenum}[\rm\quad  (A)]
		\item Let $q=p^\alpha$ with $\alpha\geq1$ be a prime power dividing $r$, for instance $q=b(r)$, let $k:=\frac{r}{q}$, and identify $[r]$ with $\FF_p^{\,\alpha}\times[k]$ by a bijection.
		Then the additive group $\FF_p^{\,\alpha}\cong(\Z/p)^\alpha$ acts on $[r]$ by translation of the first coordinate, $g\cdot(u,\iota):=(u+g,\iota)$.
		This action identifies $\FF_p^{\,\alpha}$ with an elementary abelian subgroup $G\leq\sym_r$, $G\cong(\Z/p)^\alpha$, whose orbits are the $k$ sets $\FF_p^{\,\alpha}\times\{\iota\}$, $\iota\in[k]$, each with exactly $q$ elements.
		\item Every elementary abelian subgroup $G\leq\sym_r$ has an orbit with at most $b(r)$ elements.
	\end{compactenum}
\end{lemma}
\begin{proof}
	(A) Translation by $g\neq0$ moves every label, so only $g=0$ acts as the identity permutation of $[r]$. In other words, the action is faithful.
	Hence the action identifies $\FF_p^{\,\alpha}$ with a subgroup $G$ of $\sym_r$ isomorphic to $(\Z/p)^\alpha$.
	The orbit of $(u,\iota)$ is $\FF_p^{\,\alpha}\times\{\iota\}$.
	This is the group of~\cite[Sec.\,4]{BlagojevicSierksma} with $s=0$.

	\smallskip
	(B) Let $G\cong(\Z/p)^\beta$.
	Every orbit of $G$ has a power of $p$ as its size, so the smallest orbit size $p^\gamma$ divides all orbit sizes and hence their sum $r$.
	Therefore $p^\gamma\leq b(r)$.
\end{proof}

\medskip
By Lemma~\ref{lem : group}(B), every elementary abelian group $G\leq\sym_r$ has an orbit $O$ with at most $b(r)$ labels, and the $G$-box $U$ with $U\cap Q_j=O$ for all $j$ gives every vertex of $\Delta_n$ at most $b(r)$ labels.
So, with elementary abelian groups, $b(r)$ is the largest number of labels which can be guaranteed in every column of every $G$-box.
In particular, the fixed labels $\infty_1,\dots,\infty_s$ of~\cite[Thm.\,1.7]{BlagojevicSierksma} do not help for Tverberg graphs in general.
There the group $\FF_p^{\,\alpha}$ acts on the labels $(\FF_p^{\,\alpha}\times[k])\sqcup\{\infty_1,\dots,\infty_s\}$ by translation of the first coordinate and fixes every $\infty_\ell$, so every $\{\infty_\ell\}$ is an orbit with one element.

\medskip
\subsection{The refined \texorpdfstring{$G$}{G}-box criterion}
\label{sec : refined}

The $G$-box criterion asks for the connectivity of $\LL[U]$ for all $G$-boxes $U$, but in its proof, outlined in Section~\ref{sec : deleted join}, only the $G$-boxes $U_{a,y}$ with $y\in S(\E^G)$ occur, in Step~5.
For an affine map $a\colon\Delta_n\longrightarrow\R^d$ and a $G$-box $U$ we call
\[
	N_G(U)\ :=\ \{y\in S(\E^G) : U_{a,y}=U\}
\]
the {\em stratum} of $U$.
Clearly, $N_G(U)$ depends on $a$.

\medskip
The next theorem shows that the full connectivity is needed only for the $G$-boxes whose strata are closed.

\medskip
\begin{theorem}[Refined $G$-box criterion]
	\label{th : refined}
	Let $r\geq2$, $d\geq1$ and $n\geq0$ be integers, $N:=(d+1)(r-1)$, $p$ a prime, $G\leq\sym_r$ an elementary abelian $p$-group with at least two orbits, and let $a\colon\Delta_n\longrightarrow\R^d$ be an affine map.
	Assume that $\LL\subseteq\K$ is a $G$-invariant sub-complex such that 
	\begin{compactenum}[\rm\quad  (A)]
		\item $\LL$ is homologically $(N-1)$-connected over $\FF_p$,
		\item for every $G$-box $U\neq V$ the complex $\LL[U]$ is connected and homologically $(N-3)$-connected over $\FF_p$, and
		\item  $\LL[U]$ is homologically $(N-2)$-connected over $\FF_p$ for every $G$-box $U$ whose stratum $N_G(U)$ is non-empty and closed.
	\end{compactenum}
	Then $\LL$ contains a Tverberg face of $a$.
\end{theorem}
\begin{proof}
	All cohomology in this proof has coefficients in $\FF_p$.
	The incidence space $T_{a,\LL}$ is a compact semialgebraic set, hence $T_{a,\LL}$ and $T'$ are triangulable~\cite[Thm.\,9.2.1]{BochnakCosteRoy1998}.
	For $\mathrm{E}G$ we can take a CW complex with a free cellular $G$-action, so that both $\mathrm{B}G$ and $\mathrm{E}G\times T'$ are CW complexes.
	The quotient map $\mathrm{E}G\times T'\longrightarrow T'_G$ by the free action of the finite group $G$ is a covering map, and it is a closed map with finite fibres.
	Hence, $T'_G$ is locally contractible, and it is paracompact as the image of a paracompact Hausdorff space under a closed continuous map~\cite{Michael1957}.
	So $T'$, its Borel construction $T'_G$ and the base $\mathrm{B}G\times S(\E^G)$ are paracompact, Hausdorff and locally contractible.
	 On such spaces singular, \v{C}ech, Alexander--Spanier and sheaf cohomology with constant coefficients coincide~\cite[Ch.\,III]{Bredon1997}.
	We may therefore use any of them while the sheaf-theoretic description is needed only for the Leray spectral sequence below.
	
	\medskip
	Suppose that $\LL$ contains no Tverberg face of $a$, and let $T'\subseteq T_{a,\LL}$ and $\rho\colon T'\longrightarrow S(\E^G)$ be as in Step~4 of the outline in Section~\ref{sec : deleted join}.
	Steps~1--4 of that proof use only the $(N-1)$-connectivity of $\LL$, and they give classes $e\in H^*(\mathrm{B}G)$ and $\omega\in H^{N-1}(\mathrm{B}G)$ with $\rho_G^*(e\times[S(\E^G)]^*)=\rho_G^*(\omega\times1)$ in $H^{N-1}_G(T')$, while by Step~6 we have $e\times[S(\E^G)]^*\neq\omega\times1$ in $H^{N-1}(\mathrm{B}G\times S(\E^G))$.
	\ul{So it suffices to show that $\rho_G^*$ is injective in degree $N-1$.} 
	This is Step~5 of the outline, which uses the Vietoris--Begle theorem, and here we use the Leray spectral sequence instead.

	\medskip 
	Since $G$ acts trivially on $S(\E^G)$, the Borel construction of $S(\E^G)$ is $\mathrm{B}G\times S(\E^G)$, and $\rho$ induces the map $\rho_G\colon T'_G\longrightarrow\mathrm{B}G\times S(\E^G)$, $[\xi,z]\longmapsto([\xi],\rho(z))$.
	The map $\rho_G$ is proper, that is, it is closed and all its fibres are compact.
	Here we use the notion of a proper map from~\cite[Ch.\,I, \S10.1, Def.\,1]{Bourbaki1989}: a continuous map $f\colon X\longrightarrow Y$ is proper if $f\times\mathrm{id}_Z\colon X\times Z\longrightarrow Y\times Z$ is closed for every space $Z$.
	By~\cite[Ch.\,I, \S10.2, Thm.\,1]{Bourbaki1989} this holds if and only if $f$ is closed and all fibres $f^{-1}(y)$ are compact.
	
	\medskip
	{\em $\rho_G$ is closed.}
	The space $T'=\pr_1^{-1}(S(\E^G))$ is a closed subspace of the compact space $T_{a,\LL}$, so it is compact.
	Hence $\rho$ is a continuous map from a compact space to the Hausdorff space $S(\E^G)$, and so it is proper~\cite[Ch.\,I, \S10.2, Cor.\,2]{Bourbaki1989}.
	By the definition of a proper map we have that the map $\mathrm{id}\times\rho\colon\mathrm{E}G\times T'\longrightarrow\mathrm{E}G\times S(\E^G)$ is closed.
	Denote by $\eta_{T'}\colon\mathrm{E}G\times T'\longrightarrow T'_G$ and $\eta_{S(\E^G)}\colon\mathrm{E}G\times S(\E^G)\longrightarrow\mathrm{B}G\times S(\E^G)$ the quotient maps.
	Then the following diagram commutes:
	\[
	\begin{tikzcd}[column sep=large]
	\mathrm{E}G\times T' \arrow[r,"\mathrm{id}\times\rho"] \arrow[d,"\eta_{T'}"'] & \mathrm{E}G\times S(\E^G) \arrow[d,"\eta_{S(\E^G)}"]\\
	T'_G \arrow[r,"\rho_G"'] & \mathrm{B}G\times S(\E^G).
	\end{tikzcd}
	\]
	Both quotient maps are closed, because for a closed set $C$ its saturation $\eta_{T'}^{-1}(\eta_{T'}(C))=\bigcup_{g\in G}gC$ is a finite union of closed sets, and similarly for $\eta_{S(\E^G)}$.
	If $C\subseteq T'_G$ is closed, then $\eta_{T'}^{-1}(C)$ is closed, and since $\eta_{T'}$ is surjective, the commutativity of the diagram gives that
	\[
	\rho_G(C)=\eta_{S(\E^G)}\big((\mathrm{id}\times\rho)(\eta_{T'}^{-1}(C))\big)
	\]
	is closed.
	Thus, $\rho_G$ is indeed closed.
	
	\medskip
	Moreover, a finite group acting continuously on a Hausdorff space acts properly, and the orbit space of a proper action is Hausdorff~\cite[Ch.\,III, \S4.1 and \S4.2]{Bourbaki1989}.
	Applied to the Hausdorff space $\mathrm{E}G\times T'$ this shows that $T'_G$ is Hausdorff, as claimed above.
	
	\medskip
	{\em The fibres of $\rho_G$ are compact.}
	Let $b=([\xi],y)\in \mathrm{B}G\times S(\E^G)$ be a point of the base and consider the continuous map
	\[
	j_\xi\colon\rho^{-1}(y)\longrightarrow\rho_G^{-1}(b),\qquad z\longmapsto[\xi,z].
	\]
	
	\smallskip
	It is surjective.
	Indeed, if $[\xi',z]\in\rho_G^{-1}(b)$, then $\xi'=g\xi$ for some $g\in G$, so $[\xi',z]=[\xi,g^{-1}z]$, and $\rho(g^{-1}z)=g^{-1}\rho(z)=g^{-1}y=y$ because $\rho$ is equivariant and $G$ fixes $S(\E^G)$ pointwise.
	
	\smallskip
	It is injective, since $[\xi,z]=[\xi,z']$ means $g\xi=\xi$ and $gz=z'$ for some $g\in G$. 
	Then $g$ is the neutral element because $G$ acts freely on $\mathrm{E}G$.
	The domain $\rho^{-1}(y)=\{y\}\times\LL[U_{a,y}]$ is a closed subset of the compact space $T'$, hence compact.
	Therefore the fibre $\rho_G^{-1}(b)=j_\xi(\rho^{-1}(y))$ is compact, as a continuous image of a compact space.
	Since $T'_G$ is Hausdorff, the continuous bijection $j_\xi$ from a compact space is a homeomorphism.
	So the fibre $\rho_G^{-1}(b)$ is homeomorphic to the finite simplicial complex $\LL[U_{a,y}]$, compare the proof of~\cite[Lem.\,3.5]{BlagojevicSierksma}.
	
	\medbreak
	Together with the closedness of $\rho_G$ this shows that $\rho_G$ is proper according to~\cite[Ch.\,I, \S10.2, Thm.\,1]{Bourbaki1989}.
	
	\medskip
	The Leray spectral sequence of the map $\rho_G\colon T'_G\longrightarrow\mathrm{B}G\times S(\E^G)$ has the form
	\[
		E_2^{s,t}=H^s\big(\mathrm{B}G\times S(\E^G);\mathscr R^t\big)\ \Longrightarrow\ H^{s+t}_G(T').
	\]
	Here $\mathscr R^t$ is the Leray sheaf of $\rho_G$, that is, the sheaf generated by the presheaf $W\longmapsto H^t(\rho_G^{-1}(W))$ on the base $\mathrm{B}G\times S(\E^G)$.
	For the spectral sequence see~\cite[Ch.\,IV, Sec.\,6]{Bredon1997}.
	
	\medskip
	In the following we are going to use two special properties of the sheaf $\mathscr R^t$.
	
	\medskip
	\emph{The stalk of $\mathscr R^t$ at a point $b=([\xi],y)$ of the base is the cohomology $H^t(\LL[U_{a,y}])$ of the fibre $\rho_G^{-1}(b)$}.
	By definition this stalk is the direct limit
	\begin{equation}\label{eq : stalk}
	\mathscr R^t_b=\varinjlim_{W\ni b}H^t\big(\rho_G^{-1}(W)\big)
	\end{equation}
	over the open neighbourhoods $W$ of $b$ in the base.
	Since $\rho_G$ is closed, the sets $\rho_G^{-1}(W)$ are cofinal among all open neighbourhoods of the fibre.
	Indeed, if $O\supseteq\rho_G^{-1}(b)$ is open, then $W:=(\mathrm{B}G\times S(\E^G))-\rho_G(T'_G-O)$ is an open neighbourhood of $b$ with $\rho_G^{-1}(W)\subseteq O$.
	A cofinal subsystem has the same direct limit, and so
	\begin{equation}\label{eq : cofinal}
	\varinjlim_{W\ni b}H^t\big(\rho_G^{-1}(W)\big)=\varinjlim_{O\supseteq\rho_G^{-1}(b)}H^t(O),
	\end{equation}
	where $O$ runs over the open neighbourhoods of $\rho_G^{-1}(b)$ in $T'_G$.
	
	\smallskip
	Recall that a subspace $A$ of a space $X$ is {\em taut} if the restriction maps induce an isomorphism
	\[
	\varinjlim_{O\supseteq A}H^t(O)\xrightarrow{\ \cong\ }H^t(A)
	\]
	for all $t$, where $O$ runs over the open neighbourhoods of $A$ in $X$, see~\cite[Ch.\,II, Sec.\,10]{Bredon1997} and~\cite[Sec.\,6.6]{Spanier1966}, and compare~\cite[Prop.\,2.5.1]{KashiwaraSchapira1990}.
	In other words, every class on $A$ extends to some neighbourhood, and two classes on neighbourhoods which agree on $A$ agree on a smaller neighbourhood.
	
	\smallskip
	Now, the fibre $\rho_G^{-1}(b)$ is compact, as we have just shown, and the space $T'_G$ is Hausdorff.
	Compact subspaces of Hausdorff spaces are taut~\cite[Ch.\,II, Sec.\,10]{Bredon1997}, \cite[Prop.\,2.5.1]{KashiwaraSchapira1990}, so the fibre $\rho_G^{-1}(b)$ is taut in $T'_G$.
	Alternatively, the fibre is closed in the paracompact space $T'_G$, and closed subspaces of paracompact spaces are taut~\cite[Ch.\,II, Sec.\,10]{Bredon1997}, \cite[Sec.\,6.6]{Spanier1966}.
	Combining \eqref{eq : stalk}, the cofinality \eqref{eq : cofinal}, the tautness of the fibre and the homeomorphism $j_\xi$ we get
	\[
	\mathscr R^t_b
	=\varinjlim_{W\ni b}H^t\big(\rho_G^{-1}(W)\big)
	=\varinjlim_{O\supseteq\rho_G^{-1}(b)}H^t(O)
	\cong H^t\big(\rho_G^{-1}(b)\big)
	\cong H^t\big(\LL[U_{a,y}]\big).
	\]
	This is the description of the stalks of the Leray sheaf of a closed map with taut fibres given in~\cite[Ch.\,IV, Sec.\,4]{Bredon1997}.
	
	\medskip
	\emph{For a subspace $Z$ of the base, the restriction of $\mathscr R^t$ to $Z$ is the Leray sheaf of the restriction $\rho_G^{-1}(Z)\longrightarrow Z$ of $\rho_G$}.
	Indeed, this restriction is again a closed map, since for a closed set $C\subseteq\rho_G^{-1}(Z)$ we have $C=C'\cap\rho_G^{-1}(Z)$ with $C'$ closed in $T'_G$ and hence $\rho_G(C)=\rho_G(C')\cap Z$ is closed in $Z$.
	It has the same fibres, which are compact, and so by the first property both sheaves have the stalk $H^t(\rho_G^{-1}(z))$ at every $z\in Z$, and the natural map from the first sheaf to the second is an isomorphism on all stalks.
	
	\medskip
	We now analyse the spectral sequence in more detail.
	By assumption (B) the fibres $\LL[U_{a,y}]$ of $\rho_G$ with $U_{a,y}\neq V$ are non-empty and connected.
	For $U_{a,y}=V$ the fibre is $\LL$, which is non-empty and connected by assumption (A), because $N\geq2$.
	Therefore, $\mathscr R^0$ is the constant sheaf $\FF_p$, $E_2^{N-1,0}=H^{N-1}(\mathrm{B}G\times S(\E^G))$, and so the edge homomorphism $E_2^{N-1,0}\longrightarrow H^{N-1}_G(T')$ coincides with $\rho_G^*$.

	\smallskip
	Let us recall how the edge homomorphism factors.
	The group $E_\alpha^{N-1,0}$ lies in the row $t=0$, so every differential which starts from it ends in the row $t=1-\alpha<0$ and is therefore zero.
	Hence, each $E_{\alpha+1}^{N-1,0}$ is the quotient of $E_\alpha^{N-1,0}$ by the image of the only differential which ends there,
	\[
	d_\alpha\colon E_\alpha^{N-1-\alpha,\alpha-1}\longrightarrow E_\alpha^{N-1,0}.
	\]
	For $\alpha\geq N$ the source of $d_\alpha$ lies in the column $s=N-1-\alpha<0$, so it is zero, and so $E_N^{N-1,0}=E_\infty^{N-1,0}$.
	Finally, $E_\infty^{N-1,0}$ is the last step of the filtration of $H^{N-1}_G(T')$, and so it is a subgroup of $H^{N-1}_G(T')$.
	Consequently, the edge homomorphism is the composition
	\[
	E_2^{N-1,0}\twoheadrightarrow E_3^{N-1,0}\twoheadrightarrow\cdots\twoheadrightarrow E_N^{N-1,0}=E_\infty^{N-1,0}\hookrightarrow H^{N-1}_G(T')
	\]
	of the quotient maps by the images of $d_2,\dots,d_{N-1}$ followed by an inclusion.
	It is injective if and only if each of these quotient maps is an isomorphism, that is, if and only if $d_\alpha=0$ for all $2\leq\alpha\leq N-1$.
	This certainly holds when all sources of these differentials vanish.
	The source $E_\alpha^{N-1-\alpha,\alpha-1}$ of $d_\alpha$ is a subquotient of $E_2^{N-1-\alpha,\alpha-1}$, and with $t:=\alpha-1$ these are the groups $E_2^{N-2-t,t}$ with $1\leq t\leq N-2$.
	\ul{Thus, it suffices to show that $E_2^{N-2-t,t}=0$ for $1\leq t\leq N-2$.}
	
	\smallskip
	For $1\leq t\leq N-3$ all stalks of $\mathscr R^t$ vanish by the assumptions (A) and (B), so $E_2^{N-2-t,t}=0$.
	If $N=2$ there are no such differentials.
	On the other hand, if $N\geq3$, it remains to show that $E_2^{0,N-2}=H^0(\mathrm{B}G\times S(\E^G);\mathscr R^{N-2})$, the group of global sections of $\mathscr R^{N-2}$, is zero.
	For an illustration of the $E_2$-term of this spectral sequence see Figure~\ref{fig : E2}.

	\begin{figure}[ht]
		\centering
		\begin{tikzpicture}[x=13mm,y=9mm,lab/.style={font=\scriptsize},e/.style={font=\small,fill=white,inner sep=1.5pt}]
			\draw[->] (-0.5,-0.5)--(6.3,-0.5) node[lab,right] {$s$};
			\draw[->] (-0.5,-0.5)--(-0.5,5.0) node[lab,above] {$t$};
			\foreach \x/\l in {0/0,1/1,2/2,4/N-2,5/N-1}{\node[lab] at (\x,-0.85) {$\l$};}
			\node[lab] at (3,-0.85) {$\cdots$};
			\foreach \y/\l in {0/0,1/1,3/N-3,4/N-2}{\node[lab,anchor=east] at (-0.6,\y) {$\l$};}
			\node[lab,anchor=east] at (-0.6,2) {$\vdots$};
			\draw[dashed,gray] (-0.45,4.6)--(6.0,4.6);
			\draw[->,thick,gray] (0.25,3.8)--(4.72,0.2);
			\node[lab,fill=white,inner sep=1pt] at (2.55,2.3) {$d_{N-1}$};
			\fill[gray!20] (-0.3,3.7) rectangle (0.3,4.3);
			\node[e,fill=gray!20] at (0,4) {$\ast$};
			\foreach \x in {0,1,2,4,5}{\node[e] at (\x,0) {$\ast$};}
			\foreach \x in {0,1,2,4,5}{\node[e] at (\x,1) {$0$}; \node[e] at (\x,3) {$0$};}
			\foreach \x in {0,1,2,4,5}{\node[e] at (\x,2) {$\vdots$};}
			\foreach \y in {0,1,3,4}{\node[e] at (3,\y) {$\cdots$};}
			\foreach \x in {1,2,4,5}{\node[e] at (\x,4) {$\ast$};}
			\draw[thick] (5,0) circle (0.25);
		\end{tikzpicture}
		\caption{\small The $E_2$-term of the Leray spectral sequence of $\rho_G$ in the rows $0\leq t\leq N-2$.
			The rows $1\leq t\leq N-3$ vanish, since all stalks of $\mathscr R^t$ vanish there.
			The differential $d_\alpha$ which hits the circled group $E_\alpha^{N-1,0}$ starts in the row $t=\alpha-1$, which vanishes for $2\leq\alpha\leq N-2$. Only the source $E_2^{0,N-2}$ of $d_{N-1}$, shaded, is not known to vanish a priori, and it is shown to vanish in the proof.}
		\label{fig : E2}
	\end{figure}
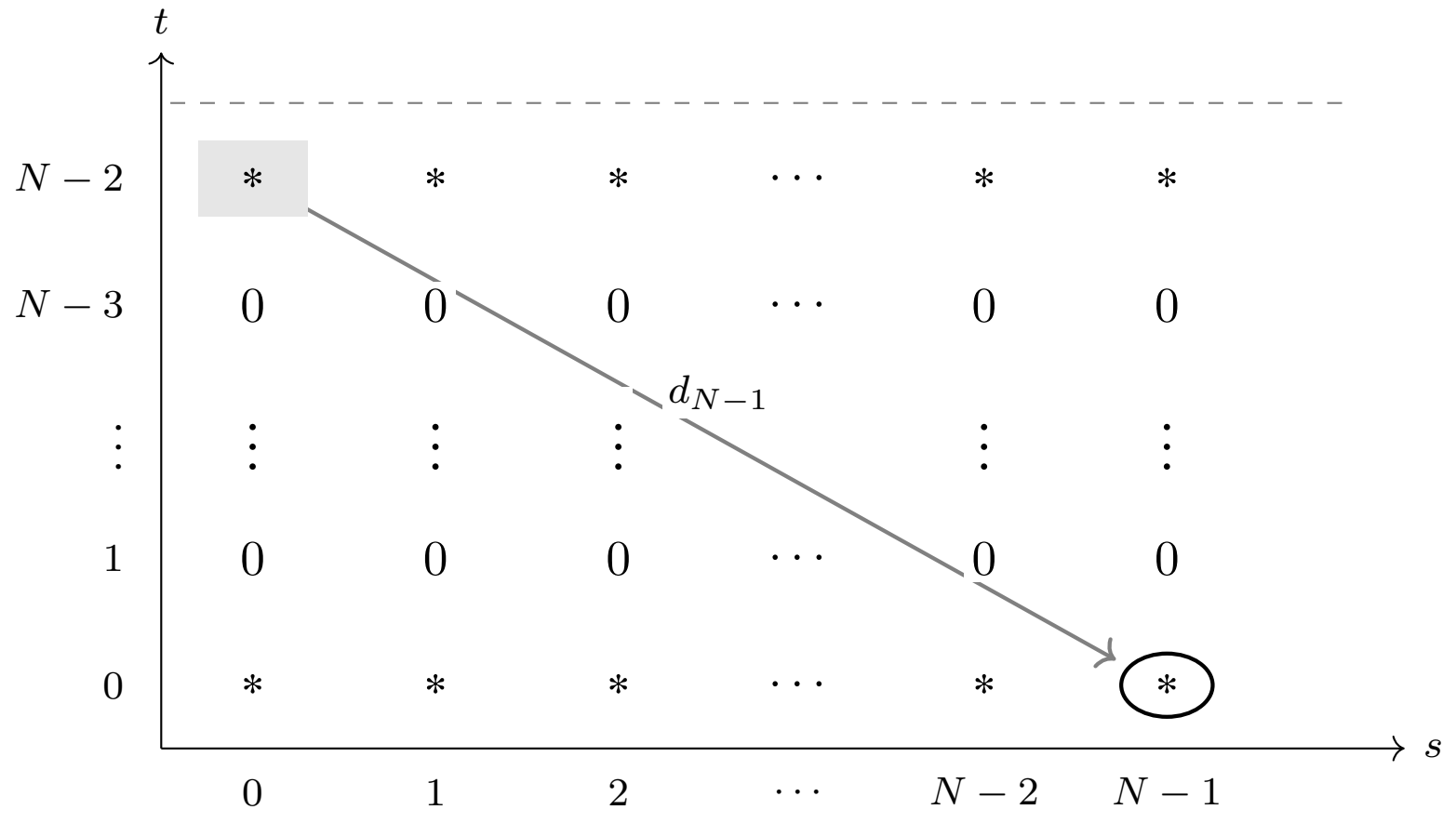
	
	\medskip
	For that we use two properties of strata.
	First, a non-empty stratum of a $G$-box $U\neq V$ is connected: it is the intersection of $S(\E^G)$ with the closed half-spaces $\{y\in\E : \langle y,A(u)\rangle\leq0\}$, $u\in U$, and the open half-spaces $\{y\in\E :\langle y,A(v)\rangle>0\}$, $v\in V-U$. Since $U\neq V$, at least one open half-space occurs, so the stratum contains no pair of antipodal points and with any two of its points the shortest arc between them.
	Second, if $y''$ lies in the closure of $N_G(U)$ in $S(\E^G)$ but not in $N_G(U)$, then $U\subsetneq U_{a,y''}$.
	Indeed, $y''$ is the limit of a sequence $y_1,y_2,\ldots$ in $N_G(U)$, and for every $u\in U$ we have $\langle y_m,A(u)\rangle\leq0$ for all $m$, hence $\langle y'',A(u)\rangle\leq0$, that is $u\in U_{a,y''}$.
	So $U\subseteq U_{a,y''}$, and the inclusion is strict because $U_{a,y''}\neq U$ when $y''\notin N_G(U)$.
	Thus $U_{a,y''}$ contains at least one vertex $v\in V-U$.
	For such a vertex the strict inequalities $\langle y_m,A(v)\rangle>0$ turn into the equality $\langle y'',A(v)\rangle=0$ in the limit.
	
	\medskip
	Now let $\varsigma$ be a global section of $\mathscr R^{N-2}$, and let $\Sigma$ be its support, that is the set of points of the base at which the germ of $\varsigma$ is not zero.
	The support of a section is always closed, since a germ which is zero at a point is zero in a neighbourhood of this point.
	We show that $\Sigma=\emptyset$, which means $\varsigma=0$.

	\medskip
	Suppose that $\Sigma\neq\emptyset$.
	For every point $([\xi'],y')\in\Sigma$ the set $U_{a,y'}$ is a $G$-box, since $y'\in S(\E^G)$.
	These $G$-boxes form a non-empty family of subsets of the finite set $V$, and a non-empty finite family of sets, ordered by inclusion, has a maximal element.
	So we may choose a point $([\xi],y)\in\Sigma$ such that no $G$-box $U_{a,y'}$ with $([\xi'],y')\in\Sigma$ strictly contains $U:=U_{a,y}$.
	By the first property of $\mathscr R^{N-2}$, its stalk at $([\xi],y)$ is $H^{N-2}(\LL[U])$.
	This stalk is not zero, since it contains the germ of $\varsigma$ at $([\xi],y)$, which is not zero because $([\xi],y)\in\Sigma$.
	Hence, 
	\begin{compactitem}[\quad --]
		\item $U\neq V$, because $\LL$ is homologically $(N-1)$-connected, and
		\item the stratum $N_G(U)$, which contains $y$, is not closed, by assumption (C).
	\end{compactitem}

	\medskip
	Next we restrict $\varsigma$ to $Y:=\{[\xi]\}\times N_G(U)$.
	The preimage $\rho_G^{-1}(Y)$ consists of the classes $[\xi,z]$ with $z\in T'$ and $\rho(z)\in N_G(U)$.
	Since $G$ acts freely on $\mathrm{E}G$, every such class has a unique representative with first coordinate $\xi$.
	Moreover, the map $T'\longrightarrow T'_G$, $z\longmapsto[\xi,z]$, is a continuous injection of a compact space into a Hausdorff space, hence an embedding.
	So $\rho_G^{-1}(Y)$ is homeomorphic to 
	\[
	\rho^{-1}(N_G(U))=\{(y',\lambda) : y'\in N_G(U),\ \supp(\lambda)\subseteq U\}=N_G(U)\times\LL[U],
	\] 
	and over $Y$ the map $\rho_G$ is the projection of this product onto $N_G(U)$.
	By the second property of $\mathscr R^{N-2}$, its restriction to $Y$ is the Leray sheaf of a product projection, that is, the constant sheaf with stalk $H^{N-2}(\LL[U])$.
	The sections of a constant sheaf with stalk $M$ over a space are the continuous maps from this space to $M$ with the discrete topology, that is, the locally constant maps~\cite[Ch.\,I, Sec.\,1]{Bredon1997}, \cite[Ch.\,II, Ex.\,1.0.3]{Hartshorne1977}.
	In particular, a section of a constant sheaf over a connected space is constant, and $Y$ is connected by the first property of strata.
	So the germs of $\varsigma$ at all points of $Y$ coincide with its germ at $([\xi],y)$, and they are not zero.
	Consequently, $Y\subseteq\Sigma$.

	\medskip
	Since $\Sigma$ is closed, it also contains the closure $\{[\xi]\}\times\overline{N_G(U)}$ of $Y$.
	As $N_G(U)$ is not closed, there is a point $y''\in\overline{N_G(U)}-N_G(U)$, and $U\subsetneq U_{a,y''}$ by the second property of strata.
	Then $([\xi],y'')\in\Sigma$ and $U_{a,y''}$ is strictly larger than $U$, which contradicts the maximality of $U$.
	Hence, $\varsigma=0$, so $E_2^{0,N-2}=0$, and $\rho_G^*$ is injective in degree $N-1$.
	This concludes the proof.
\end{proof}

\medskip
Note that for groups with at least two orbits, Theorem~\ref{th : G-box} is the special case of Theorem~\ref{th : refined} in which nothing is known about the strata.

\medskip
\subsection{Closed strata}
\label{sec : strata}

In this subsection we describe the strata on the fixed sphere with respect to the action of elementary abelian groups.
Specifically, we show that when a stratum 
\[
	N_G(U)\ =\ \{y\in S(\E^G) : U_{a,y}=U\}
\]
is closed, $U$ is forced to have many orbits in each column.

\medskip
Let $a\colon\Delta_n\longrightarrow\R^d$ be an affine map, $q=p^\alpha$ a prime power dividing $r$, $k:=\frac{r}{q}$, and let $G=\FF_p^{\,\alpha}$ act on $[r]=\FF_p^{\,\alpha}\times[k]$ as in Lemma~\ref{lem : group}\,(A).
Then the $G$-orbits of this action on $[r]$ are the sets $\FF_p^{\,\alpha}\times\{\iota\}$, indexed by $\iota\in[k]$.

\medskip
The strata can be described  by affine functions in the following way, compare with~\cite[Sec.\,2 and Lem.\,3.1]{BlagojevicSierksma}.
For $y\in\E$ denote by $y_1,\dots,y_r\in\R^{d+1}$ the columns of the matrix $y$.
Since the rows of $y$ sum to zero, we have $y_1+\cdots+y_r=0$.
For $1\leq i\leq r$ let $f^y_i\colon \R^d\longrightarrow\R$ be the affine function on $\R^d$ given by
\[
	f^y_i(x)\ :=\ \langle y_i,(1,x)\rangle , 
\]
for $x\in\R^d$.
Thus, $f^y_i(a(p_j))=\langle y_i,\widehat{a}(p_j)\rangle$ for $0\leq j\leq n$, where $\widehat{a}(x)=(1,a(x))$ is the map used in the definition of the test map $A$ in the proof outline in Section~\ref{sec : deleted join}.
Recall from there that the $\ell$-th column of $A(p_{j,i})\in\E$ is $(\delta_{i\ell}-\frac1r)\,\widehat{a}(p_j)$, where $\delta_{i\ell}$ denotes the Kronecker delta. 
Therefore, for all $0\leq j\leq n$ and $1\leq i\leq r$ we have
\begin{equation}
	\label{eq : values}
	\langle y,A(p_{j,i})\rangle\ =\ \sum_{\ell=1}^r\big(\delta_{i\ell}-\tfrac1r\big)\langle y_\ell,\widehat{a}(p_j)\rangle\ =\ \langle y_i,\widehat{a}(p_j)\rangle-\tfrac1r\Big\langle\sum_{\ell=1}^ry_\ell,\widehat{a}(p_j)\Big\rangle\ =\ f^y_i(a(p_j)).
\end{equation}
Thus, the vertex $p_{j,i}$ lies in $U_{a,y}$ if and only if $f^y_i(a(p_j))\leq0$.
Moreover, $f^y_1+\cdots+f^y_r=0$, because $y_1+\cdots+y_r=0$.
The affine function $f^y_i$ corresponds to $y_i$ under the usual linear isomorphism $z\longmapsto\langle z,(1,\cdot\,)\rangle$ between $\R^{d+1}$ and the space of affine functions on $\R^d$.
Hence, the map $y\longmapsto(f^y_1,\dots,f^y_r)$ is a linear isomorphism from $\E$ onto the space of $r$-tuples of affine functions on $\R^d$ with sum zero.

\medskip
An element $g\in G$ moves the column $i$ of $y$ to the column $g\cdot i$, so $f^{g\cdot y}_{g\cdot i}=f^y_i$.
Hence, $y\in\E^G$ if and only if $f^y_i$ depends only on the orbit $\iota$ of $i$, and then we write $\phi^y_\iota$ for this function.
Every orbit has $q$ elements, so $0=f^y_1+\cdots+f^y_r=q\,(\phi^y_1+\cdots+\phi^y_k)$, and $\phi^y_1+\cdots+\phi^y_k=0$.
Conversely, a $k$-tuple of affine functions with sum zero defines, by repeating each function $q$ times, an $r$-tuple with sum zero, and hence a unique point of $\E^G$.
Therefore, $y\longmapsto(\phi^y_1,\dots,\phi^y_k)$ is a linear isomorphism from $\E^G$ onto the space of $k$-tuples of affine functions on $\R^d$ with sum zero.
In particular, $\dim\E^G=(d+1)(k-1)$.

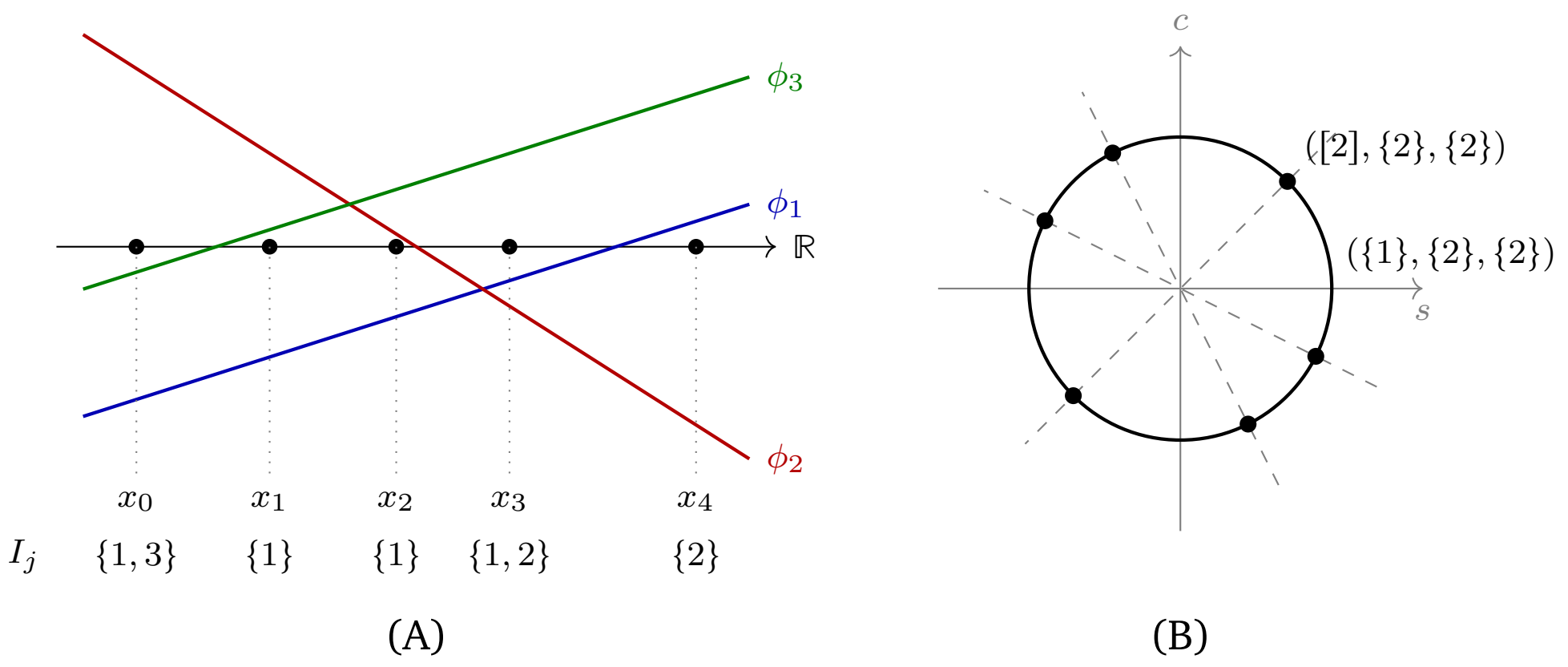
\begin{figure}[ht]
	\centering
	\begin{tikzpicture}[x=11mm,y=3.5mm,lab/.style={font=\scriptsize},pt/.style={circle,fill=black,inner sep=1.3pt}]
		\draw[->] (-2.2,0)--(3.2,0) node[lab,right] {$\R$};
		\draw[thick,blue!70!black] (-2,-4)--(3,1) node[lab,right] {$\phi_1$};
		\draw[thick,red!70!black] (-2,5)--(3,-5) node[lab,right] {$\phi_2$};
		\draw[thick,green!50!black] (-2,-1)--(3,4) node[lab,right] {$\phi_3$};
		\foreach \x/\j/\I in {-1.6/0/{\{1,3\}},-0.6/1/{\{1\}},0.35/2/{\{1\}},1.2/3/{\{1,2\}},2.6/4/{\{2\}}}{
		\node[pt] at (\x,0) {};
		\draw[gray,dotted] (\x,0)--(\x,-5.4);
		\node[lab] at (\x,-6.0) {$x_{\j}$};
		\node[lab] at (\x,-7.3) {$\I$};}
		\node[lab,anchor=east] at (-2.2,-7.3) {$I_j$};
		\node[font=\small] at (0.5,-9.2) {(A)};
	\end{tikzpicture}
	\qquad
	\begin{tikzpicture}[scale=1.25,lab/.style={font=\scriptsize},pt/.style={circle,fill=black,inner sep=1.4pt}]
		\draw[->,gray] (-1.6,0)--(1.6,0) node[lab,below] {$s$};
		\draw[->,gray] (0,-1.6)--(0,1.6) node[lab,above] {$c$};
		\foreach \a in {45,-26.565,-63.435}{\draw[gray,dashed] (\a:1.45)--(\a+180:1.45);}
		\draw[thick] (0,0) circle (1);
		\foreach \a in {45,-26.565,-63.435,225,153.435,116.565}{\node[pt] at (\a:1) {};}
		\node[lab,anchor=south west] at (45:1) {$([2],\{2\},\{2\})$};
		\node[lab,anchor=west,fill=white,inner sep=1pt] at (12:1.08) {$(\{1\},\{2\},\{2\})$};
		\node[font=\small] at (0,-2.3) {(B)};
	\end{tikzpicture}
	\caption{\small The description of strata by affine functions for $d=1$.
		(A) The case $k=3$: a triple $(\phi_1,\phi_2,\phi_3)$ of affine functions on $\R$ with sum zero, and five points $x_0,\dots,x_4$. Below every point $x_j$ is the set $I_j=\{\iota : \phi_\iota(x_j)\leq0\}$ of orbits in the column $Q_j$ of the corresponding $G$-box.
		The orbits $2$ and $3$ are missing in the columns of $x_1$ and $x_2$, which consist of a single orbit.
		(B) The case $k=2$: for the points $x_0=-1$, $x_1=\frac12$, $x_2=2$ the fixed sphere is the circle of pairs $(\phi,-\phi)$ with $\phi(x)=sx+c$ and $s^2+c^2=1$. The dashed lines are the sets $\{\phi(x_j)=0\}$.
		Its strata are the six open arcs and the six marked points, labelled by the families $(I_0,I_1,I_2)$. Exactly the six points are closed strata, and each of them has a column $[2]$.}
	\label{fig : strata}
\end{figure}

\medskip
Now, let $x_j:=a(p_j)$ for $0\leq j\leq n$.  
For a $G$-box $U$, with $I_j\subseteq[k]$ denoting the set of orbits contained in $U\cap Q_j$, we have that a point $y\in S(\E^G)$ lies in $N_G(U)$ if and only if $I_j=\{\iota : \phi^y_\iota(x_j)\leq0\}$ for all $0\leq j\leq n$.
Figure~\ref{fig : strata} illustrates this description for $d=1$.

The closed strata have the following dual characterisation.

\medskip
\begin{lemma}
	\label{lem : closed strata}
	Let $r\geq2$, $d\geq1$ and $n\geq0$ be integers, $G\leq\sym_r$ an elementary abelian group, $a\colon\Delta_n\longrightarrow\R^d$ an affine map, and $U\neq V$ a $G$-box whose stratum is non-empty.
	Then $N_G(U)$ is closed if and only if every $y\in\E^G-\{0\}$ with $\langle y,A(u)\rangle\leq0$ for all $u\in U$ satisfies that $\langle y,A(v)\rangle>0$ for all $v\in V-U$, that is
	\begin{equation}\label{eq:condition}
		(\forall y\in \E^G-\{0\})\ (\forall u\in U)\ \langle y,A(u)\rangle\leq0 \ \Longrightarrow \ (\forall v\in V-U)\ \langle y,A(v)\rangle>0. 
	\end{equation}
\end{lemma}
\begin{proof}
	Set $C:=\{y\in\E^G : \langle y,A(u)\rangle\leq0\ \text{for all}\ u\in U\}$.

	\smallskip
	Assume that the condition \eqref{eq:condition} holds.
	Every $y\in N_G(U)$ satisfies $\langle y,A(u)\rangle\leq0$ for all $u\in U$, so $N_G(U)\subseteq C\cap S(\E^G)$.
	Conversely, a point $y\in C\cap S(\E^G)$ is non-zero, and \eqref{eq:condition} gives $\langle y,A(v)\rangle>0$ for all $v\in V-U$, that is, $U_{a,y}=U$ and $y\in N_G(U)$.
	Hence $N_G(U)=C\cap S(\E^G)$, which is closed as the intersection of the sphere $S(\E^G)$ with finitely many closed half-spaces.

	\smallskip
	Conversely, let $N_G(U)$ be closed, $y_0\in N_G(U)$ and $y_1\in C-\{0\}$.
	The segment from $y_0$ to $y_1$ passes through $0$ only if $y_1=-cy_0$ for some $c>0$, and we show that this is impossible.
	Suppose $y_1=-cy_0$ with $c>0$, and let $u\in U$.
	On the one hand, $\langle y_0,A(u)\rangle\leq0$, as $y_0\in N_G(U)\subseteq C$.
	On the other hand, $\langle y_0,A(u)\rangle=-\frac1c\langle y_1,A(u)\rangle\geq0$, as $y_1\in C$.
	So $\langle y_0,A(u)\rangle=0$ for all $u\in U$. 
	Since $U\neq V$, there is a column $Q_j\not\subseteq U$.
	The vectors $A(p_{j,1}),\dots,A(p_{j,r})$ sum to zero, because $\sum_{i=1}^r(\delta_{i\ell}-\frac1r)=0$ for every $\ell$, and therefore the values $\langle y_0,A(p_{j,i})\rangle$, $1\leq i\leq r$, sum to zero as well.
	Those with $p_{j,i}\in U$ vanish, and the others, of which there is at least one, are positive, because $y_0\in N_G(U)$. This is impossible.
	So the segment from $y_0$ to $y_1$ avoids $0$, and its normalisation is a path in $C\cap S(\E^G)$ from $y_0$ to $\frac{y_1}{\|y_1\|}$.
	Denote this path by $\gamma\colon[0,1]\longrightarrow C\cap S(\E^G)$, and let $J:=\{t\in[0,1] : \gamma(t)\in N_G(U)\}$.
	Since $\gamma(t)\in C$ for every $t$, a parameter $t$ lies in $J$ if and only if $\langle\gamma(t),A(v)\rangle>0$ for all $v\in V-U$.
	These are finitely many strict inequalities between continuous functions of $t$, so $J$ is open in $[0,1]$.
	On the other hand, $J=\gamma^{-1}(N_G(U))$ is closed in $[0,1]$, because $N_G(U)$ is closed.
	Moreover, $0\in J$, as $\gamma(0)=y_0\in N_G(U)$.
	Since the interval $[0,1]$ is connected, $J=[0,1]$, and in particular $\gamma(1)=\frac{y_1}{\|y_1\|}\in N_G(U)$, that is, $\langle y_1,A(v)\rangle>0$ for all $v\in V-U$.
	Thus, condition  \eqref{eq:condition} is satisfied.
\end{proof}

\medskip
The previous lemma is closely related to~\cite[Lem.\,2.21]{Blagojevic2025}.
For $y\in\E^G$ we have $\langle y,A(v)\rangle=\langle y,A^G(v)\rangle$, where $A^G$ denotes the composition of $A$ with the orthogonal projection $\E\longrightarrow\E^G$.
Hence, the stratum $N_G(U)$ is the set of points of $S(\E^G)$ which lie in the open hemispheres $\{y : \langle y,A^G(v)\rangle>0\}$ for all $v\in V-U$, and in none of the open hemispheres $\{y : \langle y,A^G(u)\rangle>0\}$ for $u\in U$.
In the notation of~\cite[Sec.\,2.2 and Sec.\,2.4]{Blagojevic2025} this is the set 
\[
N(V,U)=B(A^G(V-U))-B\langle A^G(U)\rangle
\]
for the map $A^G$, after the vectors $A^G(v)$ are normalised to unit length.
Since $S(\E^G)$ is compact, the stratum $N_G(U)$ is closed if and only if it is compact.
By~\cite[Lem.\,2.21]{Blagojevic2025}, the homology of the Alexander dual of the zero-avoiding complex induced on $U$ is concentrated in a single degree when $N(V,U)$ is compact, and it vanishes in a range of degrees otherwise.
Lemma~\ref{lem : closed strata} characterises the same dichotomy, compact versus non-compact, by the condition \eqref{eq:condition} on the linear functionals $\langle\,\cdot\,,A(v)\rangle$ on $\E^G$.
This parallels Theorem~\ref{th : refined}, where the strongest connectivity is required only for the $G$-boxes whose strata are closed, that is, compact.

\medskip
We also use a separation property of convex sets in the following way.

\medskip
\begin{lemma}
	\label{lem : separation}
	Let $d\geq1$ and $m\geq2$ be integers and let $C_1,\dots,C_m$ be non-empty convex subsets of $\R^d$ such that $\operatorname{relint}(C_1)\cap\cdots\cap\operatorname{relint}(C_m)=\emptyset$.
	Then there are affine functions $\ell_1,\dots,\ell_m$ on $\R^d$, not all zero, with $\ell_1+\cdots+\ell_m=0$ and $\ell_i(z)\leq0$ for all $z\in C_i$ and all $1\leq i\leq m$.
\end{lemma}
\begin{proof}
	The relative interior of $C:=C_1\times\cdots\times C_m\subseteq(\R^d)^m$ is $\operatorname{relint}(C_1)\times\cdots\times\operatorname{relint}(C_m)$, which is disjoint from the diagonal $D:=\{(z,\dots,z):z\in\R^d\}$.
	By~\cite[Thm.\,11.3]{Rockafellar1970} the two sets are properly separated.
	This means that there are linear functions $u_1,\dots,u_m$ on $\R^d$, not all zero, and a constant $c_0\in\R$ with $\sum_{i=1}^mu_i(z_i)\leq c_0\leq\sum_{i=1}^mu_i(z)$ for all $(z_1,\dots,z_m)\in C$ and $(z,\dots,z)\in D$.
	The right side is a linear function of $z$ which is bounded below, so it vanishes. 
	Hence, $u_1+\cdots+u_m=0$ and $c_0\leq0$.
	With $\mu_i:=\sup_{C_i}u_i$ we get $\mu_1+\cdots+\mu_m\leq c_0\leq0$.
	Then, $\ell_i:=u_i-\mu_i+\frac1m(\mu_1+\cdots+\mu_m)$, for $1\leq i\leq m$, are affine functions having the required properties.
\end{proof}

\begin{figure}[ht]
	\centering
	\begin{tikzpicture}[x=12mm,y=12mm,lab/.style={font=\small}]
		\fill[gray!25] (-0.18,-0.78)--(0.18,0.78)--(-1.2,1.1)--(-2,0)--(-1.3,-1.0)--cycle;
		\fill[gray!55] (-0.18,-0.78)--(1.6,-1.1)--(1.9,0.6)--(0.18,0.78)--cycle;
		\draw (-0.18,-0.78)--(0.18,0.78)--(-1.2,1.1)--(-2,0)--(-1.3,-1.0)--cycle;
		\draw (-0.18,-0.78)--(1.6,-1.1)--(1.9,0.6)--(0.18,0.78)--cycle;
		\draw[thick] (-0.42,-1.82)--(0.42,1.82);
		\node[lab] at (-1.05,0.1) {$C_1$};
		\node[lab] at (1.0,-0.15) {$C_2$};
		\node[lab,anchor=south] at (0.42,1.82) {$\{\ell_1=0\}$};
		\node[lab] at (-1.9,-1.55) {$\ell_1<0$};
		\node[lab] at (1.9,-1.55) {$\ell_1>0$};
		\node[font=\small] at (0,-2.6) {(A)};
	\end{tikzpicture}
	\qquad\qquad
	\begin{tikzpicture}[x=10mm,y=10mm,lab/.style={font=\small}]
		\foreach \a in {0,120,240}{\fill[gray,opacity=0.22,rotate=\a] (-1.732,-1)--(1.732,-1)--(1.1547,0)--(-1.1547,0)--cycle;}
		\draw (0,2)--(-1.732,-1)--(1.732,-1)--cycle;
		\foreach \a in {0,120,240}{\draw[thick,rotate=\a] (-1.55,0)--(1.55,0);}
		\fill (0,0) circle (1.4pt);
		\node[lab,anchor=south] at (0,2) {$v_1$};
		\node[lab,anchor=north east] at (-1.732,-1) {$v_2$};
		\node[lab,anchor=north west] at (1.732,-1) {$v_3$};
		\node[lab] at (0,-0.72) {$C_1$};
		\node[lab] at (0.62,0.36) {$C_2$};
		\node[lab] at (-0.62,0.36) {$C_3$};
		\node[font=\small] at (0,-2.0) {(B)};
	\end{tikzpicture}
	\caption{\small Lemma~\ref{lem : separation} in two cases.
		(A) For $d=2$ and $m=2$ the convex sets $C_1$ and $C_2$ meet along a segment, but their relative interiors are disjoint. The line $\{\ell_1=0\}$ separates them, and $\ell_2=-\ell_1$.
		(B) For $d=2$ and $m=3$ let $\lambda_i$ be the barycentric coordinate of the vertex $v_i$ of the triangle, and let $C_i$ be the part of the triangle where $\lambda_i\leq\frac13$.
		Any two of these sets have overlapping interiors, so no two of them can be separated, but only the barycentre lies in all three, and it is not in their interiors.
		The affine functions $\ell_i:=\lambda_i-\frac13$ have sum zero and $\ell_i\leq0$ on $C_i$, and the thick lines are the lines $\{\ell_i=0\}$.}
	\label{fig : separation}
\end{figure}
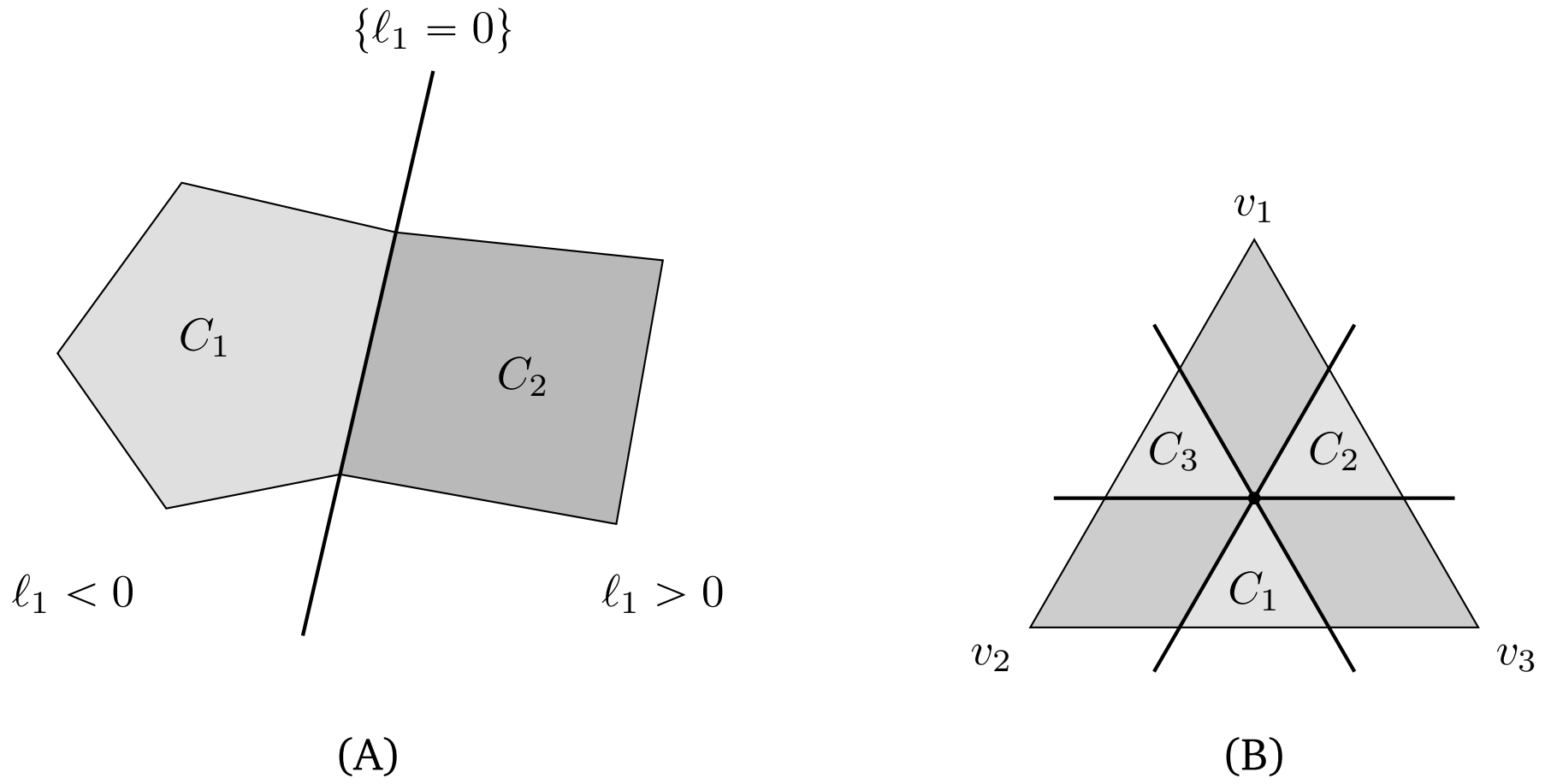

\medskip
For $m=2$, Lemma~\ref{lem : separation} is the usual proper separation of two convex sets by a hyperplane, with $\ell_2=-\ell_1$.
For $m\geq3$ the relative interiors of the sets may intersect pairwise, and then no two of them can be separated, see Figure~\ref{fig : separation}.

\medskip
\begin{theorem}
	\label{th : closed strata}
	Let $r\geq2$, $d\geq1$ and $n\geq d$ be integers, let $q=p^\alpha$ be a prime power dividing $r$, let $k:=\frac{r}{q}$, identify $[r]$ with $\FF_p^{\,\alpha}\times[k]$, and let $G=\FF_p^{\,\alpha}$ act on $[r]$ by translation of the first coordinate, so that the orbits of $G$ are the $k$ sets $\FF_p^{\,\alpha}\times\{\iota\}$, $\iota\in[k]$.
	
	\smallskip\noindent
	Let $a\colon\Delta_n\longrightarrow\R^d$ be an affine map such that the points $a(p_0),\dots,a(p_n)$ affinely span $\R^d$.
	If the stratum of a $G$-box $U$ is non-empty and closed, then every set $U\cap Q_j$ is the union of at least $k-d$ orbits.
\end{theorem}
\begin{proof}
	If $U=V$ there is nothing to prove, so let $U\neq V$.
	We identify $y\in\E^G$ with the tuple $(\phi^y_1,\dots,\phi^y_k)$. 
	As above, write $x_j:=a(p_j)$ and denote by
	\[
	I_j:=\{\iota\in[k] : \FF_p^{\,\alpha}\times\{\iota\}\subseteq U\cap Q_j\}
	\]
	the set of orbits contained in the column $U\cap Q_j$ of the $G$-box $U$, for $0\leq j\leq n$.
	Since $U$ is a $G$-box, each set $U\cap Q_j$ is the union of the orbits in $I_j$, and $I_j\neq\emptyset$.
	We say that the orbit $\iota$ is {\em missing} in the column $j$ if $\iota\notin I_j$, and we collect the orbits which are missing in at least one column in the set
	\[
	D:=\{\iota\in[k] : \iota\notin I_j\text{ for some }0\leq j\leq n\}=[k]-(I_0\cap\cdots\cap I_n).
	\]
	Since $U\neq V$, some set $I_j$ is a proper subset of $[k]$, and so $D\neq\emptyset$.
	For $\iota\in[k]$ we put 
	\[
	P_\iota:=\conv\{x_j : \iota\in I_j\},
	\]
	the convex hull of the points $x_j$ in whose columns the orbit $\iota$ is present. An example is shown in Figure~\ref{fig : P}.

	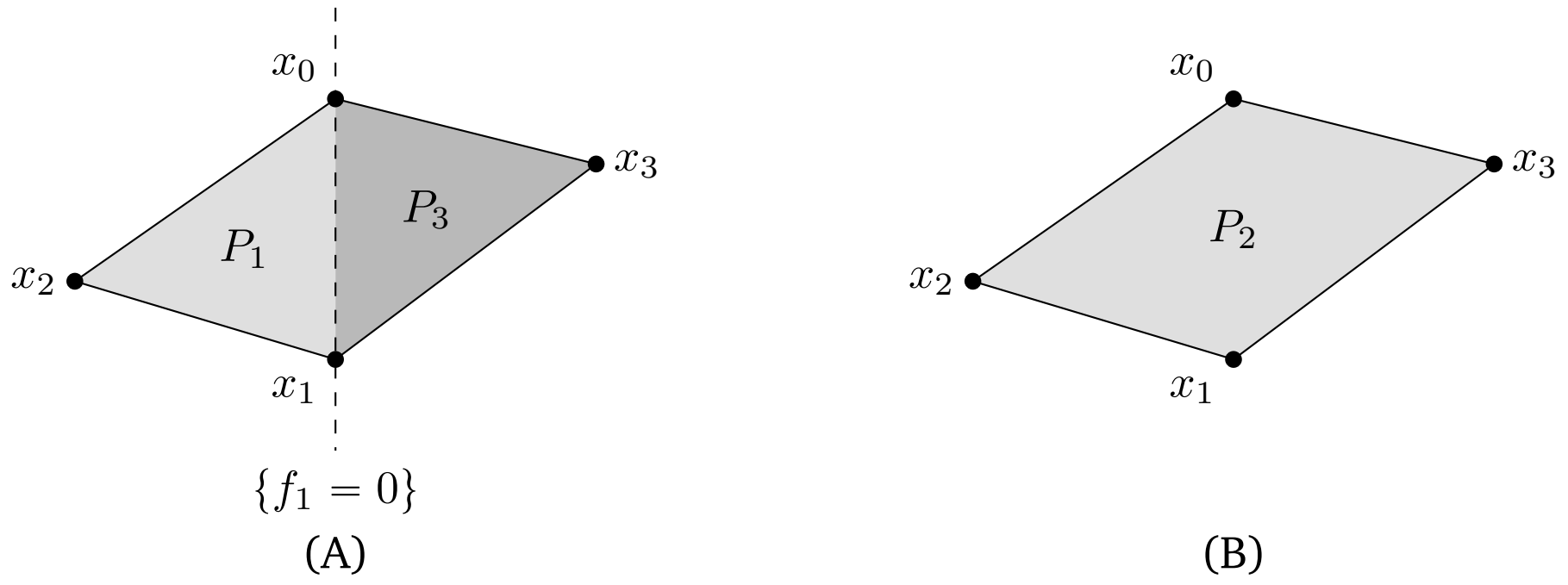
\begin{figure}[ht]
		\centering
		\begin{tikzpicture}[x=10mm,y=10mm,lab/.style={font=\small},pt/.style={circle,fill=black,inner sep=1.3pt}]
			\fill[gray!25] (0,1)--(0,-1)--(-2,-0.4)--cycle;
			\fill[gray!55] (0,1)--(0,-1)--(2,0.5)--cycle;
			\draw (0,1)--(-2,-0.4)--(0,-1)--(2,0.5)--cycle;
			\draw[dashed] (0,1.7)--(0,-1.7) node[lab,below] {$\{f_1=0\}$};
			\foreach \p/\n/\a in {(0,1)/0/north west,(0,-1)/1/south west,(-2,-0.4)/2/east,(2,0.5)/3/west}{\node[pt] at \p {};}
			\node[lab,anchor=south east] at (0,1) {$x_0$};
			\node[lab,anchor=north east] at (0,-1) {$x_1$};
			\node[lab,anchor=east] at (-2,-0.4) {$x_2$};
			\node[lab,anchor=west] at (2,0.5) {$x_3$};
			\node[lab] at (-0.7,-0.15) {$P_1$};
			\node[lab] at (0.7,0.15) {$P_3$};
			\node[font=\small] at (0,-2.5) {(A)};
		\end{tikzpicture}
		\qquad\qquad
		\begin{tikzpicture}[x=10mm,y=10mm,lab/.style={font=\small},pt/.style={circle,fill=black,inner sep=1.3pt}]
			\fill[gray!25] (0,1)--(-2,-0.4)--(0,-1)--(2,0.5)--cycle;
			\draw (0,1)--(-2,-0.4)--(0,-1)--(2,0.5)--cycle;
			\foreach \p in {(0,1),(0,-1),(-2,-0.4),(2,0.5)}{\node[pt] at \p {};}
			\node[lab,anchor=south east] at (0,1) {$x_0$};
			\node[lab,anchor=north east] at (0,-1) {$x_1$};
			\node[lab,anchor=east] at (-2,-0.4) {$x_2$};
			\node[lab,anchor=west] at (2,0.5) {$x_3$};
			\node[lab] at (0,0) {$P_2$};
			\node[font=\small] at (0,-2.5) {(B)};
		\end{tikzpicture}
		\caption{\small The sets $P_\iota$ for $k=3$, $d=2$ and four points $x_0,\dots,x_3$, where $x_0$ and $x_1$ lie on the line $\{s=0\}$ in $\R^2$ with coordinates $(s,t)$, and $x_2$ and $x_3$ lie on different sides of it.
			The triple $(f_1,f_2,f_3)$ with $f_1(s,t)=s$, $f_2=0$ and $f_3=-f_1$ gives $I_0=I_1=[3]$, $I_2=\{1,2\}$ and $I_3=\{2,3\}$.
			(A) The triangles $P_1=\conv\{x_0,x_1,x_2\}$ and $P_3=\conv\{x_0,x_1,x_3\}$, which are separated by the line $\{f_1=0\}$.
			(B) The quadrilateral $P_2=\conv\{x_0,x_1,x_2,x_3\}$.
			Here $D=\{1,3\}$, the stratum of the corresponding $G$-box consists of a single point, so it is closed, and the relative interiors of $P_1$, $P_2$ and $P_3$ have no common point, as claimed in Step~1.}
		\label{fig : P}
	\end{figure}

	\medskip
	Let $y\in\E^G$ correspond to the tuple $(g_1,\dots,g_k)=(\phi^y_1,\dots,\phi^y_k)$.
	According to \eqref{eq : values} we have that 
	\[
	\langle y,A(p_{j,i})\rangle=g_\iota(x_j)
	\] 
	whenever $i$ lies in the orbit $\iota$. 
	Consequently,  the inequalities $\langle y,A(u)\rangle\leq0$, for $u\in U$, of Lemma~\ref{lem : closed strata} read $g_\iota(x_j)\leq0$ whenever $\iota\in I_j$, and the inequalities $\langle y,A(v)\rangle>0$ for $v\in V-U$, read $g_\iota(x_j)>0$ whenever $\iota\notin I_j$.
	Since the stratum $N_G(U)$ is closed, Lemma~\ref{lem : closed strata} therefore says that for every non-zero tuple $(g_1,\dots,g_k)$ of affine functions with sum zero the following implication holds:
	\begin{equation}
		\label{eq : strict}
		(\forall j)\,(\forall \iota\in I_j)\ g_\iota(x_j)\leq0 \ \Longrightarrow \ (\forall j)\,(\forall \iota\notin I_j) \ g_\iota(x_j)>0.
	\end{equation}
	Fix $y\in N_G(U)$ and write $f_\iota:=\phi^y_\iota$.

	\medskip
\textit{Step 1: The relative interiors of the sets $P_\iota$, for $\iota\in[k]$, have no common point.}
	Suppose $z\in\relint(P_1)\cap\cdots\cap\relint(P_k)$.
	Since $y\in N_G(U)$, we have $I_j=\{\iota : f_\iota(x_j)\leq0\}$ for all $j$, as observed before Figure~\ref{fig : strata}.
	Thus, $f_\iota(x_j)\leq0$ for every point $x_j$ spanning $P_\iota$, that is, whenever $\iota\in I_j$.
	The function $f_\iota$ is affine, so it is non-positive on the convex hull $P_\iota$ of these points, and in particular $f_\iota(z)\leq0$ because  $z\in P_\iota$.
	This holds for every $\iota\in[k]$, and since $\sum_{\iota\in[k]} f_\iota(z)=0$, we get that $f_\iota(z)=0$ for all $\iota$.
	An affine function which is non-positive on a convex set and vanishes at a relative interior point vanishes on the set, so $f_\iota(x_j)=0$ whenever $\iota\in I_j$.
	For a column $j$ with a missing orbit, which exists as $U\neq V$, we get that  $\sum_{\iota\notin I_j}f_\iota(x_j)=-\sum_{\iota\in I_j}f_\iota(x_j)=0$, while all terms on the left are positive.
	This is a contradiction, so no such $z$ exists.

	\medskip
	{\em Step 2: $|D|\leq d+1$.}
	Suppose $|D|\geq d+2$.
	Then in particular $k\geq3$.
	Take $\iota_0\in D$, and let $j_0$ be a column in which $\iota_0$ is missing.
	Any non-zero tuple $g$ of affine functions with sum zero, with $g_\iota(x_j)\leq0$ whenever $\iota\in I_j$ and with $g_{\iota_0}=0$ violates~\eqref{eq : strict}, because $g_{\iota_0}(x_{j_0})=0$.
	We use this argument twice.
	
	First, suppose that $P_\iota=\emptyset$ for some $\iota\neq\iota_0$, that is, $\iota\notin I_j$ for all $j$.
	Choose any $\iota'\notin\{\iota,\iota_0\}$, which is possible as $k\geq3$, and let $g$ be the tuple with $g_\iota:=1$, $g_{\iota'}:=-1$ and all other entries $0$.
	This tuple is non-zero, has sum zero and satisfies $g_{\iota_0}=0$.
	We check that it also satisfies the inequalities $g_\kappa(x_j)\leq0$ for all $0\leq j\leq n$ and $\kappa\in I_j$:
	\begin{compactitem}[\quad --]
		\item the entry $g_\iota=1$ does not appear in any of these inequalities, because $\iota\notin I_j$ for all $j$,
		\item the entry $g_{\iota'}=-1$ is negative at every point, so $g_{\iota'}(x_j)\leq0$  for every $j$, independently of whether $\iota'\in I_j$ or not, and
		\item all other entries are $0$, and so $0\leq0$.
	\end{compactitem}
	Hence, $g$ violates~\eqref{eq : strict}, as explained above.
	This contradiction shows that all sets $P_\iota$, $\iota\neq\iota_0$, are non-empty.
	
	Second, suppose that the relative interiors of the sets $P_\iota$, $\iota\neq\iota_0$, have no common point.
	Lemma~\ref{lem : separation}, applied to these $k-1\geq2$ non-empty convex sets, gives affine functions $g_\iota$, $\iota\neq\iota_0$, not all zero, with sum zero and $g_\iota\leq0$ on $P_\iota$, and together with $g_{\iota_0}:=0$ we again obtain a contradiction.
	
	Since $\iota_0\in D$ was arbitrary, we have shown that for every $\iota_0\in D$ the sets $P_\iota$, $\iota\neq\iota_0$, are non-empty and their relative interiors have a common point.
	As $|D|\geq2$, every $P_\iota$ is non-empty, because we can choose $\iota_0\in D$ different from $\iota$.

	Let $\iota_1,\dots,\iota_{d+1}\in[k]$.
	Since $|D|\geq d+2$, there is an orbit $\iota_0\in D$ different from $\iota_1,\dots,\iota_{d+1}$.
	By what we have shown, the relative interiors of the sets $P_\iota$, $\iota\neq\iota_0$, have a common point, and in particular $\relint(P_{\iota_1}),\dots,\relint(P_{\iota_{d+1}})$ have a point in common.
	Thus, every $d+1$ members of the family have a common point, and by Helly's theorem~\cite{DanzerGrunbaumKlee1963} all $k$ relative interiors have a point in common.
	This contradicts Step~1, and so $|D|\leq d+1$.

	\medskip
	{\em Step 3: If $|D|\geq2$, then $I_j\cap D\neq\emptyset$ for every column $j$.}
	Suppose that $I_{j_0}\cap D=\emptyset$ for some column $j_0$.
	We show that then every tuple $g=(g_1,\dots,g_k)$ of affine functions with sum zero and $g_\iota(x_j)\leq0$ whenever $\iota\in I_j$ is zero.
	This is a contradiction, because the tuple $(f_1,\dots,f_k)=(\phi^y_1,\dots,\phi^y_k)$ of the point $y\in N_G(U)$ fixed before Step~1 is such a tuple, as $I_j=\{\iota : f_\iota(x_j)\leq0\}$, and it is non-zero, as $y\neq0$ and $y\longmapsto(\phi^y_1,\dots,\phi^y_k)$ is a linear isomorphism.

	First, we record what the assumption means.
	The condition $I_{j_0}\cap D=\emptyset$ says that every orbit of $D$ is missing in the column $j_0$.
	On the other hand, an orbit $\nu\notin D$ is, by the definition of $D$, missing in no column, that is, $\nu\in I_j$ for all $0\leq j\leq n$.
	Consequently, $g_\nu(x_j)\leq0$ for all $\nu\notin D$ and all $j$, and so the affine function $h:=\sum_{\nu\notin D}g_\nu$ satisfies $h(x_j)\leq0$ for all $j$.

	Next, fix $\iota_1\in D$ and move the entries of $g$ outside $D$ into the entry $\iota_1$, that is, define the tuple $g'$ by
	\[
	g'_{\iota_1}:=g_{\iota_1}+h,\qquad g'_\iota:=g_\iota\ \text{ for }\iota\in D-\{\iota_1\},\qquad g'_\nu:=0\ \text{ for }\nu\notin D.
	\]
	The tuple $g'$ has the same sum as $g$, namely zero.
	It also satisfies $g'_\iota(x_j)\leq0$ whenever $\iota\in I_j$: for $\iota=\iota_1$ because $g_{\iota_1}(x_j)\leq0$ and $h(x_j)\leq0$, for the other $\iota\in D$ because these entries are unchanged, and for $\nu\notin D$ because $g'_\nu=0$.
	Since the entries of $g'$ outside $D$ vanish, at the column $j_0$ we get $\sum_{\iota\in D}g'_\iota(x_{j_0})=0$.
	If $g'$ were non-zero, then~\eqref{eq : strict} would give $g'_\iota(x_{j_0})>0$ for all $\iota\in D$, because all these orbits are missing in the column $j_0$, and the sum would be positive.
	Hence, $g'=0$, that is, $g_\iota=0$ for $\iota\in D-\{\iota_1\}$ and $g_{\iota_1}=-h$.

	Finally, since $|D|\geq2$, we can repeat this with a second orbit $\iota_2\in D$ in place of $\iota_1$, and obtain also $g_{\iota_1}=0$.
	Thus, $g_\iota=0$ for all $\iota\in D$, and $h=-g_{\iota_1}=0$.
	The functions $g_\nu$, $\nu\notin D$, are therefore non-positive at all points $x_j$ and have sum $h=0$, so they vanish at all points $x_j$.
	As the points $x_0,\dots,x_n$ affinely span $\R^d$, the affine functions $g_\nu$ vanish identically.
	So $g=0$, which gives the announced contradiction.

	\medskip
	{\em Step 4: Conclusion.}
	Fix a column $j$.
	The orbits missing in this column form the set $[k]-I_j$, and by the definition of $D$ every such orbit belongs to $D$, that is, $[k]-I_j\subseteq D$.
	We show that $|[k]-I_j|\leq d$, which is the same as $|I_j|\geq k-d$.

	If $|D|\leq1$, then $|[k]-I_j|\leq|D|\leq1\leq d$.
	If $|D|\geq2$, then Step~3 gives $I_j\cap D\neq\emptyset$, that is, some orbit of $D$ is present in the column $j$ and so is not missing there.
	Hence $[k]-I_j$ is a proper subset of $D$, and with Step~2 we get $|[k]-I_j|\leq|D|-1\leq d$.
\end{proof}

\medskip
So if $k\geq d+2$, no column of a $G$-box with a non-empty closed stratum consists of a single orbit.

\bigskip
\section{Complexes of partial list colourings}
\label{sec : list colourings}
\bigskip

In this section we study complexes of partial proper list colourings, mainly of graphs in $\mathscr B$, and at the end of arbitrary graphs under a local condition on the list sizes.
Throughout this section $\FF$ is an arbitrary field, all homology groups are reduced and have coefficients in $\FF$, and all complexes contain the empty face, so that reduced homology is computed with augmented chains and is defined in degree $-1$.
 
\medskip
Let $F$ be a graph on a finite set $Y$, and let $L=(L_w : w\in Y)$ be a family of finite sets, the {\em lists}.
We regard every list $L_w$ as a $0$-dimensional simplicial complex with the vertex set $\{w\}\times L_w$, and consider the join $\bigast_{w\in Y}L_w$.
Its vertices are the pairs $(w,c)$ with $w\in Y$ and $c\in L_w$, and its faces are the sets $\sigma$ of such pairs which contain at most one pair $(w,c)$ for every $w\in Y$.
The complex of {\em partial proper list colourings} is the sub-complex
\[
\LL_F(L):=\big\{\sigma\in \bigast_{w\in Y}L_w : \{(x,c),(y,c)\}\not\subseteq\sigma\text{ for every edge }xy\text{ of }F\text{ and every colour }c\big\}.
\]
A face $\sigma$ of the join is a colouring of some of the vertices of $F$ from their lists, where $w$ gets the colour $c$ if $(w,c)\in\sigma$, and $\sigma$ lies in $\LL_F(L)$ if and only if adjacent coloured vertices get different colours.

\medskip
Both conditions in the definition of $\LL_F(L)$ concern pairs of vertices, so $\LL_F(L)$ is a flag complex, that is, an independence complex.
Recall that the {\em independence complex} $\Ind(\Lambda)$ of a graph $\Lambda$ is the simplicial complex whose faces are the sets of pairwise non-adjacent vertices of $\Lambda$.
Let $F\cart  L$ be the graph on the vertex set $\{(w,c) : w\in Y,\ c\in L_w\}$ in which $(w,c)$ and $(w',c')$ are adjacent if either $w=w'$ and $c\neq c'$, or $ww'$ is an edge of $F$ and $c=c'$.
Then
\[
\LL_F(L)=\Ind(F\cart  L).
\]
For every $w\in Y$ the {\em fibre} $\{w\}\times L_w$ is a clique of $F\cart  L$.
If all lists are equal to a set $C$, then $F\cart  L$ is the Cartesian product $F\cart K_C$ of $F$ with the complete graph on $C$.

\medskip
For a graph $H$ on the vertex set of $\Delta_n$ we take all lists equal to $[r]$.
Identifying the vertex $(p_j,i)$ of the join $[r]^{*(n+1)}$ with the vertex $p_{j,i}$, this join becomes $\K=Q_0*\cdots*Q_n$.
The {\em colouring complex} of $H$ is
\[
\LL_H:=\LL_H([r],\dots,[r])=\big\{J\in\K : \{p_{j,i},p_{j',i}\}\not\subseteq J\text{ for every edge }p_jp_{j'}\text{ of }H\text{ and every }i\in[r]\big\}.
\]
So $\LL_H$ consists of all faces of $\K$ which, viewed as partial labellings, do not give the same label to the two vertices of an edge of $H$.
A typical face of $\LL_H$ is $J=\{p_{j,c(j)} : j\in S\}$, where $S\subseteq\{0,\dots,n\}$ and $c\colon S\longrightarrow[r]$ is a map with $c(j)\neq c(j')$ for every edge $p_jp_{j'}$ of $H$ with $j,j'\in S$.
In other words, $\LL_H=\Ind(H\cart K_r)$ is the complex used by Engstr\"om~\cite{Engstrom2011}.
His criterion~\cite[Thm.\,2.2]{Engstrom2011} assumes that this complex is $(N-1)$-connected in the sense of homotopy.
For prime powers $r=q$ and the group $G$ with one orbit the only $G$-box is $V$, and the $G$-box criterion, Theorem~\ref{th : G-box}, asks only for the homological $(N-1)$-connectivity of $\LL_H$ over $\FF_p$.
Since every $(N-1)$-connected complex is homologically $(N-1)$-connected over every field, while the converse fails, Theorem~\ref{th : G-box} is more general than the criterion of Engstr\"om already for prime powers.
Beyond prime powers, Theorem~\ref{th : G-box} asks in addition for the connectivity of the independence complexes $\LL_H[U]=\Ind\big((H\cart K_r)[U]\big)$ of the induced sub-graphs on $G$-boxes $U$.

\begin{figure}[ht]
	\centering
	\begin{tikzpicture}[x=13mm,y=5.5mm,
			pt/.style={circle,fill=black,inner sep=1.5pt},
			hollow/.style={circle,draw,fill=white,inner sep=1.5pt},
			gv/.style={circle,fill=black,inner sep=1.7pt},
			lab/.style={font=\scriptsize}]
		\foreach \a/\b in {1/3,4/6}{\draw[dashed,rounded corners=2pt] (-0.45,-\a+0.45) rectangle (4.45,-\b-0.45);}
		\foreach \i in {1,...,6}{\node[lab,anchor=east] at (-0.55,-\i) {$\i$};}
		\foreach \j/\A/\B in {0/1/0,1/1/1,2/0/1,3/1/0,4/1/0}{
				\node[lab] at (\j,0.15) {$Q_{\j}$};
				\foreach \i in {1,2,3}{\ifnum\A=1 \node[pt] at (\j,-\i) {}; \else \node[hollow] at (\j,-\i) {}; \fi}
				\foreach \i in {4,5,6}{\ifnum\B=1 \node[pt] at (\j,-\i) {}; \else \node[hollow] at (\j,-\i) {}; \fi}}
		\foreach \j in {0,...,4}{\node[gv] (v\j) at (\j,-8) {}; \node[lab,below] at (\j,-8.2) {$p_{\j}$};}
		\draw[thick] (v0)--(v1)--(v2);
		\draw[thick] (v3)--(v4);
		\node[lab,anchor=west] at (4.6,-8) {$H$};
		\foreach \j/\L in {0/{\{1,2,3\}},1/{[6]},2/{\{4,5,6\}},3/{\{1,2,3\}},4/{\{1,2,3\}}}{\node[lab] at (\j,-9.6) {$\L$};}
		\node[lab,anchor=west] at (4.6,-9.6) {lists};
	\end{tikzpicture}
	\caption{\small A $G$-box $U$ for $r=6$, $n=4$ and the group $G\cong\Z/3$ of Lemma~\ref{lem : group} with $q=3$, whose orbits $\{1,2,3\}$ and $\{4,5,6\}$ are framed, and filled vertices belong to $U$.
		For the graph $H$ with the edges $p_0p_1$, $p_1p_2$ and $p_3p_4$, the induced sub-complex $\LL_H[U]$ is the complex of partial proper list colourings of $H$ with the lists $U\cap Q_j$ shown below the vertices.
		At the edge $p_3p_4$ both lists are the same orbit. Such columns are the reason for the deficits in the connectivity bounds of this section.}
	\label{fig : lists}
\end{figure}
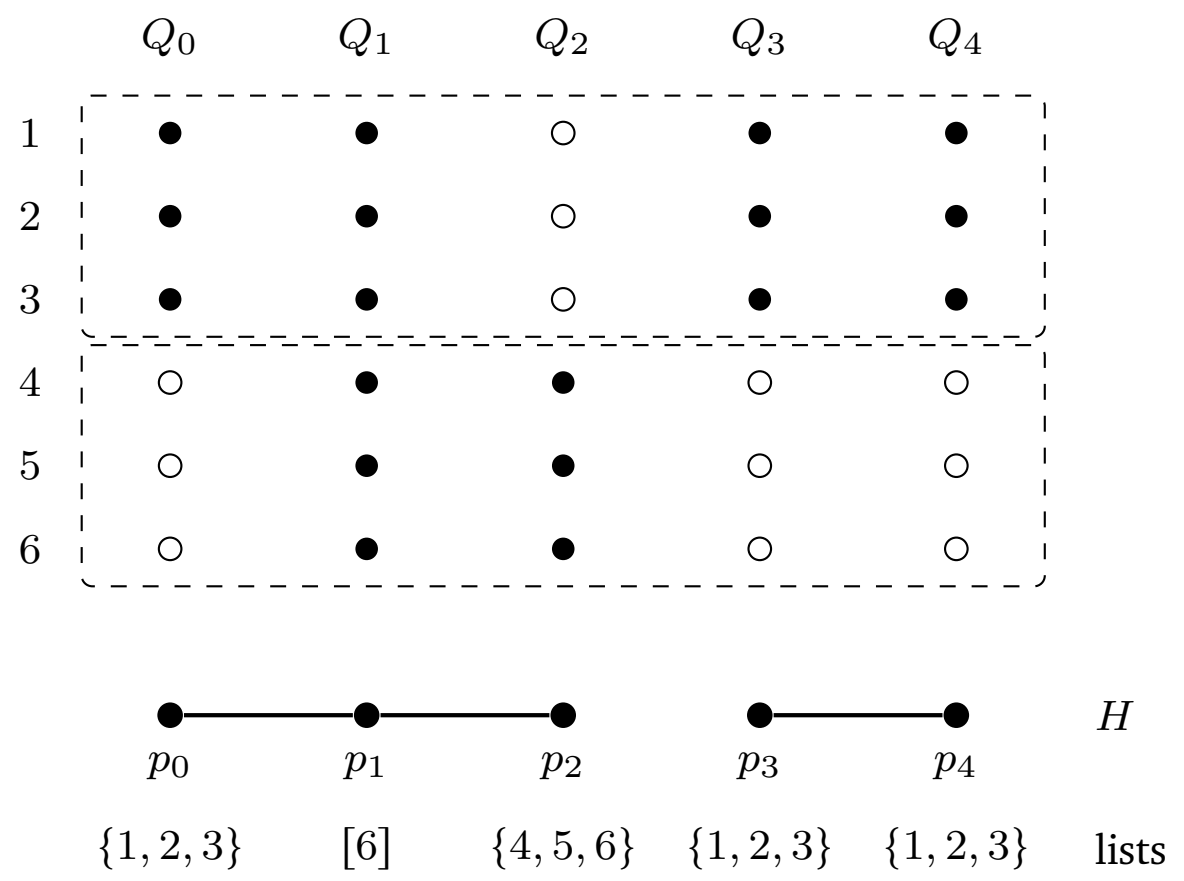

\medskip
If $U\subseteq V$ and $U_j:=U\cap Q_j$ is viewed as a subset of $[r]$, then the induced sub-complex $\LL_H[U]$ is the complex $\LL_H(U_0,\dots,U_n)\subseteq U_0*\cdots*U_n$ of partial proper list colourings of $H$ with the lists $U_0,\dots,U_n$.
Figure~\ref{fig : lists} shows an example.

\begin{figure}[ht]
	\centering
	\begin{tikzpicture}[x=12mm,y=10mm,pt/.style={circle,fill=black,inner sep=1.5pt},lab/.style={font=\small}]
		\foreach \x/\n in {0/y_1,1/w,2/y_2,3/z}{\node[pt] (a\x) at (\x,0) {}; \node[lab,above] at (\x,0.1) {$\n$};}
		\draw (a0)--(a1)--(a2)--(a3);
		\foreach \x/\L in {0/{\{1,2\}},1/{\{1,2,3\}},2/{\{1,3\}},3/{\{1,2\}}}{\node[lab] at (\x,-0.6) {$\L$};}
		\node[font=\small] at (1.5,-1.5) {(A)};
	\end{tikzpicture}
	\qquad
	\begin{tikzpicture}[x=12mm,y=10mm,pt/.style={circle,fill=black,inner sep=1.5pt},lab/.style={font=\small}]
		\foreach \x/\n in {0/y_1,1/w,2/y_2,3/z}{\node[pt] (a\x) at (\x,0) {}; \node[lab,above] at (\x,0.1) {$\n$};}
		\draw (a0)--(a1)--(a2)--(a3);
		\foreach \x/\L in {0/{\{1,2\}},1/{\{2,3\}},2/{\{1,3\}},3/{\{1,2\}}}{\node[lab] at (\x,-0.6) {$\L$};}
		\node[font=\small] at (1.5,-1.5) {(B)};
	\end{tikzpicture}
	\qquad
	\begin{tikzpicture}[x=12mm,y=10mm,pt/.style={circle,fill=black,inner sep=1.5pt},lab/.style={font=\small}]
		\foreach \x/\n in {0/y_1,2/y_2,3/z}{\node[pt] (a\x) at (\x,0) {}; \node[lab,above] at (\x,0.1) {$\n$};}
		\draw (a2)--(a3);
		\foreach \x/\L in {0/{\{2\}},2/{\{3\}},3/{\{1,2\}}}{\node[lab] at (\x,-0.6) {$\L$};}
		\node[font=\small] at (1.5,-1.5) {(C)};
	\end{tikzpicture}
	\caption{\small Deletion and link of the vertex $(w,1)$ in a colouring complex.
		(A) A path $F$ with the lists $L$ below its vertices.
		(B) The deletion $\LL_F(L)-(w,1)$ is the complex of $F$ in which $w$ loses the colour $1$.
		(C) The link of $(w,1)$ is the complex of $F-w$ in which the neighbours $y_1$ and $y_2$ of $w$ lose the colour $1$.}
	\label{fig : deletion and link}
\end{figure}
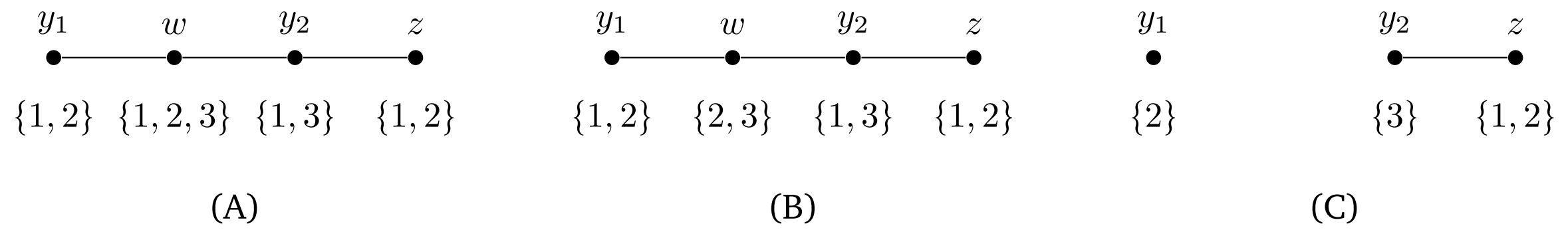

\medskip
Now we recall the notions which we use to study the connectivity of colouring complexes.
For a vertex $u$ of a graph $F$ we denote by $F-u$ the graph obtained by deleting $u$ and the edges incident to $u$, and for lists $L$ of $F$ we keep the notation $L$ for the lists $(L_y : y\in Y-\{u\})$ of $F-u$.
Similarly, for a vertex $x$ of a simplicial complex $K$ let $K-x:=\{\sigma\in K : x\notin\sigma\}$, which is the sub-complex induced on all vertices except $x$, and let $\operatorname{lk}_Kx:=\{\sigma\in K : x\notin\sigma,\ \sigma\cup\{x\}\in K\}$ be the link of $x$.
If $K=\Ind(\Lambda)$ is an independence complex, then $K-x=\Ind(\Lambda-x)$ and $\operatorname{lk}_Kx=\Ind(\Lambda-N[x])$, where $N[x]$ denotes the set consisting of $x$ and its neighbours in $\Lambda$.
The closed star of $x$ is the cone $x*\operatorname{lk}_Kx$, which has vanishing reduced homology, and $(x*\operatorname{lk}_Kx)\cap(K-x)=\operatorname{lk}_Kx$.
Hence the Mayer--Vietoris sequence~\cite[\S25]{Munkres1984}, \cite[Sec.\,2.2]{Hatcher2002} of the decomposition $K=(K-x)\cup(x*\operatorname{lk}_Kx)$, which is exact also in degree $-1$ in augmented homology,  gives the exact sequences
\[
\widetilde H_i(K-x)\longrightarrow\widetilde H_i(K)\longrightarrow\widetilde H_{i-1}(\operatorname{lk}_Kx)
\qquad\text{and}\qquad
\widetilde H_i(\operatorname{lk}_Kx)\longrightarrow\widetilde H_i(K-x)\longrightarrow\widetilde H_i(K).
\]
They imply the following two statements, which we use repeatedly:
\begin{compactenum}[\rm\quad (D1)]
	\item if $K-x$ is homologically $\tau$-connected and $\operatorname{lk}_Kx$ is homologically $(\tau-1)$-connected, then $K$ is homologically $\tau$-connected, compare~\cite[Lem.\,3.1]{Engstrom2011}, and
	\item if $K$ and $\operatorname{lk}_Kx$ are homologically $\tau$-connected, then $K-x$ is homologically $\tau$-connected.
\end{compactenum}
For a colouring complex and a vertex $(w,c)$ both complexes are colouring complexes again:
\begin{compactitem}[\quad --]
	\item $\LL_F(L)-(w,c)$ is the complex of $F$ with the list of $w$ replaced by $L_w-\{c\}$, and
	\item $\operatorname{lk}_{\LL_F(L)}(w,c)$ is the complex of $F-w$ with the lists $L_y-\{c\}$ for the neighbours $y$ of $w$ and the lists $L_y$ for the other vertices $y$ of $F-w$.
\end{compactitem}
Figure~\ref{fig : deletion and link} shows an example.

\medskip
Next we decompose a colouring complex along the colours of one vertex.
Let $u$ be a vertex of $F$.
We consider the following complexes:
\begin{compactitem}[\quad --]
	\item $Z:=\LL_{F-u}(L)$, the complex of partial proper list colourings of $F-u$ with the lists $L_y$, $y\in Y-\{u\}$, and
	\item for $c\in L_u$, the complex $Z_c:=\LL_{F-u}(L^c)$, where $L^c_y:=L_y-\{c\}$ for the neighbours $y$ of $u$ in $F$ and $L^c_y:=L_y$ for the other vertices $y$ of $F-u$.
\end{compactitem}
So $Z_c\subseteq Z$, and $Z_c$ consists of the faces of $Z$ which do not give the colour $c$ to a neighbour of $u$.
For $\Gamma\subseteq L_u$ set $Z_\Gamma:=\bigcap_{c\in\Gamma}Z_c$.
This is the complex of $F-u$ in which the neighbours of $u$ lose all colours of $\Gamma$.

\begin{figure}[ht]
	\centering
	\begin{tikzpicture}[scale=1.1,lab/.style={font=\small},
			zp/.style={circle,fill=black,inner sep=2.2pt},
			ap/.style={circle,draw=black,thick,fill=white,inner sep=2.2pt}]
		\foreach \k/\n/\t in {0/{(u,1)}/ap,1/{(y,2)}/zp,2/{(u,3)}/ap,3/{(y,1)}/zp,4/{(u,2)}/ap,5/{(y,3)}/zp}{\coordinate (h\k) at (90+60*\k:1.5);}
		\draw (h0)--(h1)--(h2)--(h3)--(h4)--(h5)--cycle;
		\draw[line width=1.6pt] (h5)--(h0)--(h1);
		\foreach \k/\n/\t in {0/{(u,1)}/ap,1/{(y,2)}/zp,2/{(u,3)}/ap,3/{(y,1)}/zp,4/{(u,2)}/ap,5/{(y,3)}/zp}{\node[\t] at (h\k) {}; \node[lab] at (90+60*\k:1.95) {$\n$};}
	\end{tikzpicture}
	\caption{\small The complex $\LL_F(L)$ for the graph $F$ with one edge $uy$ and the lists $L_u=L_y=[3]$ is a hexagon.
		Here $Z$ consists of the three filled vertices $(y,1)$, $(y,2)$, $(y,3)$, and $Z_c$ of the two filled vertices $(y,c')$ with $c'\neq c$.
		The hexagon is the union of the three cones $(u,c)*Z_c$, and the cone $(u,1)*Z_1$ is drawn thick.
		Any two of the cones intersect in a single vertex $Z_{\{c,c'\}}$ of $Z$.}
	\label{fig : cones}
\end{figure}
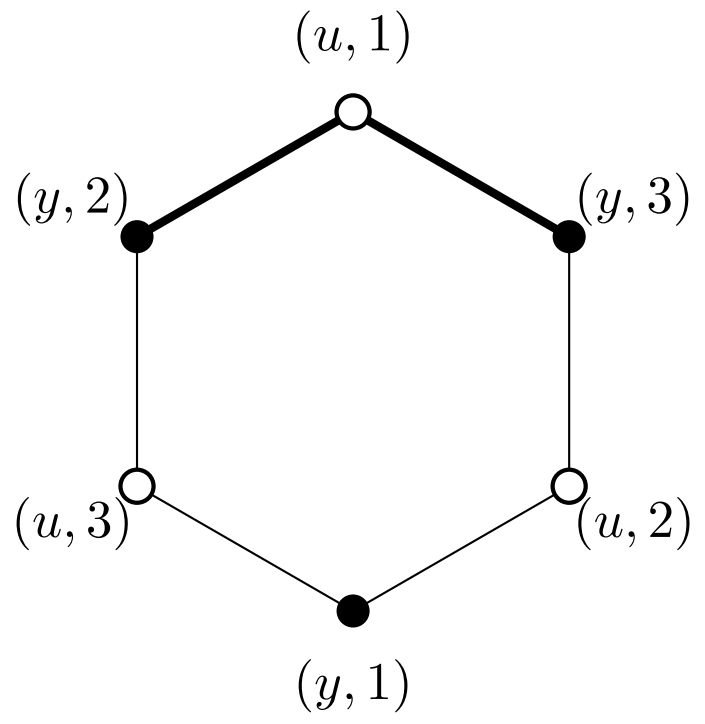

\medskip
A face of $\LL_F(L)$ either does not colour the vertex $u$ of $F$, and then it lies in $Z$, or it contains exactly one vertex $(u,c)$, and then the remaining vertices form a face of $Z_c$.
Hence
\begin{equation}
	\label{eq : decomposition}
	\LL_F(L)\ =\ Z\cup\bigcup_{c\in L_u}(u,c)*Z_c,
\end{equation}
the cone $(u,c)*Z_c$ intersects $Z$ in $Z_c$, and a face which lies in two different cones contains neither of their apices.
Figure~\ref{fig : cones} shows this decomposition in the smallest interesting case.

\medskip
\begin{lemma}
	\label{lem : cones}
	Let $u$ be a vertex of $F$, let $\tau\geq0$, and let $S\subseteq L_u$ be a non-empty set such that
	\begin{compactenum}[\rm \quad (A)]
		\item $Z=\bigcup_{c\in S}Z_c$,
		\item for every $\Gamma\subseteq S$ with $|\Gamma|\geq2$ the complex $Z_\Gamma$ has a vertex and is homologically $(\tau-|\Gamma|+1)$-connected, and
		\item for every $c\in L_u$ the complex $Z_c$ is homologically $(\tau-1)$-connected.
	\end{compactenum}
	Then $\LL_F(L)$ is homologically $\tau$-connected.
\end{lemma}
\begin{proof}
	The proof has two steps.
	First, we show with Meshulam's nerve theorem that the union $W$ of the cones $(u,c)*Z_c$ with $c\in S$ is homologically $\tau$-connected.
	Second, we add the remaining cones one by one and show with (D1) that this preserves homological $\tau$-connectivity.
	By the decomposition~\eqref{eq : decomposition}, the result is $\LL_F(L)$.

	\smallskip
	{\em Step 1.}
	Let $W:=\bigcup_{c\in S}(u,c)*Z_c$.
	By (A), $W$ contains $Z$.
	The cones $(u,c)*Z_c$, $c\in S$, cover $W$.
	A face which lies in two different cones contains neither apex, and so for $\Gamma\subseteq S$ with $|\Gamma|\geq2$ the intersection of the cones with $c\in\Gamma$ is $Z_\Gamma$.
	
	We use Meshulam's homological nerve theorem~\cite[Thm.\,2.1]{Meshulam2001} in the following special case: if a finite simplicial complex $X$ is the union of sub-complexes $X_c$, $c\in S$, such that for every non-empty $\Gamma\subseteq S$ the intersection $X_\Gamma:=\bigcap_{c\in\Gamma}X_c$ has a vertex and is homologically $(\tau-|\Gamma|+1)$-connected, then the nerve of the cover is the full simplex on $S$, and the complex $X$ is homologically $\tau$-connected.
	
	For the cones $X_c:=(u,c)*Z_c$ these assumptions hold: for $|\Gamma|=1$ because cones are contractible, and for $|\Gamma|\geq2$ by (B).
	So $W$ is homologically $\tau$-connected.
	
	\smallskip
	{\em Step 2.}
	Let $c_1,\dots,c_t$ be the elements of $L_u-S$, and for $0\leq s\leq t$ let
	\[
		W_s:=W\cup(u,c_1)*Z_{c_1}\cup\cdots\cup(u,c_s)*Z_{c_s}.
	\]
	Then $W_0=W$, and $W_t=\LL_F(L)$ by the decomposition~\eqref{eq : decomposition}, since $Z\subseteq W$.
	Let $1\leq s\leq t$.
	The vertex $(u,c_s)$ of $W_s$ lies only in the cone $(u,c_s)*Z_{c_s}$.
	Hence its link in $W_s$ is $Z_{c_s}$, and its deletion is $W_s-(u,c_s)=W_{s-1}\cup Z_{c_s}=W_{s-1}$, because $Z_{c_s}\subseteq Z\subseteq W_{s-1}$.
	By (C) the link is homologically $(\tau-1)$-connected, so by (D1) the complex $W_s$ is homologically $\tau$-connected if $W_{s-1}$ is.
	Induction on $s$ shows that the complex $W_t=\LL_F(L)$ is homologically $\tau$-connected.
\end{proof}

\medskip
The next lemma gives a second way to show connectivity.
We fix a vertex $(w,c^*)$ of $\LL_F(L)$ and delete, one at a time, the vertices which prevent it from being a cone point: the vertices $(y,c^*)$ for the neighbours $y$ of $w$, and the vertices $(w,c)$ with $c\neq c^*$.
After these deletions we obtain a cone with apex $(w,c^*)$, and (D1) transports its connectivity back to $\LL_F(L)$, provided the links of the deleted vertices are sufficiently connected.

\medskip
\begin{lemma}
	\label{lem : chain}
	Let $\tau\geq-1$, let $w$ be a vertex of $F$ and $c^*\in L_w$, let $y_1,\dots,y_j$ be all neighbours of $w$ whose lists contain $c^*$, in any order, and let $c_1,\dots,c_m$ be all elements of $L_w-\{c^*\}$, so $m=|L_w|-1$.
	Assume that the following complexes are homologically $(\tau-1)$-connected:
	\begin{compactenum}[\rm\quad (A)]
		\item for $1\leq i\leq j$, the complex of $F-y_i$ in which $c^*$ is removed from the lists of the neighbours of $y_i$ and of $y_1,\dots,y_{i-1}$, and
		\item for $1\leq h\leq m$, the complex of $F-w$ in which $c^*$ is removed from the lists of $y_1,\dots,y_j$ and $c_h$ from the lists of all neighbours of $w$.
	\end{compactenum}
	Then $\LL_F(L)$ is homologically $\tau$-connected.
\end{lemma}
\begin{proof}
	Let $x_1,\dots,x_{j+m}$ be the vertices
	\[
		(y_1,c^*),\ \dots,\ (y_j,c^*),\ (w,c_1),\ \dots,\ (w,c_m)
	\]
	of $\LL_F(L)$, in this order.
	Let $A_0:=\LL_F(L)$, and for $1\leq i\leq j+m$ let $A_i:=A_{i-1}-x_i$ be the complex obtained by deleting $x_1,\dots,x_i$.
	By the description of deletions, $A_i$ is the colouring complex of $F$ in which $c^*$ is removed from the lists of $y_1,\dots,y_{\min\{i,j\}}$ and, for $i>j$, the colours $c_1,\dots,c_{i-j}$ are removed from the list of $w$.

	\smallskip
	First, $A_{j+m}$ is homologically $\tau$-connected.
	Indeed, in $A_{j+m}$ the list of $w$ is $\{c^*\}$ and no neighbour of $w$ has $c^*$ in its list.
	So a face of $A_{j+m}$ which does not contain $(w,c^*)$ does not colour $w$, and adding $(w,c^*)$ to it gives again a proper colouring.
	Thus, $A_{j+m}$ is a cone with apex $(w,c^*)$, and it has vanishing reduced homology.

	\smallskip
	Second, by the description of links, the link of $x_i$ in $A_{i-1}$ is the complex from the assumption (A) of the lemma for $1\leq i\leq j$, and the complex from the assumption (B) of the lemma for $i=j+h$ with $1\leq h\leq m$.
	So, by assumption, these links are homologically $(\tau-1)$-connected.
	Since $A_i=A_{i-1}-x_i$, statement (D1) shows that $A_{i-1}$ is homologically $\tau$-connected if $A_i$ is.
	Applying this for $i=j+m,j+m-1,\dots,1$, we obtain that $A_0=\LL_F(L)$ is homologically $\tau$-connected.
\end{proof}

\begin{figure}[ht]
	\centering
	\begin{tikzpicture}[x=13mm,y=9mm,pt/.style={circle,fill=black,inner sep=1.4pt},lab/.style={font=\small}]
		\node[pt] (z) at (-2,0) {}; \node[pt] (y1) at (-1,0) {}; \node[pt] (w) at (0,0) {}; \node[pt] (y2) at (1,0) {}; \node[pt] (y3) at (0,-1) {};
		\draw (z)--(y1)--(w)--(y2) (w)--(y3);
		\node[lab,above] at (z) {$z$}; \node[lab,above] at (y1) {$y_1$}; \node[lab,above] at (w) {$w$}; \node[lab,above] at (y2) {$y_2$}; \node[lab,left] at (y3) {$y_3$};
		\node[lab,below] at (-2,-0.15) {$\{1,2\}$}; \node[lab,below] at (-1,-0.15) {$\{1,2,3\}$}; \node[lab,above] at (0,0.45) {$\{1,2\}$}; \node[lab,below] at (1,-0.15) {$\{1,3\}$}; \node[lab,right] at (0.1,-1) {$\{2,3\}$};
		\node[font=\small] at (-0.5,-2) {(A)};
	\end{tikzpicture}
	\qquad
	\begin{tikzpicture}[x=13mm,y=9mm,pt/.style={circle,fill=black,inner sep=1.4pt},lab/.style={font=\small}]
		\node[pt] (z) at (-2,0) {}; \node[pt] (w) at (0,0) {}; \node[pt] (y2) at (1,0) {}; \node[pt] (y3) at (0,-1) {};
		\draw (w)--(y2) (w)--(y3);
		\node[lab,above] at (z) {$z$}; \node[lab,above] at (w) {$w$}; \node[lab,above] at (y2) {$y_2$}; \node[lab,left] at (y3) {$y_3$};
		\node[lab,below] at (-2,-0.15) {$\{2\}$}; \node[lab,above] at (0,0.45) {$\{2\}$}; \node[lab,below] at (1,-0.15) {$\{1,3\}$}; \node[lab,right] at (0.1,-1) {$\{2,3\}$};
		\node[font=\small] at (-0.5,-2) {(B)};
	\end{tikzpicture}

	\medskip
	\begin{tikzpicture}[x=13mm,y=9mm,pt/.style={circle,fill=black,inner sep=1.4pt},lab/.style={font=\small}]
		\node[pt] (z) at (-2,0) {}; \node[pt] (y1) at (-1,0) {}; \node[pt] (w) at (0,0) {}; \node[pt] (y3) at (0,-1) {};
		\draw (z)--(y1)--(w) (w)--(y3);
		\node[lab,above] at (z) {$z$}; \node[lab,above] at (y1) {$y_1$}; \node[lab,above] at (w) {$w$}; \node[lab,left] at (y3) {$y_3$};
		\node[lab,below] at (-2,-0.15) {$\{1,2\}$}; \node[lab,below] at (-1,-0.15) {$\{2,3\}$}; \node[lab,above] at (0,0.45) {$\{2\}$}; \node[lab,right] at (0.1,-1) {$\{2,3\}$};
		\node[font=\small] at (-0.5,-2) {(C)};
	\end{tikzpicture}
	\qquad
	\begin{tikzpicture}[x=13mm,y=9mm,pt/.style={circle,fill=black,inner sep=1.4pt},lab/.style={font=\small}]
		\node[pt] (z) at (-2,0) {}; \node[pt] (y1) at (-1,0) {}; \node[pt] (y2) at (1,0) {}; \node[pt] (y3) at (0,-1) {};
		\draw (z)--(y1);
		\node[lab,above] at (z) {$z$}; \node[lab,above] at (y1) {$y_1$}; \node[lab,above] at (y2) {$y_2$}; \node[lab,left] at (y3) {$y_3$};
		\node[lab,below] at (-2,-0.15) {$\{1,2\}$}; \node[lab,below] at (-1,-0.15) {$\{3\}$}; \node[lab,below] at (1,-0.15) {$\{3\}$}; \node[lab,right] at (0.1,-1) {$\{3\}$};
		\node[font=\small] at (-0.5,-2) {(D)};
	\end{tikzpicture}
	\caption{\small The complexes of Lemma~\ref{lem : chain} for $c^*=1$.
		(A) The graph $F$ with its lists. Here $y_1,y_2$ are the neighbours of $w$ whose lists contain $1$, so $j=2$, and $L_w-\{1\}=\{2\}$, so $m=1$ and $c_1=2$.
		(B) The complex in assumption (A) of the lemma for $i=1$: the graph $F-y_1$, in which the neighbours $z$ and $w$ of $y_1$ lose the colour $1$.
		(C) The complex in assumption (A) of the lemma for $i=2$: the graph $F-y_2$, in which the neighbour $w$ of $y_2$ and the vertex $y_1$ lose the colour $1$.
		(D) The complex in assumption (B) of the lemma for $h=1$: the graph $F-w$, in which $y_1$ and $y_2$ lose the colour $1$ and all neighbours $y_1,y_2,y_3$ of $w$ lose the colour $2$.}
	\label{fig : chain}
\end{figure}
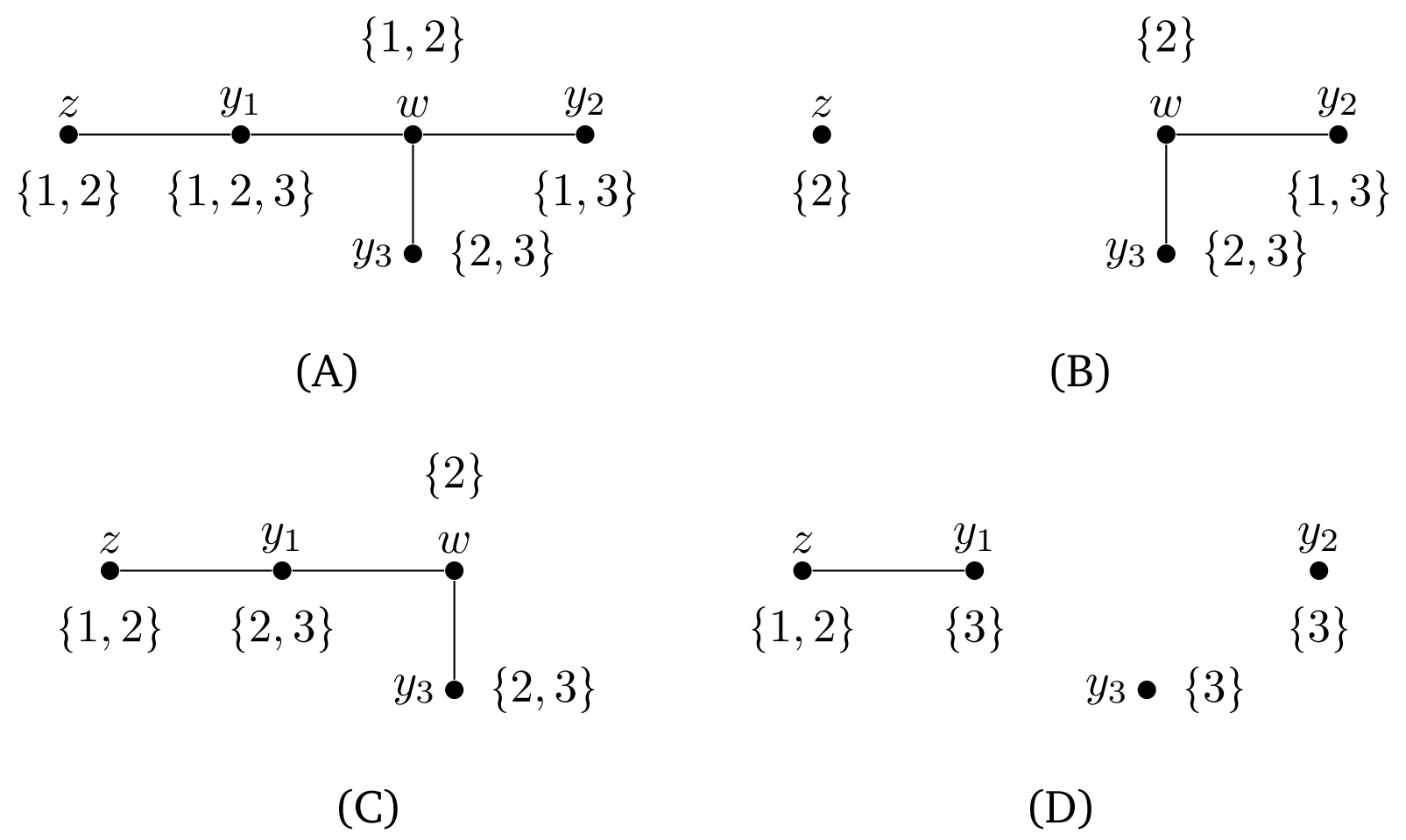

\medskip
Figure~\ref{fig : chain} shows the complexes in the assumptions (A) and (B) of Lemma~\ref{lem : chain} in an example.

\medskip
Recall that a vertex of a graph is {\em simplicial} if its neighbours are pairwise adjacent.
By Dirac's theorem~\cite{Dirac1961} a chordal graph with at least two vertices has two simplicial vertices, which can be chosen non-adjacent if the graph is not complete. For a complete graph every vertex is simplicial.

\medskip 
Let $F\in\mathscr B$.
Every block of an induced sub-graph of $F$ is an induced sub-graph of a block of $F$.
Next, an induced sub-graph of a chordal graph is chordal.
Further, a block of an induced sub-graph of a cycle is a vertex, an edge or the whole cycle.
In summary,  induced sub-graphs of $F$ belong to $\mathscr B$.
An induced cycle of length at least four has no cut vertex, so it lies in one block.
It cannot lie in a chordal block, and in a block which is a cycle it is the whole block.
Thus, the induced cycles of length at least four of $F$ are exactly its blocks which are cycles of length at least four.
We now bound from below the homological connectivity of the colouring complexes of graphs in $\mathscr B$, in terms of the degrees of the vertices and the sizes of their lists.
For this we use rooted graphs, as defined in Section~\ref{sec : intro}.
The root of a component is the last vertex which is removed in the induction below, and it needs one more colour.

\medskip
\begin{theorem}
	\label{th : colourings}
	Let $F\in\mathscr B$ be a rooted graph on a set $Y$ with $v$ vertices, and let $L=(L_w : w\in Y)$ be lists.
	For $w\in Y$ let $\chi(w):=1$ if $w$ is a root or lies on an induced cycle of length at least four, $\chi(w):=0$ otherwise, and let
	\[
		\kappa:=\sum_{w\in Y}\max\{0,\,2\deg_F(w)+\chi(w)-|L_w|\}.
	\]
	Then $\LL_F(L)$ is homologically $(v-2-\kappa)$-connected.
	In particular, if $|L_w|\geq2\deg_F(w)+1$ for all $w\in Y$, then $\LL_F(L)$ is homologically $(v-2)$-connected.
\end{theorem}
\begin{proof}
	For $w\in Y$ we set
	\[
		\varphi(w):=2\deg_F(w)+\chi(w)\qquad\text{and}\qquad\kappa_w:=\max\{0,\,\varphi(w)-|L_w|\}.
	\]
	Thus $\kappa=\sum_{w\in Y}\kappa_w$, and $\kappa_w=0$ if and only if $|L_w|\geq\varphi(w)$.
	In particular, $\kappa=0$ means that every list $L_w$ has at least $\varphi(w)$ elements.

	\medskip
	The idea of the proof is to express $\LL_F(L)$ through colouring complexes of graphs with fewer vertices, for which the theorem holds by induction, and to transfer their connectivity to $\LL_F(L)$ with (D2) and with Lemma~\ref{lem : cones} and Lemma~\ref{lem : chain}.
	In Step~1 we reduce to the case $\kappa=0$: we add a new colour to the list of a vertex $w$ with $\kappa_w\geq1$, and recover the old complex as a deletion with (D2).
	For $\kappa=0$ we remove a vertex: a component with one vertex is split off as a join factor (Step~2).
	Otherwise, a simplicial vertex of $F$, which is not a root, is removed with Lemma~\ref{lem : chain} (Step~4).
	If there is no such vertex, a vertex of degree two on a cycle is removed with Lemma~\ref{lem : cones} (Step~5). Step~3 shows that one of these two vertices exists.
	The assumptions of the two lemmas concern colouring complexes of graphs $F-x$ with modified lists, and the main point is to check that their sums $\kappa$ are small enough for the induction.

	\medskip
	When a vertex $x$ is deleted from $F$, we keep the roots of $F$ different from $x$, and in every connected component of $F-x$ without a root we choose a neighbour of $x$ as its root. We call these roots {\em new}.
	The rooted graph $F-x$ belongs to $\mathscr B$, and we denote by $\chi'$ and $\varphi'(y):=2\deg_{F-x}(y)+\chi'(y)$, $y\in Y-\{x\}$, the functions $\chi$ and $\varphi$ of the rooted graph $F-x$.
	Then the following holds.
	\begin{compactitem}[\quad(O)]
		\item For every vertex $y\neq x$ we have $\varphi'(y)\leq\varphi(y)$.
		Moreover, $\varphi'(y)\leq\varphi(y)-2$ if $y$ is a neighbour of $x$ which is not a new root, and $\varphi'(y)\leq\varphi(y)-1$ if $y$ is a new root.
	\end{compactitem}
	Indeed, deleting $x$ lowers the degrees of the neighbours of $x$ by one, every induced cycle of $F-x$ is an induced cycle of $F$, and a new root is a neighbour of $x$.

	\medskip
	We prove the theorem by induction on the number of vertices $v$ of $F$, and for fixed $v$ by induction on $\kappa$.
	So we assume that the theorem holds for all rooted graphs in $\mathscr B$ with fewer than $v$ vertices and all lists, and for the graph $F$ with all lists whose sum is smaller than $\kappa$.
	For $v=0$ the claim is empty.
	Let $v\geq1$ and $\tau:=v-2-\kappa$.

	\medskip
	{\em Step 1: We may assume $\kappa=0$.}
	Let $\kappa_w\geq1$, and add a new colour $c$, which lies in no list, to the list of $w$.
	For the new complex $A$ the sum $\kappa$ drops by one, so $A$ is homologically $(\tau+1)$-connected by induction.
	The link of $(w,c)$ in $A$ is $\LL_{F-w}(L)$, whose sum is at most $\kappa-\kappa_w$ by (O), and it is homologically $((v-1)-2-\kappa+1)$-connected, that is $\tau$-connected, by induction.
	By (D2), $\LL_F(L)=A-(w,c)$ is homologically $\tau$-connected.

	\medskip
	{\em Step 2: A connected component with one vertex.}
	Let $\kappa=0$.
	If a component of $F$ consists of one vertex $w$, then $\LL_F(L)=L_w*\LL_{F-w}(L)$, where $L_w\neq\emptyset$ since $\varphi(w)=1$.
	The complex $\LL_{F-w}(L)$ is homologically $(v-3)$-connected by induction.
	By the K\"unneth formula for joins~\cite[Lem.\,2.1]{Milnor1956}, over a field the join of a homologically $\tau_1$-connected complex and a homologically $\tau_2$-connected complex is homologically $(\tau_1+\tau_2+2)$-connected, and with $\tau_1=-1$ and $\tau_2=v-3$ this gives the claim.

	\medskip
	{\em Step 3: The choice of a vertex.}
	Let $\kappa=0$ and let all components of $F$ have at least two vertices.
	We claim that $F$ has a simplicial vertex which is not a root, or a vertex $u$ which is not a root, lies on a block $B$ which is a cycle of length at least four, and whose two neighbours lie on $B$ and have degree two in $F$.
	Let $C$ be a component of $F$ with root $\rho$.
	If $C$ is a single block, then either it is chordal and has two simplicial vertices, one of which is not $\rho$, or it is a cycle of length at least four, all of whose vertices have degree two, and we take $u\neq\rho$ on it.
	Otherwise $C$ has at least two end blocks, that is blocks containing exactly one cut vertex of $C$.
	Since $\rho$ is a non-cut vertex of at most one block, there is an end block $B$, with cut vertex $s$, such that $\rho$ is not a vertex of $B$ other than $s$.
	If $B$ is chordal, it has a simplicial vertex different from $s$, and this vertex is simplicial in $F$, since all its neighbours lie in $B$.
	When $B$ is a cycle of length at least four, we take for $u$ a vertex of $B$ which is neither $s$ nor a neighbour of $s$.

	\medskip
	{\em Step 4: A simplicial vertex $w$ which is not a root.}
	Let $c^*\in L_w$ and let $y_1,\dots,y_j$ and $c_1,\dots,c_m$ be as in Lemma~\ref{lem : chain}.
	The complexes in (A) and (B) are complexes of graphs with $v-1$ vertices, and by induction it suffices to show that their sums are $0$.

	\smallskip
	In (A) the vertices $y_1,\dots,y_{i-1}$ are neighbours of $y_i$, since $w$ is simplicial, so exactly the neighbours of $y_i$ lose $c^*$.
	A neighbour of $y_i$ loses at most one colour, and by (O) its value of $\varphi$ drops by at least one. The other lists are unchanged.

	\smallskip
	In (B) a neighbour of $w$ loses at most two colours, and all other lists are unchanged.
	The vertex $w$ is not a root, and deleting it does not disconnect its component, since its neighbours are pairwise adjacent. So $F-w$ has the same roots as $F$, and by (O) the value of $\varphi$ at a neighbour of $w$ drops by two.

	\medskip
	{\em Step 5: A vertex $u$ on a cycle $B$ as in Step 3.}
	Let $y_1,y_2$ be the neighbours of $u$.
	The vertex $u$ is not a root, and deleting it does not disconnect its component, since $u$ lies on a cycle. So $F-u$ has the same roots as $F$.
	In $F-u$ the vertices $y_1,y_2$ have degree one and lie on no cycle, so their values of $\varphi$ drop from $5$ to $2$, or to $3$ for the root, and at most one of them is a root. All other values of $\varphi$ do not increase.
	We apply Lemma~\ref{lem : cones} with $\tau=v-2$.
	Choose $S\subseteq L_u$ with three elements, which is possible since $|L_u|\geq\varphi(u)=5$.
	A face of $Z$ colours $y_1,y_2$ with at most two colours, so it lies in $Z_c$ for some $c\in S$, and (A) holds.
	In $Z_\Gamma$, $\Gamma\subseteq L_u$, the lists of $y_1,y_2$ have at least $5-|\Gamma|$ elements, and all other lists are unchanged.
	For $|\Gamma|=1$ the sum $\kappa$ of $Z_\Gamma$ is $0$, so every $Z_c$ is homologically $(v-3)$-connected by induction, which is (C).
	For $\Gamma\subseteq S$ with $t:=|\Gamma|\in\{2,3\}$ the sum of $Z_\Gamma$ is $0$ for $t=2$ and at most $1$ for $t=3$, so $Z_\Gamma$ is homologically $(v-1-t)$-connected by induction, and it has a vertex, since the list of $y_1$ is not empty. This is (B).
	Lemma~\ref{lem : cones} gives the claim.
\end{proof}

\medskip
For a chordal graph the theorem says that lists of size $2\deg_F(w)$ suffice, except for one vertex in every component.
The bound cannot be lowered further: for the triangle $K_3$ with equal lists of size $4$ the complex $\LL_{K_3}(L)$ is the chessboard complex $\Delta_{3,4}$, which is a torus, and for the edge $K_2$ with equal lists of size $2$ the complex is disconnected.
Applied to the complete graph $K_\ell$ and lists of size at least $2\ell-1$, the theorem gives the $(\ell-2)$-connectivity which corresponds to Hell's condition.

\medskip
\begin{proposition}
	\label{prop : stars}
	Let $k\geq1$, and let $F$ be the star $K_{1,k}$ with centre $u$.
	If $|L_u|\geq k+1$ and $|L_w|\geq k+2$ for the leaves $w$, then $\LL_F(L)$ is homologically $(k-1)$-connected.
\end{proposition}
\begin{proof}
	We apply Lemma~\ref{lem : cones} with $\tau=k-1$.
	The complex $Z$ is the join of the $k$ lists of the leaves, viewed as $0$-dimensional complexes, and $Z_\Gamma$ is the join of these lists with the colours of $\Gamma$ removed.
	Choose $S\subseteq L_u$ with $k+1$ elements.
	A face of $Z$ uses at most $k$ colours, so it lies in $Z_c$ for some $c\in S$, and (A) holds.
	For $c\in L_u$ and $\Gamma\subseteq S$ the lists of the leaves in $Z_c$ and $Z_\Gamma$ are non-empty, since they lose at most $k+1$ of their at least $k+2$ elements.
	By the K\"unneth formula for joins~\cite[Lem.\,2.1]{Milnor1956}, a join of $k$ non-empty $0$-dimensional complexes is homologically $(k-2)$-connected.
	Hence every $Z_c$ is homologically $(\tau-1)$-connected, which is (C), and for $|\Gamma|\geq2$ the complex $Z_\Gamma$ has a vertex and is homologically $(k-2)$-connected, hence $(\tau-|\Gamma|+1)$-connected, which is (B).
	Lemma~\ref{lem : cones} gives the claim.
\end{proof}

\medskip
For complete graphs with short lists we use the bound of Bj\"orner, Lov\'asz, Vre\'cica and \v{Z}ivaljevi\'c~\cite{BjornerLovaszVrecicaZivaljevic}, see also~\cite[Thm.\,1.1]{ShareshianWachs2007}: the chessboard complex $\Delta_{\ell,b}$, which is $\LL_{K_\ell}(L)$ for equal lists $L_1=\dots=L_\ell$ of size $b$, is $(\nu(\ell,b)-2)$-connected.
To reach lists which are not equal we split the list of one vertex.

\medskip
\begin{lemma}
	\label{lem : splitting}
	Let $F$ be a graph, $u$ a vertex of $F$, and $L_u=S_1\sqcup\dots\sqcup S_k$ a partition of its list into non-empty sets.
	For $1\leq j\leq k$ let $L^{(j)}$ be the family of lists $L$ with $L_u$ replaced by $S_j$, and let $\tau\geq-1$.
	If all complexes $\LL_F(L^{(j)})$ are homologically $\tau$-connected and $\LL_{F-u}(L)$ is homologically $(\tau-1)$-connected, then $\LL_F(L)$ is homologically $\tau$-connected.
\end{lemma}
\begin{proof}
	The complex $\LL_F(L^{(j)})$ consists of the faces of $\LL_F(L)$ which do not colour $u$ or colour it by an element of $S_j$.
	These sub-complexes cover $\LL_F(L)$, and any $t\geq2$ of them intersect in the faces which do not colour $u$, that is in $\LL_{F-u}(L)$.
	For $\tau=-1$ the claim holds since $\LL_F(L^{(1)})\neq\emptyset$.
	For $\tau\geq0$ the complex $\LL_{F-u}(L)$ has a vertex, as it is homologically $(\tau-1)$-connected, and it is homologically $(\tau-t+1)$-connected for $t\geq2$.
	So Meshulam's nerve theorem~\cite[Thm.\,2.1]{Meshulam2001}, in the special case recalled in the proof of Lemma~\ref{lem : cones}, applies to this cover and gives the claim.
\end{proof}

\medskip
\begin{proposition}
	\label{prop : complete graphs}
	Let $\ell\geq0$ and $b\geq1$ be integers, and let $C$ be a finite set of colours with a fixed partition into parts of size $b$.
	Let $L=(L_1,\dots,L_\ell)$ be lists for the vertices $1,\dots,\ell$ of the complete graph $K_\ell$ such that every $L_i\subseteq C$ is a non-empty union of parts.
	Then the complex $\LL_{K_\ell}(L)$ is homologically $(\nu(\ell,b)-2)$-connected over every field.
\end{proposition}
\begin{proof}
	We first show that
	\begin{equation}
		\label{eq : nu subadditive}
		\nu(\ell_1,b)+\nu(\ell_2,b)\ \geq\ \nu(\ell_1+\ell_2,b)\qquad\text{for all }\ell_1,\ell_2\geq1.
	\end{equation}
	If $\nu(\ell_1,b)=b$ or $\nu(\ell_2,b)=b$, this holds since $\nu(\ell_1+\ell_2,b)\leq b$, and this covers $b=1$.
	Otherwise, $b\geq2$, and for $i=1,2$ the number $\nu(\ell_i,b)$ is $\ell_i$ or $\lfloor\frac{\ell_i+b+1}{3}\rfloor$.
	Since $\nu(\ell_1+\ell_2,b)$ is at most $\ell_1+\ell_2$ and at most $\lfloor\frac{\ell_1+\ell_2+b+1}{3}\rfloor$, it suffices to check the following inequalities
	\begin{gather*}
		\ell_1+\ell_2\ \geq\ \nu(\ell_1+\ell_2,b),\qquad\quad
		\ell_1+\Big\lfloor\frac{\ell_2+b+1}{3}\Big\rfloor\ \geq\ \Big\lfloor\frac{\ell_1+\ell_2+b+1}{3}\Big\rfloor,\\
		\Big\lfloor\frac{\ell_1+b+1}{3}\Big\rfloor+\Big\lfloor\frac{\ell_2+b+1}{3}\Big\rfloor\ \geq\ \Big\lfloor\frac{\ell_1+\ell_2+2b+2}{3}\Big\rfloor-1\ \geq\ \Big\lfloor\frac{\ell_1+\ell_2+b+1}{3}\Big\rfloor,
	\end{gather*}
	together with the inequality symmetric to the second one.
	The first two are clear, in the third we use $\lfloor x\rfloor+\lfloor y\rfloor\geq\lfloor x+y\rfloor-1$ and then $b\geq2$.

	\medskip
	We now prove the proposition by induction on the number of pairs $(i,P)$ with $1\leq i\leq\ell$ and $P$ a part contained in $L_i$.
	For $\ell=0$ there is nothing to prove.
	Assume first that every list is a single part.
	Different parts are disjoint, so two vertices with different lists can never get the same colour, and the edge between them imposes no condition.
	Hence $\LL_{K_\ell}(L)$ is the join, over the parts $P$ which occur as lists, of the complexes of the complete graphs on the vertices with the list $P$.
	If $\ell_P\geq1$ vertices have the list $P$, the corresponding factor is the chessboard complex $\Delta_{\ell_P,b}$.
	By~\cite{BjornerLovaszVrecicaZivaljevic} the factor $\Delta_{\ell_P,b}$ is $(\nu(\ell_P,b)-2)$-connected, hence homologically $(\nu(\ell_P,b)-2)$-connected, and by the K\"unneth formula for joins~\cite[Lem.\,2.1]{Milnor1956} the join is homologically $(\sum_P\nu(\ell_P,b)-2)$-connected.
	Here and below the sums run over the parts $P$ which occur as lists.
	Since $\sum_P\ell_P=\ell$, repeated use of~\eqref{eq : nu subadditive} gives
	\[
		\sum_P\nu(\ell_P,b)\ \geq\ \nu(\ell,b).
	\]

	\smallskip
	Otherwise the list of some vertex $u$ is the union of $k\geq2$ parts $P_1,\dots,P_k$, and we apply Lemma~\ref{lem : splitting} with $S_j:=P_j$ and $\tau:=\nu(\ell,b)-2$.
	The lists $L^{(j)}$ are again unions of parts, with $k-1$ fewer pairs, so $\LL_{K_\ell}(L^{(j)})$ is homologically $\tau$-connected by induction.
	Moreover, $K_\ell-u$ is a complete graph on $\ell-1$ vertices whose lists are unions of parts, with $k$ fewer pairs, so $\LL_{K_\ell-u}(L)$ is homologically $(\nu(\ell-1,b)-2)$-connected by induction, where
	\[
		\nu(\ell-1,b)\ \geq\ \nu(\ell,b)-1,
	\]
	since each of the three terms in the minimum decreases by at most one.
	So the assumptions of Lemma~\ref{lem : splitting} hold, and $\LL_{K_\ell}(L)$ is homologically $\tau$-connected.
\end{proof}

\medskip
We do not know whether Proposition~\ref{prop : complete graphs} holds for arbitrary lists of size at least $b$.
For $b\geq2\ell-2$ it does, by Theorem~\ref{th : colourings} with $\kappa\leq1$, since $\nu(\ell,2\ell-2)=\ell-1$.

\medskip
The numbers $\varepsilon_b$ of Section~\ref{sec : introduction} come from the connectivity bounds of Theorem~\ref{th : colourings}, Proposition~\ref{prop : stars} and Proposition~\ref{prop : complete graphs}, applied with lists of at least $b$ colours. They measure how far these bounds fall short of the connectivity $v-2$ of a complex of partial colourings of a graph with $v$ vertices.
For a vertex $w$ and a root $u$, the summand $\max\{0,\,2\deg_H(w)+\chi_u(w)-b\}$ is the number of colours which are missing at $w$ compared with $2\deg_H(w)+\chi_u(w)$, the list size which Theorem~\ref{th : colourings} requires at $w$. One extra colour is needed at the root of the induction and at the vertices of induced cycles of length at least four.
In the case of a complete graph $K_\ell$ the number $\ell-\nu(\ell,b)$ compares the connectivity $\nu(\ell,b)-2$ of the chessboard complex $\Delta_{\ell,b}$ with $\ell-2$, and it vanishes exactly when $b\geq2\ell-1$.
Consequently, $\varepsilon_b(H)=0$ whenever $b\geq2\Delta+1$, where $\Delta$ is the maximum degree of $H$.
The next proposition turns these observations into a connectivity bound.

\medskip
\begin{proposition}
	\label{prop : condition B}
	Let $r\geq2$ and $n\geq0$ be integers, let $H\in\mathscr B$ be a graph on the vertex set of $\Delta_n$, let $b\geq1$ be a divisor of $r$, and fix a partition of $[r]$ into parts of size $b$.
	If $U\subseteq V$ is a set such that every $U\cap Q_j$, viewed as a subset of $[r]$, is a non-empty union of parts, then $\LL_H[U]$ is homologically $(n-1-\varepsilon_b(H))$-connected over every field.
	When moreover every $U\cap Q_j$ has at least $m\geq b$ elements, the complex $\LL_H[U]$ is homologically $(n-1-\varepsilon_{b,m}(H))$-connected over every field.
\end{proposition}
\begin{proof}
	The complex $\LL_H[U]$ is the complex of partial proper list colourings of $H$ with lists of size at least $b$ which are unions of parts, and it is the join of the corresponding complexes of the connected components of $H$, since there are no restrictions between vertices of different components.
	Let $C$ be a component with $v$ vertices.
	If $C$ is a complete graph, Proposition~\ref{prop : complete graphs} shows that its complex is homologically $(v-2-\varepsilon_b(C))$-connected.
	For a star $K_{1,k}$ with $k\leq b-2$ this follows from Proposition~\ref{prop : stars}.
	Otherwise we choose as root a vertex $u$ of $C$ with $\varepsilon_b(C,u)=\varepsilon_b(C)$, and Theorem~\ref{th : colourings} gives the bound, since the sum $\kappa$ of $C$ is at most $\varepsilon_b(C,u)$ when all lists have at least $b$ elements.
	The components have $n+1$ vertices altogether, and the K\"unneth formula for joins~\cite[Lem.\,2.1]{Milnor1956} shows that the join is homologically $(n+1-\varepsilon_b(H)-2)$-connected.
	For the second statement the same argument applies with $\varepsilon_{b,m}(C)$ in place of $\varepsilon_b(C)$: for a complete graph both Proposition~\ref{prop : complete graphs} and Theorem~\ref{th : colourings} apply, and we take the better bound, for a star $K_{1,k}$ with $k\leq m-2$ the lists have at least $k+2$ elements, and otherwise the sum $\kappa$ of Theorem~\ref{th : colourings} is at most $\varepsilon_m(C,u)$.
\end{proof}

\medskip
We close this section with a connectivity bound for arbitrary graphs which depends only on the local structure of the graph and on the sizes of the lists.
It is a version for lists of Engstr\"om's bound~\cite[Thm.\,3.5]{Engstrom2011}, and the proof follows his.
We need the following lemma, which is~\cite[Lem.\,3.2]{Engstrom2011}, based on~\cite[Lem.\,10.4(ii)]{Bjorner1995}. Its short proof is included here in our homological setting.

\medskip
\begin{lemma}
	\label{lem : clique neighbourhood}
	Let $\tau\geq-1$, let $\Lambda$ be a graph, and let $x$ be a vertex of $\Lambda$ whose neighbours are pairwise adjacent.
	If $\Ind(\Lambda-(N[x]\cup N[u]))$ is homologically $(\tau-1)$-connected for every neighbour $u$ of $x$, then $\Ind(\Lambda)$ is homologically $\tau$-connected.
\end{lemma}
\begin{proof}
	For a vertex $z$ let $\St(z):=\{\sigma\in\Ind(\Lambda) : \sigma\cup\{z\}\in\Ind(\Lambda)\}$ be its closed star, a cone with apex $z$.
	Every maximal face of $\Ind(\Lambda)$ contains a vertex of $N[x]$, as otherwise $x$ could be added, and so $\Ind(\Lambda)$ is the union of the stars $\St(z)$, $z\in N[x]$.
	As $N[x]$ is a clique, $N[x]\subseteq N[u]$ for every neighbour $u$ of $x$.
	Hence, for different neighbours $u$ and $u'$ of $x$, a face of $\St(u)\cap\St(u')$ avoids $N[u]\cup N[u']\supseteq N[x]$, and so it lies in $\St(x)$.
	Let $u_1,\dots,u_s$ be the neighbours of $x$, and let $K_t:=\St(x)\cup\St(u_1)\cup\cdots\cup\St(u_t)$.
	Then
	\[
		K_{t-1}\cap\St(u_t)\ =\ \St(x)\cap\St(u_t)\ =\ \Ind\big(\Lambda-(N[x]\cup N[u_t])\big),
	\]
	which is homologically $(\tau-1)$-connected.
	The cone $K_0=\St(x)$ has vanishing reduced homology, and the Mayer--Vietoris sequence
	\[
	\begin{tikzcd}[column sep=2.2em]
		\widetilde H_i(K_{t-1})\oplus\widetilde H_i(\St(u_t)) \arrow[r] & \widetilde H_i(K_t) \arrow[r] & \widetilde H_{i-1}\big(K_{t-1}\cap\St(u_t)\big)
	\end{tikzcd}
	\]
	shows by induction on $t$ that $K_t$ is homologically $\tau$-connected.
	For $t=s$ we obtain $\Ind(\Lambda)$.
\end{proof}

\begin{figure}[ht]
	\centering
	\begin{tikzpicture}[x=9mm,y=5.5mm,
			pt/.style={circle,fill=black,inner sep=1.5pt},
			hollow/.style={circle,draw,fill=white,inner sep=1.5pt},
			gv/.style={circle,fill=black,inner sep=1.7pt},
			lab/.style={font=\scriptsize}]
		\foreach \X/\T in {0/(A),7.2/(B)}{
			\begin{scope}[xshift=\X*9mm]
				\draw[dashed,rounded corners=2pt] (1.6,-0.5) rectangle (2.4,-5.5);
				\foreach \i in {1,...,5}{\node[lab,anchor=east] at (-0.5,-\i) {$\i$};}
				\foreach \j in {0,...,4}{\foreach \i in {1,...,5}{\node[hollow] at (\j,-\i) {};}}
				\foreach \i in {1,...,5}{\node[pt] at (2,-\i) {};}
				\foreach \j/\n in {0/x_1,1/v',2/w,3/x_2,4/x_3}{\node[gv] (v\j) at (\j,-7) {}; \node[lab,below] at (\j,-7.2) {$\n$};}
				\draw[thick] (v0)--(v1)--(v2)--(v3)--(v4);
				\node[font=\small] at (2,-9.3) {\T};
			\end{scope}}
		\foreach \j in {0,1,3}{\node[pt] at (\j,-2) {};}
		\begin{scope}[xshift=7.2*9mm]
			\foreach \j in {1,3}{\foreach \i in {2,4}{\node[pt] at (\j,-\i) {};}}
		\end{scope}
		\node[lab,anchor=west] at (4.4,-7) {$F$};
		\node[lab,anchor=west] at (11.6,-7) {$F$};
	\end{tikzpicture}
	\caption{\small Squids with the body $w$ for the path $F$ with the vertices $x_1,v',w,x_2,x_3$ and all lists equal to $[5]$.
		The vertex $(x,c)$ of $F\cart L$ is drawn in the column of $x$ and the row of the colour $c$. The fibre over $w$ is framed, and filled vertices belong to the squid.
		(A) The largest squid of type (A) for the neighbour $v'$ of $w$ and the colour $i=2$: the fibre over $w$ together with the vertices $(x,2)$ for $x\in N_F(v')\cup N_F(w)=\{x_1,w,v',x_2\}$.
		(B) The largest squid of type (B) for the colours $i=2$ and $i'=4$: the fibre over $w$ together with the vertices $(x,2)$ and $(x,4)$ for $x\in N_F(w)=\{v',x_2\}$.}
	\label{fig : squids}
\end{figure}
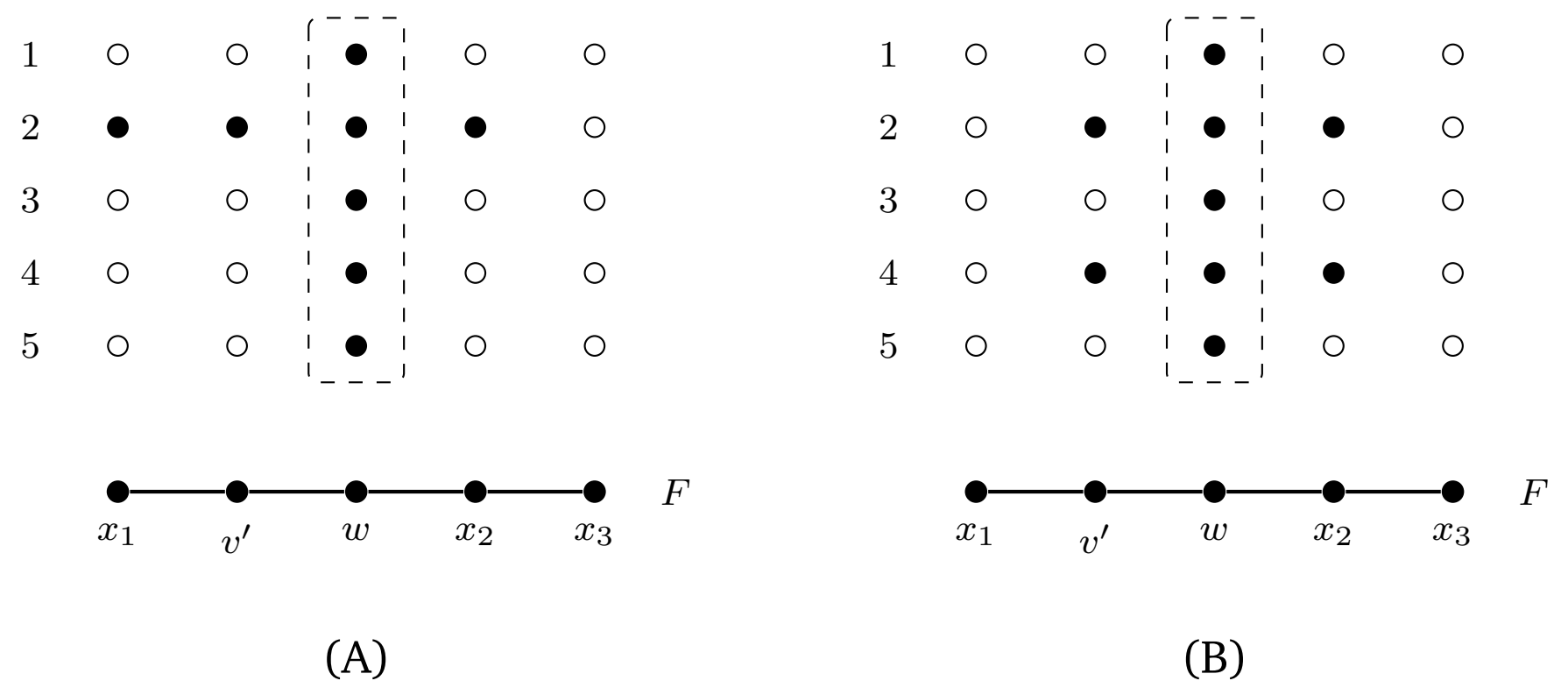

\medskip
For a vertex $v$ of a graph $F$ let $N_F(v)$ be the set of its neighbours and $N^2_F(v)$ the set of vertices at distance two from $v$.
For $X\subseteq Y$ and a colour $i$ we write $X\times\{i\}$ for the set of vertices $(x,i)$ of $F\cart  L$ with $x\in X$.
The following sets of vertices are the squids of~\cite[Def.\,3.3]{Engstrom2011}, adapted to lists.
A {\em squid with body} $w\in Y$ is a set $S$ of vertices of $F\cart  L$ which contains the fibre $\{w\}\times L_w$ and satisfies one of the following:
\begin{compactenum}[\rm\quad (A)]
	\item $S\subseteq(N_F(v')\cup N_F(w))\times\{i\}\cup\{w\}\times L_w$ for a neighbour $v'$ of $w$ and a colour $i$, or
	\item $S\subseteq N_F(w)\times\{i,i'\}\cup\{w\}\times L_w$ for two colours $i\neq i'$.
\end{compactenum}
Removing a squid from $F\cart  L$ amounts to removing its body from $F$ and one or two colours from the lists of some vertices near the body.
Figure~\ref{fig : squids} shows the largest squids of both types for a path.

\medskip
\begin{lemma}
	\label{lem : squid count}
	Let $F$ be a graph on $Y$ and $L$ lists with $|L_v|>|N^2_F(v)|+2|N_F(v)|$ for every $v\in Y$.
	Let $S_1,\dots,S_m$ be squids with pairwise different bodies.
	Then for every vertex $u\in Y$ which is not one of the bodies, the fibre over $u$ is not contained in $S_1\cup\cdots\cup S_m$.
\end{lemma}
\begin{proof}
	Let $S$ be a squid with body $w\neq u$ which meets the fibre over $u$.
	If $S$ is of type (A), then $u\in N_F(w)$, or $u\in N_F(v')$ for a neighbour $v'$ of $w$, so $w\in N_F(u)\cup N^2_F(u)$, and $S$ contains at most one vertex of the fibre over $u$.
	If $S$ is of type (B), then $w\in N_F(u)$, and $S$ contains at most two vertices of the fibre over $u$.
	Since the bodies are different, the squids remove at most $2|N_F(u)|+|N^2_F(u)|<|L_u|$ vertices of the fibre over $u$.
\end{proof}

\medskip
\begin{theorem}
	\label{th : squids}
	Let $F$ be a graph on $Y$ and $L$ lists with $|L_v|>|N^2_F(v)|+2|N_F(v)|$ for every $v\in Y$.
	If $S_1,\dots,S_m$ are squids with pairwise different bodies and $m<|Y|$, then the complex $\Ind\big((F\cart  L)-(S_1\cup\cdots\cup S_m)\big)$ is homologically $(|Y|-m-2)$-connected.
	In particular, $\LL_F(L)$ is homologically $(|Y|-2)$-connected.
\end{theorem}
\begin{proof}
	We use induction on $|Y|-m$.
	Let $\Lambda:=(F\cart  L)-(S_1\cup\cdots\cup S_m)$ and $\tau:=|Y|-m-2$.
	If $|Y|-m=1$, then $\tau=-1$, and $\Lambda$ has a vertex by Lemma~\ref{lem : squid count}.
	Let $|Y|-m\geq2$.
	Choose a vertex $u\in Y$ which is not a body, and by Lemma~\ref{lem : squid count} a vertex $(u,i)$ of $\Lambda$.
	Let $w_1,\dots,w_t$ be the neighbours $w$ of $u$ in $F$ with $(w,i)\in\Lambda$, and for $0\leq l\leq t$ let $\Lambda_l:=\Lambda-\{(w_1,i),\dots,(w_l,i)\}$.

	\smallskip
	First, $\Ind(\Lambda_t)$ is homologically $\tau$-connected.
	In $\Lambda_t$ the neighbours of $(u,i)$ are the other vertices of the fibre over $u$, and they are pairwise adjacent.
	For such a neighbour $(u,i')$ we have
	\[
		\Lambda_t-\big(N[(u,i)]\cup N[(u,i')]\big)\ =\ (F\cart  L)-(S_1\cup\cdots\cup S_m\cup S'),\qquad S':=N_F(u)\times\{i,i'\}\cup\{u\}\times L_u,
	\]
	where $S'$ is a squid of type (B) with the body $u$, which is different from the bodies of $S_1,\dots,S_m$.
	Since $m+1<|Y|$, by induction the independence complex of this graph is homologically $(\tau-1)$-connected, and Lemma~\ref{lem : clique neighbourhood} applies.

	\smallskip
	Next, for $l=t,t-1,\dots,1$ consider the vertex $(w_l,i)$ of $\Lambda_{l-1}$.
	Its deletion is $\Lambda_{l-1}-(w_l,i)=\Lambda_l$, and its link corresponds to
	\[
		\Lambda_{l-1}-N[(w_l,i)]\ =\ (F\cart  L)-(S_1\cup\cdots\cup S_m\cup S''_l),
	\]
	where
	\[
		S''_l:=\big(\{w_1,\dots,w_{l-1}\}\cup N_F(w_l)\big)\times\{i\}\cup\{w_l\}\times L_{w_l}.
	\]
	Since $w_1,\dots,w_{l-1}$ are neighbours of $u$ and $u$ is a neighbour of $w_l$, the set $S''_l$ is a squid of type (A) with $v'=u$ and the body $w_l$, which is not a body of $S_1,\dots,S_m$, because $(w_l,i)\in\Lambda$.
	Again by induction, the link is homologically $(\tau-1)$-connected, and (D1) shows that $\Ind(\Lambda_{l-1})$ is homologically $\tau$-connected if $\Ind(\Lambda_l)$ is.
	For $l=1$ we obtain $\Ind(\Lambda_0)=\Ind(\Lambda)$.
\end{proof}

\bigskip
\section{Proofs of the main results}
\label{sec : proofs}
\bigskip

\begin{figure}[ht]
	\centering
	\begin{tikzpicture}[node distance=4mm,
			box/.style={draw,rounded corners=2pt,align=center,font=\small,inner sep=4pt,minimum height=14mm},
			arr/.style={-{Stealth[length=2mm]}}]
		\node[box,text width=24mm] (g) {\small group $G\leq\sym_r$\\[2pt] {\footnotesize Lemma~\ref{lem : group}}};
		\node[box,text width=52mm,right=of g] (c) {\small connectivity of $\LL_H$ and $\LL_H[U]$\\[2pt] {\footnotesize Prop.~\ref{prop : condition B}, based on Thm.~\ref{th : colourings} and}\\ {\footnotesize Props.~\ref{prop : stars} and~\ref{prop : complete graphs}; or Thm.~\ref{th : squids}}};
		\node[box,text width=36mm,right=of c] (t) {\small Tverberg face $J\in\LL_H$\\[2pt] {\footnotesize Thm.~\ref{th : G-box} or Thm.~\ref{th : refined}}};
		\node[box,text width=22mm,right=of t] (f) {\small greedy extension of $J$ to a facet};
		\draw[arr] (g)--(c); \draw[arr] (c)--(t); \draw[arr] (t)--(f);
	\end{tikzpicture}
	\caption{\small The pattern of the proofs in this section.
		In the last step the vertices of $\Delta_n$ which are not labelled by $J$ are labelled one after another, each with a label which none of its already labelled neighbours in $H$ has. This is possible since every vertex has at most $r-1$ neighbours in $H$.
		The resulting facet of $\K$ lies in $\LL_H$ and contains $J$, so it is a Tverberg partition which covers all vertices of $\Delta_n$ and separates the endpoints of every edge of $H$.}
	\label{fig : pattern}
\end{figure}
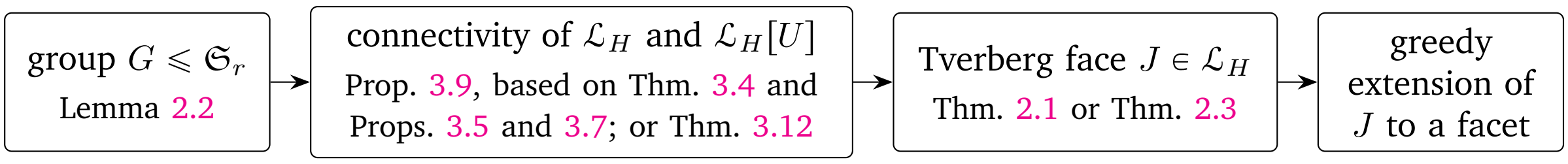

In this section we prove Theorem~\ref{th : main}, Theorem~\ref{th : small q} and Theorem~\ref{th : local} and Corollary~\ref{cor : no extra vertices}, Corollary~\ref{cor : coloured} and Corollary~\ref{cor : prime powers}.
All proofs follow the same pattern, which is shown in Figure~\ref{fig : pattern}.
First, we choose an elementary abelian group $G\leq\sym_r$ with Lemma~\ref{lem : group}.
The colouring complex $\LL_H\subseteq\K$ is invariant under every permutation of the labels, in particular under $G$, and for a $G$-box $U$ every column $U\cap Q_j$ is a non-empty union of $G$-orbits.
So $\LL_H[U]$ is a complex of partial proper list colourings of $H$ whose lists are unions of orbits, and $\LL_H$ is the complex with all lists equal to $[r]$.
Second, a connectivity bound of Section~\ref{sec : list colourings} shows that these complexes are sufficiently connected.
Third, the $G$-box criterion or its refinement gives a Tverberg face $J\in\LL_H$ of the given map.
Finally, $J$ is extended greedily to a facet of $\K$ in $\LL_H$, which is a Tverberg partition covering all vertices of $\Delta_n$ and separating the endpoints of every edge of $H$.

\medskip
The proofs differ in the choice of the group, the criterion, and the connectivity bound, as listed in Table~\ref{tab : pattern}.
Only the proof of Theorem~\ref{th : small q} uses the refined criterion: it needs the closed strata of Theorem~\ref{th : closed strata}, and it treats maps whose images do not affinely span $\R^d$ by a limit argument.
For prime powers the group has one orbit, the only $G$-box is $V$, and the proofs work for continuous maps.

\begin{table}[b]
	\centering\small
	\renewcommand{\arraystretch}{1.25}
	\setlength{\tabcolsep}{3.5pt}
	\begin{tabular}{lllll}
		\toprule
		result & orbits of $G$ & criterion & connectivity bound & maps\\
		\midrule
		Thm.~\ref{th : main} & size $b(r)$ & Thm.~\ref{th : G-box} & Prop.~\ref{prop : condition B} with $b=b(r)$ and $b=r$ & affine\\
		Thm.~\ref{th : small q} & $\frac rq$ orbits of size $q$ & Thm.~\ref{th : refined} & Prop.~\ref{prop : condition B} with $b=q$, $m=r-dq$, and $b=r$ & affine\\
		Cor.~\ref{cor : no extra vertices} & size $b(r)$ & Thm.~\ref{th : G-box} & via Thm.~\ref{th : main}, $\varepsilon_r(H)=\varepsilon_{b(r)}(H)=0$ & affine\\
		Cor.~\ref{cor : coloured} & size $b(r)$ & Thm.~\ref{th : G-box} & via the proof of Thm.~\ref{th : main} & affine\\
		Cor.~\ref{cor : prime powers} & one orbit, $r=q$ & Thm.~\ref{th : G-box} & Prop.~\ref{prop : condition B} with $b=q$ & continuous\\
		Thm.~\ref{th : local} & size $b(r)$ & Thm.~\ref{th : G-box} & Thm.~\ref{th : squids} & affine\\
		Thm.~\ref{th : local}, $r=q$ & one orbit & Thm.~\ref{th : G-box} & Thm.~\ref{th : squids} & continuous\\
		\bottomrule
	\end{tabular}\vspace{5pt}
	\caption{\small The ingredients of the proofs in this section.}
	\label{tab : pattern}
\end{table}

\medskip
\subsection{Proof of Theorem~\ref{th : main}}

Let $a\colon\Delta_n\longrightarrow\R^d$ be an affine map.
By Lemma~\ref{lem : group} there is an elementary abelian group $G\leq\sym_r$ all of whose orbits have $b(r)$ elements, and these orbits form a partition of $[r]$ into parts of size $b(r)$.
The colouring complex $\LL_H\subseteq\K$ is invariant under every permutation of the labels, in particular under $G$.
For a $G$-box $U$ every $U\cap Q_j$ is a non-empty union of orbits.
By Proposition~\ref{prop : condition B}, $\LL_H[U]$ is homologically $(n-1-\varepsilon_{b(r)}(H))$-connected, and $n-1-\varepsilon_{b(r)}(H)\geq N-2$ by assumption.
Moreover, $\LL_H=\LL_H[V]$, and every $V\cap Q_j=[r]$ is the single part of the partition of $[r]$ into one part of size $r$.
Hence Proposition~\ref{prop : condition B} with $b=r$ shows that $\LL_H$ is homologically $(n-1-\varepsilon_r(H))$-connected, and $n-1-\varepsilon_r(H)\geq N-1$ by assumption.
So $\LL_H$ satisfies the assumptions of Theorem~\ref{th : G-box}, which gives a Tverberg face $J\in\LL_H$ of $a$.

\medskip
We extend $J$ to a facet of $\K$ which lies in $\LL_H$, labelling the vertices of $\Delta_n$ which are not labelled by $J$ one after another.
At each step the vertex has at most $r-1$ neighbours in $H$, and so there is a label which none of its labelled neighbours has.
The resulting facet lies in $\LL_H$, and it is a Tverberg face of $a$, since it contains the Tverberg face $J$, see Section~\ref{sec : deleted join}.
Its faces contain all vertices of $\Delta_n$, and the two vertices of every edge of $H$ lie in different faces.
\qed

\medskip
\subsection{Proof of Theorem~\ref{th : small q}}

Let $a\colon\Delta_n\longrightarrow\R^d$ be an affine map such that $a(p_0),\dots,a(p_n)$ affinely span $\R^d$.
Put $k:=\frac{r}{q}\geq d+2\geq3$, and let $G$ be the group of Lemma~\ref{lem : group} with $k$ orbits of size $q$.
We check the assumptions (A)--(C) of Theorem~\ref{th : refined} for $\LL_H$, which is invariant under every permutation of the labels.
By Proposition~\ref{prop : condition B} with $b=r$, the complex $\LL_H$ is homologically $(n-1-\varepsilon_r(H))$-connected, and $n-1-\varepsilon_r(H)\geq N-1$.
For a $G$-box $U$, Proposition~\ref{prop : condition B} with $b=q$ shows that $\LL_H[U]$ is homologically $(n-1-\varepsilon_q(H))$-connected, and $n-1-\varepsilon_q(H)\geq N-3\geq0$, since $N\geq(d+1)(3q-1)\geq10$. In particular, $\LL_H[U]$ is connected.
If the stratum of $U$ is non-empty and closed, then by Theorem~\ref{th : closed strata} every $U\cap Q_j$ is a union of at least $k-d$ orbits, so it has at least $(k-d)q=r-dq$ elements, and by Proposition~\ref{prop : condition B} with $b=q$ and $m=r-dq$ the complex $\LL_H[U]$ is homologically $(n-1-\varepsilon_{q,r-dq}(H))$-connected, where $n-1-\varepsilon_{q,r-dq}(H)\geq N-2$.
Theorem~\ref{th : refined} gives a Tverberg face of $a$ in $\LL_H$.

\medskip
An arbitrary affine map $a$ is the limit of affine maps $a_t$ with spanning images, because $n\geq N\geq d$.
Since $\LL_H$ has finitely many faces, one of them is a Tverberg face of infinitely many $a_t$, and since the intersection condition passes to the limit, it is a Tverberg face of $a$.
It extends to a facet of $\K$ in $\LL_H$, which is a Tverberg face of $a$, exactly as in the proof of Theorem~\ref{th : main}.
\qed

\medskip
The entries of Table~\ref{tab : complete} are obtained by evaluating the conditions of Theorem~\ref{th : main} and Theorem~\ref{th : small q} for $H=K_\ell$ and all prime powers $q$ dividing $r$.
Here $\varepsilon_m(K_\ell,u)=(\ell-1)\max\{0,2\ell-2-m\}+\max\{0,2\ell-1-m\}$ for every vertex $u$.

\medskip
\subsection{Proof of Corollary~\ref{cor : no extra vertices}}

Let $b=b(r)\geq3$.
Every component $C$ of $H$ belongs to $\mathscr B$, and we show $\varepsilon_b(C)=0$.
If $C$ is a complete graph $K_\ell$ as in (A), (B) or (C), then $\ell\leq2$ or $\ell-1\leq\frac{b-1}{2}$, so $b\geq2\ell-1$, that is $\nu(\ell,b)=\ell$.
A star $K_{1,k}$ as in (A), (B) or (C) which is not complete has $2\leq k\leq b-2$, and so $\varepsilon_b(C)=0$ by definition.
For any other component $C$ as in (A) we have $2\deg_H(w)+\chi_u(w)\leq b$ for all $w$ and every $u$.
Every other tree $C$ as in (B) has a leaf $u$, and then $2\deg_H(u)+\chi_u(u)=3\leq b$ and $2\deg_H(w)+\chi_u(w)=2\deg_H(w)\leq b$ for $w\neq u$.
In all cases $\varepsilon_b(C)=0$.
Since $b\leq r$, each of the above inequalities remains true with $r$ in place of $b$, and $\nu(\ell,b)=\ell$ implies $\nu(\ell,r)=\ell$.
Hence $\varepsilon_r(C)=0$ as well, that is, $\varepsilon_r(H)=\varepsilon_b(H)=0$.
Moreover, the maximum degree of $H$ is at most $\max\{\frac{b}{2},\,b-2\}\leq r-1$.
Theorem~\ref{th : main} applies with $n\geq N$.
\qed

\medskip
\subsection{Proof of Corollary~\ref{cor : prime powers}}

Let $q=p^\alpha$ and let $G=\FF_p^{\,\alpha}$ act on $[q]=\FF_p^{\,\alpha}$ by translation.
Then $G$ has one orbit, and the only $G$-box is $V$, for which every $V\cap Q_j$ is the single part $[q]$.
By Proposition~\ref{prop : condition B}, $\LL_H=\LL_H[V]$ is homologically $(n-1-\varepsilon_q(H))$-connected, hence homologically $(N-1)$-connected, so the assumptions of Theorem~\ref{th : G-box} hold, and it gives a Tverberg face $J\in\LL_H$ of every continuous map $f\colon\Delta_n\longrightarrow\R^d$.
The extension of $J$ to a facet of $\K$ in $\LL_H$ is the same as in the proof of Theorem~\ref{th : main}, and a face of $\K$ which contains a Tverberg face of a continuous map is again a Tverberg face, see Section~\ref{sec : deleted join}.
For the graphs of Corollary~\ref{cor : no extra vertices} with $b(r)$ replaced by $q\geq3$ the proof of that corollary shows $\varepsilon_q(H)=0$ and that the maximum degree is at most $q-1$.
\qed

\medskip
\subsection{Proof of Corollary~\ref{cor : coloured}}

Let $H$ be the disjoint union of the complete graphs on the colour classes $C_1,\dots,C_m$.
Then $H\in\mathscr B$, and a face of $\Delta_n$ is rainbow if and only if it contains no edge of $H$.
By the definition of the deficit, for every integer $b\geq1$ we have
\[
	\varepsilon_b(H)\ =\ \sum_{i=1}^m\big(|C_i|-\nu(|C_i|,b)\big)\ =\ n+1-\sum_{i=1}^m\nu(|C_i|,b).
\]
Hence, the two assumptions of the corollary are the inequalities $n-1-\varepsilon_r(H)\geq N-1$ and $n-1-\varepsilon_{b(r)}(H)\geq N-2$.
The first paragraph of the proof of Theorem~\ref{th : main} uses only these two inequalities and not the bound on the maximum degree of $H$.
Therefore, it gives a Tverberg face $J\in\LL_H$ of every affine map $a\colon\Delta_n\longrightarrow\R^d$, that is, a Tverberg partition of $a$ into rainbow faces.

\medskip
If $|C_i|\leq r$ for all $i$, then the maximum degree of $H$ is at most $r-1$, and $J$ extends to a facet of $\K$ in $\LL_H$ exactly as in the proof of Theorem~\ref{th : main}.
When $r=q$ is a prime power, we have $b(r)=r$ and the first assumption reads $n-1-\varepsilon_q(H)\geq N-1$.
In this case the proof of Corollary~\ref{cor : prime powers} gives a Tverberg face $J\in\LL_H$ of every continuous map $\Delta_n\longrightarrow\R^d$, and the extension to a facet is the same.
\qed

\medskip
\subsection{Proof of Theorem~\ref{th : local}}

Let $a\colon\Delta_n\longrightarrow\R^d$ be an affine map, and let $G$ be the group of Lemma~\ref{lem : group} whose orbits have $b(r)$ elements.
For a $G$-box $U$ every $U\cap Q_j$ is a non-empty union of orbits, so the complex $\LL_H[U]=\LL_H(U_0,\dots,U_n)$ has lists with at least $b(r)$ elements.
By the assumption and Theorem~\ref{th : squids}, applied to $H$, which has $n+1$ vertices, the complex $\LL_H[U]$ is homologically $(n-1)$-connected, and $n-1\geq N-1$.
The same holds for $\LL_H=\LL_H[V]$, whose lists $[r]$ have $r\geq b(r)$ elements.
So Theorem~\ref{th : G-box} gives a Tverberg face $J\in\LL_H$ of $a$.
Every vertex $v$ has $|N(v)|<\frac{b(r)}{2}\leq r-1$ neighbours, and so $J$ extends to a facet of $\K$ in $\LL_H$, which is a Tverberg face of $a$, exactly as in the proof of Theorem~\ref{th : main}.

\medskip
If $r=q$ is a prime power, we use the group $G$ with one orbit as in the proof of Corollary~\ref{cor : prime powers}.
The only $G$-box is $V$, and $\LL_H$ is homologically $(n-1)$-connected by Theorem~\ref{th : squids}, so Theorem~\ref{th : G-box} gives a Tverberg face of every continuous map, which extends as before.
\qed

\bibliographystyle{amsplain}

\end{document}